\documentclass[a4paper,12pt]{article}

\newcommand{\red}[1]{{\colorbox{red}{{#1}}}}

\usepackage{color,makeidx,bm,latexsym,amscd,amsmath,amssymb,enumerate,stmaryrd,amsthm,float,graphicx,psfrag,geometry,comment}

\usepackage[all]{xy}
\usepackage{bm}

\usepackage{tikz,epic,eso-pic}

\tikzset{v/.style={
  circle, draw, inner sep=2pt, minimum size=6pt, fill=white}}

\theoremstyle{plain}
\newtheorem{theorem}{Theorem}[section]

\newtheorem{lemma}[theorem]{Lemma}

\newtheorem{proposition}[theorem]{Proposition}

\theoremstyle{definition}
\newtheorem{definition}[theorem]{Definition}
\newtheorem{remark}[theorem]{Remark}
\newtheorem{example}[theorem]{Example}

\def\qed{\hfill $\Box$}

\newcommand{\Div}{\operatorname{Div}}
\newcommand{\order}{\operatorname{order}}
\newcommand{\Img}{\operatorname{Im}}
\newcommand{\Ker}{\operatorname{Ker}}
\newcommand{\ch}{\operatorname{ch}}
\newcommand{\bch}{\operatorname{bch}}
\newcommand{\uch}{\operatorname{uch}}
\newcommand{\codim}{\operatorname{codim}}
\newcommand{\Rect}{\operatorname{Rect}}
\newcommand{\Int}{\operatorname{Int}}
\newcommand{\wtilde}{\widetilde}
\newcommand{\gl}{\operatorname{gl}}
\newcommand{\OS}{\operatorname{OS}}
\newcommand{\Sing}{\operatorname{Sing}}
\newcommand{\Sep}{\operatorname{Sep}}
\newcommand{\Hom}{\operatorname{Hom}}

\newcommand{\rA}{\mathrm{A}}
\newcommand{\rB}{\mathrm{B}}
\newcommand{\rC}{\mathrm{C}}
\newcommand{\rD}{\mathrm{D}}
\newcommand{\rE}{\mathrm{E}}
\newcommand{\rF}{\mathrm{F}}
\newcommand{\rG}{\mathrm{G}}
\newcommand{\rH}{\mathrm{H}}
\newcommand{\rI}{\mathrm{I}}
\newcommand{\rJ}{\mathrm{J}}
\newcommand{\rK}{\mathrm{K}}
\newcommand{\rL}{\mathrm{L}}
\newcommand{\rM}{\mathrm{M}}
\newcommand{\rN}{\mathrm{N}}
\newcommand{\rO}{\mathrm{O}}
\newcommand{\rP}{\mathrm{P}}
\newcommand{\rQ}{\mathrm{Q}}
\newcommand{\rR}{\mathrm{R}}
\newcommand{\rS}{\mathrm{S}}
\newcommand{\rT}{\mathrm{T}}
\newcommand{\rU}{\mathrm{U}}
\newcommand{\rV}{\mathrm{V}}
\newcommand{\rW}{\mathrm{W}}
\newcommand{\rX}{\mathrm{X}}
\newcommand{\rY}{\mathrm{Y}}
\newcommand{\rZ}{\mathrm{Z}}

\newcommand{\frA}{\mathfrak{A}}
\newcommand{\frB}{\mathfrak{B}}
\newcommand{\frC}{\mathfrak{C}}
\newcommand{\frD}{\mathfrak{D}}
\newcommand{\frE}{\mathfrak{E}}
\newcommand{\frF}{\mathfrak{F}}
\newcommand{\frG}{\mathfrak{G}}
\newcommand{\frH}{\mathfrak{H}}
\newcommand{\frI}{\mathfrak{I}}
\newcommand{\frJ}{\mathfrak{J}}
\newcommand{\frK}{\mathfrak{K}}
\newcommand{\frL}{\mathfrak{L}}
\newcommand{\frM}{\mathfrak{M}}
\newcommand{\frN}{\mathfrak{N}}
\newcommand{\frO}{\mathfrak{O}}
\newcommand{\frP}{\mathfrak{P}}
\newcommand{\frQ}{\mathfrak{Q}}
\newcommand{\frR}{\mathfrak{R}}
\newcommand{\frS}{\mathfrak{S}}
\newcommand{\frT}{\mathfrak{T}}
\newcommand{\frU}{\mathfrak{U}}
\newcommand{\frV}{\mathfrak{V}}
\newcommand{\frW}{\mathfrak{W}}
\newcommand{\frX}{\mathfrak{X}}
\newcommand{\frY}{\mathfrak{Y}}
\newcommand{\frZ}{\mathfrak{Z}}

\newcommand{\bA}{\mathbb{A}}
\newcommand{\bC}{\mathbb{C}}
\newcommand{\bD}{\mathbb{D}}
\newcommand{\bF}{\mathbb{F}}
\newcommand{\bG}{\mathbb{G}}
\newcommand{\bH}{\mathbb{H}}
\newcommand{\bK}{\mathbb{K}}
\newcommand{\bP}{\mathbb{P}}
\newcommand{\bQ}{\mathbb{Q}}
\newcommand{\bR}{\mathbb{R}}
\newcommand{\bS}{\mathbb{S}}
\newcommand{\bX}{\mathbb{X}}
\newcommand{\bY}{\mathbb{Y}}
\newcommand{\bZ}{\mathbb{Z}}

\newcommand{\cA}{\mathcal{A}}
\newcommand{\cB}{\mathcal{B}}
\newcommand{\cC}{\mathcal{C}}
\newcommand{\cD}{\mathcal{D}}
\newcommand{\cE}{\mathcal{E}}
\newcommand{\cF}{\mathcal{F}}
\newcommand{\cG}{\mathcal{G}}
\newcommand{\cH}{\mathcal{H}}
\newcommand{\cI}{\mathcal{I}}
\newcommand{\cJ}{\mathcal{J}}
\newcommand{\cK}{\mathcal{K}}
\newcommand{\cL}{\mathcal{L}}
\newcommand{\cM}{\mathcal{M}}
\newcommand{\cN}{\mathcal{N}}
\newcommand{\cO}{\mathcal{O}}
\newcommand{\cP}{\mathcal{P}}
\newcommand{\cQ}{\mathcal{Q}}
\newcommand{\cR}{\mathcal{R}}
\newcommand{\cS}{\mathcal{S}}
\newcommand{\cT}{\mathcal{T}}
\newcommand{\cU}{\mathcal{U}}
\newcommand{\cV}{\mathcal{V}}
\newcommand{\cW}{\mathcal{W}}
\newcommand{\cX}{\mathcal{X}}
\newcommand{\cY}{\mathcal{Y}}
\newcommand{\cZ}{\mathcal{Z}}

\newcommand{\sA}{\mathsf{A}}
\newcommand{\sB}{\mathsf{B}}
\newcommand{\sC}{\mathsf{C}}
\newcommand{\sD}{\mathsf{D}}
\newcommand{\sE}{\mathsf{E}}
\newcommand{\sF}{\mathsf{F}}
\newcommand{\sG}{\mathsf{G}}
\newcommand{\sH}{\mathsf{H}}
\newcommand{\sI}{\mathsf{I}}
\newcommand{\sJ}{\mathsf{J}}
\newcommand{\sK}{\mathsf{K}}
\newcommand{\sL}{\mathsf{L}}
\newcommand{\sM}{\mathsf{M}}
\newcommand{\sN}{\mathsf{N}}
\newcommand{\sO}{\mathsf{O}}
\newcommand{\sP}{\mathsf{P}}
\newcommand{\sQ}{\mathsf{Q}}
\newcommand{\sR}{\mathsf{R}}
\newcommand{\sS}{\mathsf{S}}
\newcommand{\sT}{\mathsf{T}}
\newcommand{\sU}{\mathsf{U}}
\newcommand{\sV}{\mathsf{V}}
\newcommand{\sW}{\mathsf{W}}
\newcommand{\sX}{\mathsf{X}}
\newcommand{\sY}{\mathsf{Y}}
\newcommand{\sZ}{\mathsf{Z}}

\newcommand{\bgu}{\bigcup}
\newcommand{\bga}{\bigcap}

\def\qed{\hfill $\Box$}

\title{Handle decompositions and Kirby diagrams for the complements of plane algebraic curves}
\author{Sakumi Sugawara \thanks{Department of Mathematics, Faculty of Science, Hokkaido University, North 10, West 8, Kita-ku, Sapporo 060-0810, JAPAN, E-mail: sugawaras@math.sci.hokudai.ac.jp} }
\date{\today}

\begin{document}
\maketitle

\begin{abstract}
The complements of plane algebraic curves are well-studied from topological and algebro-geometric viewpoints. In this paper, we will describe an explicit handle decomposition and a Kirby diagram for the complement of a plane algebraic curve. The method is based on the notion of braid monodromy. We refined this technique to obtain handle decompositions and Kirby diagrams.
\end{abstract}

\renewcommand{\thefootnote}{\fnsymbol{footnote}}
%32S50: Topological aspects of complex singularities: Lefschetz theorems, topological classification, invariants
%32Q55: Topological aspects of complex manifolds
%57K40: General topology of 4-manifolds
%57R65: Surgeries and handlebodies
\renewcommand{\thefootnote}{\arabic{footnote}}

\section{Introduction}
Let $C$ be a plane algebraic curve in $\bC^2$. The complement $M = \bC^2 \setminus C$ has been studied from topological and algebro-geometric viewpoints. 
By the Zariski--van Kampen theorem, an algorithm for computing the fundamental group $\pi_1 (M)$ is classically known \cite{zar,van-Kampen}. 
The Zariski--van Kampen's method evolved into the notion of braid monodromy \cite{moi}, and Libgober gave a presentation of the fundamental group whose associated $2$-dimensional complex is homotopy equivalent to the complement \cite{lib-hom}.
Moreover, Artal Bartolo, Carmona Ruber, and Cogolludo Agust\'in proved that braid monodromy determines the embedded topology of projective plane curves \cite{acc}.

However, as far as the author knows, there are no algorithms to describe the explicit diffeomorphism type of $M$, which is more refined information. 
Since $M$ is a Stein surface, it admits a handle decomposition with handles of index at most $2$.
We can expect that the CW complex obtained by Libgober's method is the core of the handle decomposition.
However, his method uses the deformation retraction to the link complement in $S^3$ around each singularity of the curve.
Therefore, it seems difficult to extract more information than the homotopy type using his method.

The purpose of this paper is to give explicit handle decompositions of the complements of plane algebraic curves, refining the description of the homotopy type obtained from the braid monodromy.
We will describe a handle decomposition using Kirby diagrams
(for details on handle decompositions and Kirby diagrams of $4$-manifolds, see \cite{gom-sti} for example). 
Instead of constructing the deformation retraction, we decompose the neighborhood of each singularity into some pieces to describe the diffeomorphism type.
Then, we study how these pieces are attached.
At this stage, the key point is to deal with the critical points of a certain stable map from the boundary of the regular neighborhood of the curve (Proposition \ref{prop:stable-1}, \ref{prop:stable-2}).

Recently, handle decompositions and Kirby diagrams of the complements of complexified real line arrangements have been obtained by the author and Yoshinaga \cite{sug-yos}. 
This result corresponds to a special case of our main result.
However, their method is based on the Lefschetz hyperplane section theorem and is deeply dependent on the combinatorial structure of a real line arrangement.
Since it cannot be extended to general plane algebraic curves, we need to take another strategy.
%Their method depends on the real structure of arrangements, therefore cannot extend to general complex algebraic curves.

This paper is organized as follows. 
In Section \ref{sec:braid}, we recall the notion of braid monodromy for plane algebraic curves.
In Section \ref{sec:handle}, we give an explicit handle decomposition for the complement. 
In Section \ref{sec:kirby}, we introduce the algorithm to draw Kirby diagrams and give some examples.

\vspace{3mm}
\textbf{Acknowledgement.}
The author would like to thank Professor Naohiko Kasuya for the invaluable comments on this paper and helpful discussions on this research. 
He also thanks the referee for careful reading and comments to improve the paper.
The author is supported by the JSPS KAKENHI Grant Number 22KJ0114.

\section{Braid monodromy for plane algebraic curves}\label{sec:braid}
In this section, we recall the braid monodromy for plane algebraic curves (for details, see also \cite{acc, coh-suc-braid, moi, lib-hom}, for example).

\begin{figure}[htbp]
\centering
\begin{tikzpicture}
\coordinate (P) at (3.5,2.5);
\coordinate (Q) at (7,2.5);

\draw (0,0)--(9,0)--(12,2)--(3,2)--(0,0);
%底空間

\draw[->,->=stealth] (1.5,4) -++ (0,-2);
\draw (1.5,3) node[left]{$\pi$};

\draw (P) --++ (1.5,1) --++ (0,3.5) --++ (-1.5,-1) -- cycle;
\draw (P)++(0.3,4) node[above]{$L_{p_i}$};
%p_iのファイバー

\draw (Q) --++ (1.5,1) --++ (0,3.5) --++ (-1.5,-1) --cycle;
\draw (Q)++(0.3,4) node[above]{$L_p$};
%pのファイバー

\draw [line width=5pt, white] (4,5) to[out=0, in=210] (10,7);
\draw (4,5) to[out=0, in=210] (10,7);
%曲線上側
\draw [line width=5pt, white] (4,5) to[out=0, in=150] (10,3);
\draw (4,5) to[out=0, in=150] (10,3);
\draw (4,5) to[out=0, in=210] (10,7);
%曲線下側

\draw [line width=5pt, white](Q)++(0,0.2) --++(0,3.3);
\draw (Q) --++ (0,3.5);
%消えたQのファイバー

\fill [black] (4,5) circle (0.08) ;
\fill [black] (7.5,5.68) circle (0.08) ;
\fill [black] (7.5,4.3) circle (0.08) ;
%交点のぬりつぶし

\draw [densely dashed](P)++(0.5,0) --++(0,-1.5);
\fill [black] (P)++(0.5,0)++(0,-1.5) circle(0.06);
\draw (P)++(0.5,0)++(0,-1.5) node[left]{$p_i$};

\draw [densely dashed](Q)++(0.5,0) --++(0,-1.5);
\fill [black] (Q)++(0.5,0)++(0,-1.5) circle(0.06);
\draw (Q)++(0.5,0)++(0,-1.5) node[left]{$p$};

\draw (11,7) node[above] {$C= \{f(x,y) = 0\}$};

\draw (10,0) node[above] {$\bC$};

\end{tikzpicture}
\caption{The projection $\pi:\bC^2 \rightarrow \bC$}
\label{fig:projection_pi}
\end{figure}
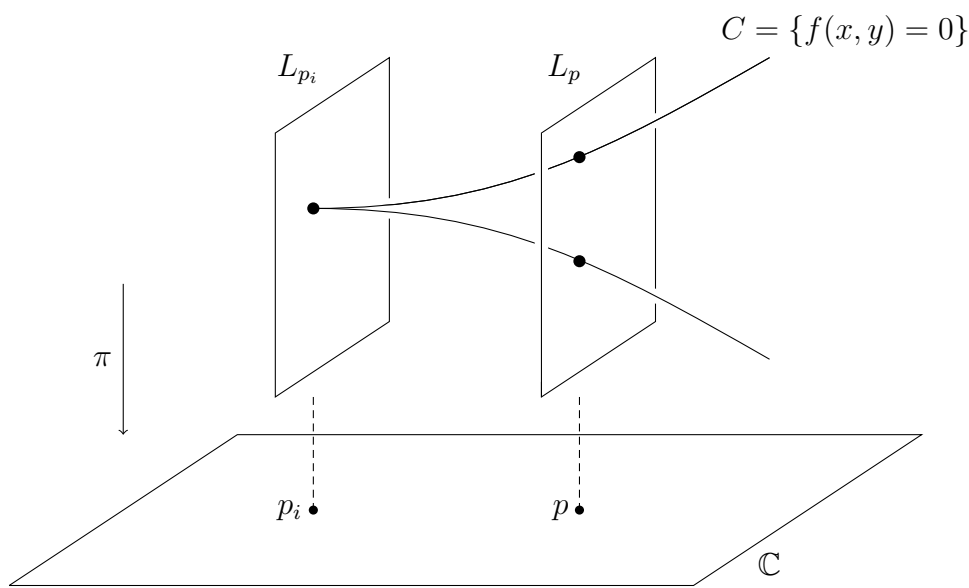

Let $C$ be a plane algebraic curve defined by a reduced polynomial $f(x,y)$ of degree $n$ such that the coefficient of $y^n$ is non-zero. 
We denote by $M = \bC^2 \setminus C$ its complement.
Let $\pi : \bC^2 \rightarrow \bC, \, \pi (x,y) = x$ be the first projection.
For $x \in \bC$, we set $L_{x} = \pi^{-1} (x)$. 
Let $X = \{p_1, \ldots , p_N\} \subset \bC$ be the set of points whose fiber $L_{p_{i}}$ contains a singularity of $C$ or $L_{p_i}$ is tangent to $C$ (see Figure \ref{fig:projection_pi}).
We can assume that the projection $\pi$ is generic, so that $L_{p_{i}}$ contains at most one singularity for each $p_i \in X$.

In this setting, we have the following proposition (see Lemma 4.3.5 in \cite{dim}).
\begin{proposition}
The restriction $\pi: \bC^{2} \setminus (C \cup \pi^{-1} (X))  \rightarrow \bC \setminus X $ is a locally trivial fibration.
\end{proposition}

\begin{definition}
Let $\Omega \subset \bC$ be a subset.
A system of arcs $s_1 , \ldots, s_{N}$ in $\Omega$ is called a \textit{non-intersecting arc system} on $\Omega$ with respect to the basepoint $p_{0} \in \Omega$ and a finite set $\{p'_1 , \ldots, p'_N\} \subset \Omega$ if the following two conditions are satisfied.
\begin{itemize}
\item Each arc $s_i$ is smooth and connects $p_0$ and $p'_{i}$.
\item There are no self-intersections in each $s_{i}$, and $s_{i} \cap s_{j} = \{p_{0}\}$ if $i \neq j$.
\end{itemize}
\end{definition}

For $1 \leq i \leq N$, let $U_i$ be a small open disk centered at $p_i$. 
Fix a point $p'_{i} \in \partial U_{i}$ for each $i$.
Take a point $p_{0} \in \bC \setminus X$ and a non-intersecting arc system $s_{1}, \ldots, s_{N}$ on $\Omega = \bC \setminus (U_{1} \cup \cdots \cup U_{N})$ with respect to the basepoint $p_{0}$ and $\{p'_1, \ldots, p'_{N}\}$.
We assume each $s_{i}$ intersects $\partial U_{i}$ transversally. 
Let $\gamma_{i} \in \pi_{1} (\bC \setminus X, p_{0})$ be the element which represents the loop based at $p_0$ obtained by concatenating $s_i$, $\partial U_i$ (going counter-clockwise) and $s_i^{-1}$.
%The fundamental group $\pi_{1} (\bC \setminus X, p_{0})$ is isomorphic to the free group $F_{N} = \langle  \gamma_{1}, \cdots, \gamma_{N} \rangle$ of rank $N$.

By trivializing the fibration $\pi: \bC^{2} \setminus (C \cup \pi^{-1} (X) ) \rightarrow \bC \setminus X$ along each loop $\gamma_{i}$, we obtain an automorphism of $L_{p_{0}}$, which preserves $L_{p_{0}} \cap C$.
The trivialization along the loops defines a homomorphism
\[
\Phi: \pi_{1} (\bC \setminus X, p_{0}) \rightarrow B(L_{p_{0}}, L_{p_{0}} \cap C),
\]
where $B(L_{p_{0}}, L_{p_{0}} \cap C)$ denotes the mapping class group of $L_{p_{0}}$ preserving the set $L_{p_{0}} \cap C$, which is known to be isomorphic to the Artin's braid group $B_{n}$ with $n$ strands.
This homomorphism $\Phi$ is called the braid monodromy of the algebraic curve.
We define $\beta_{i} = \Phi (\gamma_{i})$ as the monodromy along the loop $\gamma_{i}$.

\begin{figure}[htbp]
\centering
\begin{tikzpicture}
\coordinate (P) at (3.5,2.5);
\coordinate (Q) at (7,2.5);

\draw (0,0)--(9,0)--(12,2)--(3,2)--(0,0);
%底空間

\draw[->,->=stealth] (1.5,4) -++ (0,-2);
\draw (1.5,3) node[left]{$\pi$};
%射影$\pi$の矢印

\draw (P)++ (0.6,0) --++ (1.5,1) --++ (0,3.5) --++ (-1.5,-1) -- cycle;
\draw (P)++(0.3,4) node[above]{$L_{p'_i}$};
%p_iのファイバー

\draw (P)++(0.5,-1.5) circle (0.6cm and 0.4cm);
%$U_i$

\fill [black] (5.1,5.7) circle (0.08);
\draw (5.1,5.9) node[above] {$q_{i,1}$};
%\fill [black] (5.2,5) circle (0.08)node[left] {$q_{i,{m_{i}}}$} ;
%\fill [black] (5.3,4.5) circle (0.08) node[left] {$q_{i,{m_{i}+1}}$};
\fill [black] (5.1,3.6) circle (0.08); 
\draw (4.9,3.6) node[left] {$q_{i,n}$};
%交点のぬりつぶし

\fill[black] (4.4,4.7) circle (0.08) node[below] {$q_{i}$};
\draw (4.4,4.7) -- (5,5.42) ;
\draw (5.1,5.7) circle (0.3);
\draw (4.4,4.7) -- (5,3.87);
\draw (5.1,3.6) circle (0.3);
%基点とループ

%\draw (5.2,5.4) to[out=-90,in=-150] ++(0.8,0) node[right]{$e_1$};
%\draw (5.2,3.9) to[out=90,in=160] ++(0.8,0) node[right]{$e_n$};
%e_1,e_n名前

\fill[black] (5.3,5) circle(0.04);
\fill[black] (5.4,4.7) circle(0.04);
\fill[black] (5.3,4.4) circle(0.04);
%てんてん

\draw [densely dashed](P)++(1.1,0) --++(0,-1.5);
\fill [black] (P)++(0.5,0)++(0.6,-1.5) circle(0.06);

\fill [black] (P)++(0.5,0)++(0,-1.5) circle(0.06);
\draw (P)++(0.5,0)++(0,-1.5) node[left]{$p_i$};
%p_iぬりつぶしと名前

\draw (P)++(1.2,-1.5) node[below] {$p'_{i}$};
%p'_i名前

\draw (P) ++ (0.4,-1.3) to[out=120,in=30] ++ (-0.8,0) node[left]{$U_{i}$};
%U_i名前

%%ここから先pのファイバー

\draw (Q)++ (0.6,0) --++ (1.5,1) --++ (0,3.5) --++ (-1.5,-1) -- cycle;
\draw (Q)++(0.3,4) node[above]{$L_{p_{0}}$};
%pのファイバー

\fill [black]  (5.1,5.7) ++(3.5,0) circle (0.08);
\draw (5.1,5.9)++(3.5,0) node[above] {$q_{0,1}$};
%\fill [black] (5.2,5) circle (0.08)node[left] {$q_{i,{m_{i}}}$} ;
%\fill [black] (5.3,4.5) circle (0.08) node[left] {$q_{i,{m_{i}+1}}$};
\fill [black] (5.1,3.6) ++(3.5,0)circle (0.08); 
\draw (4.9,3.6)++(3.5,0) node[left] {$q_{0,n}$};
%交点のぬりつぶし

\fill[black] (4.4,4.7)++(3.5,0) circle (0.08) node[below] {$q_{0}$};
\draw (4.4,4.7)++(3.5,0) -- (8.5,5.42) ;
\draw (5.1,5.7)++(3.5,0) circle (0.3);
\draw (4.4,4.7)++(3.5,0) -- (8.5,3.87);
\draw (5.1,3.6)++(3.5,0) circle (0.3);
%基点とループ

\draw (5.2,5.4)++(3.5,0) to[out=-90,in=-150] ++(0.8,0) node[right]{$e_1$};
\draw (5.2,3.9)++(3.5,0) to[out=90,in=160] ++(0.8,0) node[right]{$e_n$};
%e_1,e_n名前

\fill[black] (5.3,5)++(3.5,0) circle(0.04);
\fill[black] (5.4,4.7)++(3.5,0) circle(0.04);
\fill[black] (5.3,4.4)++(3.5,0) circle(0.04);
%てんてん

\draw [densely dashed](Q)++(1.1,0) --++(0,-1.5);
\fill [black] (Q)++(0.5,0)++(0.6,-1.5) circle(0.06);

%\draw [densely dashed](Q)++(0.5,0) --++(0,-1.5);
%\fill [black] (Q)++(0.5,0)++(0,-1.5) circle(0.06);
\draw (Q)++(0.5,0)++(0.6,-1.5) node[below]{$p_{0}$};
%pまわり

\draw (P) ++ (1.1,-1.5) -- (8.1,1);
\draw (P) ++ (2.5,-1.5) node[below]{$s_i$};
%s_i

\draw (10,0) node[above] {$\bC$};

\end{tikzpicture}
\caption{Non-intersection arc systems in the fibers}
\label{fig:around_p_i}
\end{figure}
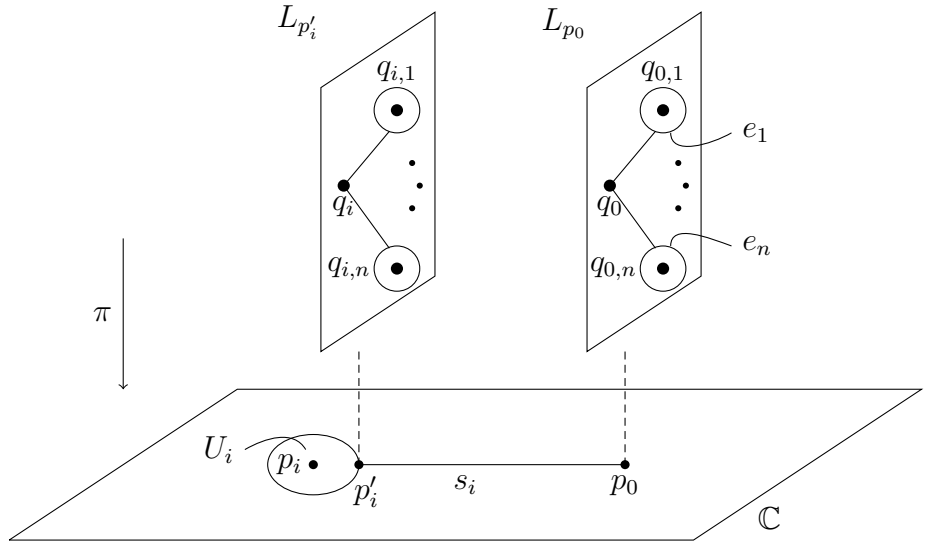

Let $L_{p_{0}} \cap C = \{q_{0,1}, \ldots, q_{0,n}\}$. 
Take a basepoint $q_{0} \in L_{p_{0}} \setminus  C$.
Choose a non-intersecting arc system with respect to the basepoint $q_{0}$ and $\{q_{0,1}, \ldots, q_{0,n}\}$.
Let $e_{1}, \ldots, e_{n} \in \pi_{1} (L_{p_{0}} \setminus  C, q_0)$ be the elements which represent the loops defined from the non-intersecting arc system (Figure \ref{fig:around_p_i}). 
The fundamental group $\pi_{1} (L_{p_{0}} \setminus  C, q_0)$ is free of rank $n$ generated by $e_{1} , \ldots, e_{n}$.

The braid group $B(L_{p_{0}}, L_{p_{0}} \cap C)$ acts on the fundamental group $\pi_{1} (L_{p_{0}} \setminus C, q_0)$ (the so-called Artin action on a free group).

%Let $m_{i}$ be the multiplicity of the singularity of $C$ on $L_{p_{i}}$. 
%If $L_{p_{i}}$ is tangent to $C$, then we define $m_{i} = 2$.
For a point $p \in \partial U_{i}$, the fiber $L_p \cap C$ is a set of $n$ points. We suppose that $m_{i}$ of them collapse into one when the point $p$ reaches $p_{i}$. 
The integer $m_{i}$ coincides with the usual multiplicity of the singularity of $C$ on $L_{p_{i}}$ if $L_{p_{i}}$ intersects $C$ transversally.
Let $L_{p'_{i}} \cap C = \{q_{i,1}, \ldots, q_{i, m_{i}}, q_{i, m_{i}+1}, \ldots ,  q_{i,n} \}$ and we assume that first $m_{i}$ points collapse.
By trivializing along $s_{i}$, each collapsing point $q_{i,j}$ reaches a point in $L_{p_{0}} \cap C$, and denote it by $q_{0,k_{i, j}} \in L_{p_{0}} \cap C$. 
Let $I_{i} = \{k_{i,1},\ldots, k_{i, m_{i} -1}\}$ be the corresponding index set from which the last one is removed.

In this setting, Libgober proved the following theorem.
\begin{theorem} {\rm \cite{lib-hom}}
The complement $M = \bC^2 \setminus C$ is homotopy equivalent to a $2$-dimensional complex associated with the following presentation of the fundamental group:
\[
\pi_{1} (\bC^{2} \setminus C) \cong \langle e_{1} , \ldots, e_{n} \mid \beta_{i} \cdot e_{j} = e_{j} \, (i=1, \ldots, N, j \in I_{i}) \rangle ,
\]
where $\cdot$ denotes the Artin action and $\beta_{i}$ is the monodromy around the loop $\gamma_{i}$.
\end{theorem}

It is convenient to consider the monodromy around each point $p_{i} \in X$ and arc $s_{i}$ separately to compute and describe the Kirby diagrams in the forthcoming sections.
Let $\overline{\beta_{i}} \in B(L_{p'_{i}}, L_{p'_{i}} \cap C)$ the local monodromy around $\partial U_{i}$ and $\phi_{i}: B(L_{p_{0}}, L_{p_{0}} \cap C) \rightarrow B(L_{p'_{i}}, L_{p'_{i}} \cap C)$ be the isomorphism induced by the trivialization along $s_{i}$.

\begin{remark}\label{rmk:braid}
Let us consider the local monodromy around $p_{i} \in X$. Let $\gamma_{i} :[0,1] \rightarrow \partial U_{i}$ be a smooth parametrization with $\gamma_{i}(0) = \gamma_{i} (1) = p'_{i}$ that orients $\partial U_i$ counterclockwise. 
Then, $L_{\gamma_{i}(t)} \cap C$ is expressed as $\{q_{i,1} (t) , \ldots, q_{i,n} (t) \}$, where $q_{i,j} (t)$ is a smooth curve with $q_{i,j} (0) = q_{i,j}$. 
Therefore, $\pi^{-1} (\partial U_{i}) \cap C = \bigcup_{j=1}^{n} \bigcup_{t \in [0,1]} q_{i,j} (t)$ forms a braid with $n$ strands in $\pi^{-1} (\partial U_{i})$.
This braid represents the element which appears as the local monodromy $\overline{\beta_{i}} \in B(L_{p'_{i}}, L_{p'_{i}} \cap C)$.
\end{remark}

\section{Handle decompositions}\label{sec:handle}

In this section, we give a handle decomposition of the complement $M = \bC^2 \setminus C$.
Throughout this section, a manifold is always smooth and, if there are corners, they are appropriately smoothed.

\subsection{1-handlebody $M_1$}

\begin{figure}[htbp]
\centering
\begin{tikzpicture}
\coordinate (P) at (0,0);
\coordinate (Q) at  (0,-3.5);

\draw (Q) ++ (1.5,1.2) --++ (1.2,0.8) --++ (0,4) --++ (-1.2,-0.8) --cycle;
%Qのファイバー

\draw [line width=5pt, white](P) to[out=5,in=210] ++(4,2);
\draw [line width=5pt, white](P) to[out=-5,in=150] ++(4,-2);
\draw [very thick](P) to[out=5,in=210] ++(4,2);
\draw [very thick](P) to[out=-5,in=150] ++(4,-2);
%曲線

\draw [line width=5pt, white](P)++(-0.3,0)to[out=90,in=180] ++(0.3,0.3)[out=0, in=210] to ++(4,2.2);
\draw [line width=5pt, white](P)++(-0.3,0) to[out=270,in=180] ++(0.3,-0.3)[out=0,in=150] to ++(4,-2.2);
\draw (P)++(-0.3,0)to[out=90,in=180] ++(0.3,0.3)[out=0, in=210] to ++(4,2.2);
\draw (P)++(-0.3,0) to[out=270,in=180] ++(0.3,-0.3)[out=0,in=150] to ++(4,-2.2);
%近傍の外側

\draw [line width=5pt, white](P)++(1.3,0) to[out=90,in=200] ++(0.1,0.1) [out=20,in=210]to ++(2.6,1.4);
\draw [line width=5pt, white](P)++(1.3,0) to[out=270,in=160] ++(0.1,-0.1)[out=-20,in=140]to ++(2.6,-1.4);
\draw (P)++(1.3,0) to[out=90,in=200] ++(0.1,0.1) [out=20,in=210]to ++(2.6,1.4);
\draw (P)++(1.3,0) to[out=270,in=160] ++(0.1,-0.1)[out=-20,in=140]to ++(2.6,-1.4);
%近傍の内側

%\fill (Q) ++(2,0)++(0,4.5) circle (0.06);
%\fill (Q) ++(2,0)++(0,3.85) circle (0.06);
%\draw [thick](Q) ++(2,0)++ (0,4.15) circle(0.15cm and 0.32cm);

\draw [densely dashed,thick](Q) ++ (2,0) ++ (0,3.87) [out=0,in=-90] to++(0.15,0.293) [out=90, in=0] to ++(-0.15,0.3);
%Qのファイバーと曲線の交差1右

\draw [line width=4pt, white](Q) ++ (2,0) ++ (0,3.87) [out=180,in=-90] to++(-0.15,0.293) [out=90, in=180] to ++(0.15,0.3);
\draw [thick](Q) ++ (2,0) ++ (0,3.87) [out=180,in=-90] to++(-0.15,0.293) [out=90, in=180] to ++(0.15,0.3);
%Qのファイバーと曲線の交差1左

%\fill (Q) ++(2,0)++(0,3.18) circle (0.06);
%\fill (Q) ++(2,0)++(0,2.5) circle (0.06);
%\draw [thick](Q) ++(2,0)++ (0,2.85) circle(0.15cm and 0.32cm);
\draw [densely dashed, thick](Q) ++ (2,0) ++ (0,2.54) [out=0,in=-90] to++(0.15,0.305) [out=90, in=0] to ++(-0.15,0.31);
%Qのファイバーと曲線の交差2右

\draw [line width=4pt,white](Q) ++ (2,0) ++ (0,2.54) [out=180,in=-90] to++(-0.15,0.305) [out=90, in=180] to ++(0.15,0.31);
\draw [thick](Q) ++ (2,0) ++ (0,2.54) [out=180,in=-90] to++(-0.15,0.305) [out=90, in=180] to ++(0.15,0.31);
%Qのファイバーと曲線の交差2左

\draw[line width=3pt, white] (Q) ++ (1.5,1) ++ (0,0.3) -- ++(0,3.5);
\draw(Q) ++ (1.5,1) ++ (0,0.3) -- ++(0,3.5);
%Qのファイバー塗り直し

\draw[line width = 1.5pt,white] (Q) ++ (1.62,2.5) --++ (0,2);
\draw[line width =1.5pt, white] (Q) ++ (1.72,2.5) --++ (0,2);
\draw[line width =1.5pt, white] (Q) ++ (1.84,3.2) --++ (0,0.7);
%背後の点線づけ

\draw [densely dashed](P) ++ (0,-1) --++ (0,-2.5);
\draw(Q)circle(2cm and 1cm);
\fill[black] (Q) circle (0.06);
\draw (Q) node [left] {$p_i$};
%Pまわり

\fill (Q) ++(2,0) circle (0.06);
\draw  (Q) ++(2,0) node[right] {$q$};
\draw [densely dashed](Q) ++(2,0) --++ (0,1);
%Q

\draw (Q) ++ (-1,0.5) to[out=120,in=30] ++(-1.5,0) node[left]{$U_{i}$};
\draw (P) ++ (4.2,2) node[right] {$\nu(C)$};
\draw (P) ++ (4.1,1.5) to[out=50,in=-50] ++(0,1); 
%C, \nu(C)U_i名前

\end{tikzpicture}
\caption{The neighborhoods $U_i$ and $\nu(C)$}
\label{fig:nbd_of_U}
\end{figure}
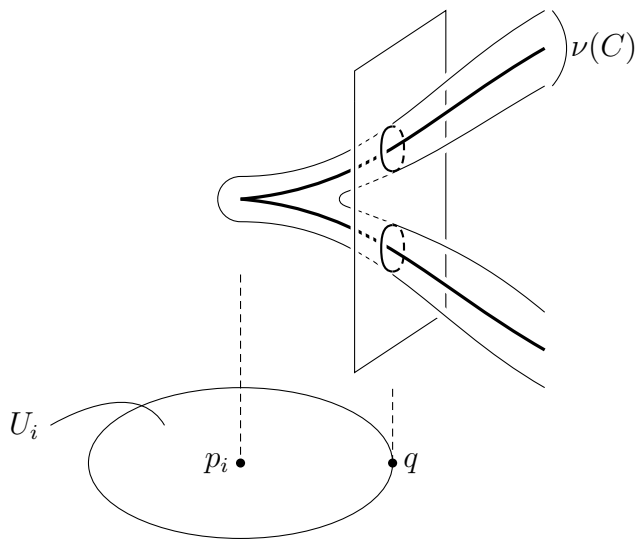

Let $D \subset \bC$ be a sufficiently large closed disk such that a non-intersecting arc system $s_1, s_2, \ldots, s_{N}$ is contained in $D$.
Take a small regular neighborhood $\nu (C)$ of $C$ such that $L_{q} \cap \nu (C)$ consists of $n$ connected components for all $q \in \partial U_{i}$ (see Figure \ref{fig:nbd_of_U}) and $\nu(C)$ deformation rectracts onto $C$ over each radius of $U_i$.
We also assume that $\nu(C)$ is generic so that the restriction of $\pi$ to the boundary  $\partial \nu (C)$ is a stable map. 
We define the disk $D_{y}$ in the fiber complex line such that $L_{p} \cap \nu (C) $ is contained in $\{p\} \times D_{y}$ for all $p \in D$.
Clearly, $D \times D_{y}$ is diffeomorphic to the standard $4$-ball.

\begin{figure}[htbp]
\centering
\begin{tikzpicture}
\coordinate (P) at (0,0);
\coordinate (P0) at (0,-1.2);

\coordinate (Q) at (6,0);
\coordinate (Q0) at (6,-1.2);

\draw (P) circle (2);
%Dの円周
\fill[black] (P0) circle (0.06);
\draw (P0) ++(0,-0.2) node[right] {$p_0$};
%p_0ぬりつぶし、名付け

\draw (P0) --++ (-1.2,1.2);
\filldraw (P0) ++ (-1.2,1.2) circle (0.06) node[right]{$p'_1$};
\draw (P0) ++ (-1.5,1.5) circle (0.4);
\draw (P0) ++ (-0.9,0.9) node[below] {$s_1$};
%s_1

\draw (P0) --++ (-0.5,1.9);
\filldraw (P0) ++ (-0.5,1.9) circle (0.06) node[right]{$p'_2$};
\draw (P0) ++ (-0.6,2.3) circle (0.4);
\draw (P0) ++ (-0.3,1.2) node[right] {$s_2$};
%s_2

\fill (P0) ++ (0.1,2) circle (0.03);
\fill (P0) ++ (0.4,1.95) circle (0.03);
\fill (P0) ++ (0.7,1.85) circle (0.03);
\fill (P0) ++ (1,1.7) circle (0.03);
%てんてん

\draw (P0) --++ (1.2,1.2);
\fill (P0) ++ (1.2,1.2) circle (0.06) ;
\draw (P0) ++ (1.2,1.15) node [below]{$p'_N$};
\draw (P0) ++ (1.5,1.5) circle (0.4);
\draw (P0) ++ (0.4,0.4) node[right] {$s_N$};
%s_N

%\draw (P0) ++ (-0.8,1.6) [out =45, in =120] to++ (2.5,-0.5);

\draw (P0) --++ (0,-0.8);
\draw (P0) ++ (0,-0.4) node[left] {$s_0$};
%s_0

\draw (P0) ++ (-1.7,1.7) to[out=160,in=-20]++ (-0.5,0.5) node[left]{$U_1$};
\draw (P0) ++ (-0.7,2.5) to[out=160,in=-20]++ (-0.5,0.5) node[left]{$U_2$};
\draw (P0) ++ (1.7,1.6) to[out=20,in=200]++ (0.5,0.5) node[right]{$U_N$};
%U_iたち名前

%%%こっから右側の円

\draw (Q) circle (2);
%Dの円周
%\fill[black] (Q0) circle (0.06);
%\draw (Q0) ++(0,-0.2) node[right] {$p_0$};
%p_0ぬりつぶし、名付け

\draw (Q0) ++ (-1.5,1.5) circle (0.4);
\draw (Q0) ++ (-0.6,2.3) circle (0.4);
\draw (Q0) ++ (1.5,1.5) circle (0.4);
%U_iたち

\draw[white,line width =5pt] (Q0) ++ (-0.3,-0.78) --++ (0.6,0);
\draw[white,line width =5pt] (Q0) ++ (-1.2,1.23) --++ (-0.16,-0.16);
\draw[white,line width =5pt] (Q0) ++ (1.18,1.28) --++ (0.19,-0.19);
\draw[white,line width =5pt] (Q0) ++ (-0.39,1.95) --++ (-0.25,-0.07);
%白ぬりつぶし

\draw (Q0) ++ (-1.4,1.11) --++(1.1,-1.1) --++ (0,-0.8);
%左側の辺
\draw (Q0) ++ (-1.2,1.23) --++(0.9,-0.9) --++ (-0.35,1.58);
%左から二番目の辺
\draw (Q0) ++ (-0.39,1.95) --++ (0.35,-1.58);

\draw (Q0) ++ (1.18,1.28) --++ (-1,-1);

\draw (Q0) ++ (1.4,1.11) --++(-1.1,-1.1) --++ (0,-0.8);
%右側の辺
%%近傍の境界たち

\fill (Q0) ++ (0.1,2) circle (0.03);
\fill (Q0) ++ (0.4,1.95) circle (0.03);
\fill (Q0) ++ (0.7,1.85) circle (0.03);
\fill (Q0) ++ (1,1.7) circle (0.03);
%てんてん

\draw (Q0) ++ (0.7,0.5) to[out=0,in=180] ++ (1.3,-0.3)node[right]{$\nu(\Gamma)$};
%nu(\Gamma)名前

\end{tikzpicture}
\caption{The non-intersecting system, the graph $\Gamma$ and its neighborhood $\nu(\Gamma)$ }
\label{fig:graphs}
\end{figure}
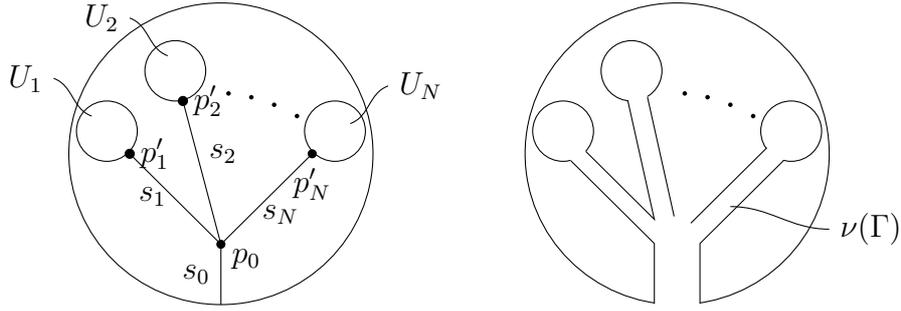

Let $M_{0} = (D \times D_{y}) \setminus \nu (C)$. The manifold $M_{0}$ is a compact $4$-manifold with boundary, and the interior of $M_{0}$ is diffeomorphic to $M$. We shall give a handle decomposition of $M_{0}$.
Let $\pi_{0} :=\pi |_{M_{0}} : M_{0} \rightarrow D$ be the restriction of the projection. 
Note that the restriction of $\pi_{0}$ to $M_{0} \setminus \bigcup_{i=1}^{N} \pi_{0}^{-1} (U_{i}) $ is again a locally trivial fibration. 
We prepare some notation (see Figure \ref{fig:graphs}). 
\begin{itemize}
\item $s_0$ : a non-intersecting arc connecting $p_0$ and a point on $\partial D$ such that $s_{i} \cap s_{0} = \{p_{0}\}$ for $i=1, \ldots, N$.
\item $\Gamma = \bigcup_{i=0}^{N} s_{i}$.
\item $\nu (s_i)$ : a small tubular neighborhood of $s_{i}$ for $i=0,1, \ldots, N$.
\item $\nu (\Gamma) = \bigcup_{i=0}^{N} \nu (s_i) \cup \bigcup_{i=1}^{N} U_i$.
\item $M_1 = \pi_{0}^{-1} (D \setminus \nu (\Gamma))$.
\item $D_{k} = D^2 \setminus \{\mbox{$k$ small disks}\}$.
\end{itemize}
%Let $s_0$ be a non self intersecting arc connecting $p_0$ and a point at $\partial U_i$ and $\Gamma = \bigcup_{i=0}^{N} s_{i}$. For each $s_i$, let us denote $\nu (s_i)$ be a small tubular neighborhood. 

Under this setting, $D^2 \setminus \nu (\Gamma)$ is diffeomorphic to $D^2$, and thus contractible. Therefore, the restriction $\pi_{0}: M_1 \rightarrow D \setminus \nu (\Gamma)$ is a trivial fibration. Each fiber is diffeomorphic to $D_n$, and hence, $M_1 \cong D_{n} \times D^2 $. %(recall that $n= \deg C$). 
Thus, we have the following proposition.

\begin{proposition}\label{prop:1-handle}
$M_1$ is diffeomorphic to a $4$-manifold obtained from a $4$-ball $B^4$ by attaching $n$ $1$-handles.
\end{proposition}

\subsection{The stable map on $\overline{U_i \setminus V_i}$ and the attachment of $\pi^{-1}_0 (A \cup V_{i})$}

Since $M_0 = M_1 \cup \pi_{0}^{-1} (\nu(\Gamma) )$, we need to see how $\pi_{0}^{-1} (\nu(\Gamma))$ is attached to $M_1$. We will describe the attachment of $\pi_{0}^{-1} (\nu(\Gamma))$ by decomposing it into some pieces. 
Here, we provide the following lemma, which is needed to describe the change of topology when each part is attached. 
The proof is easy, and we can check it by hand.

\begin{figure}[htbp]
\centering
\begin{tikzpicture}

\coordinate (A) at (0,0);
\coordinate (B) at (5,0);
\coordinate (C) at (10,0);
\draw (A) arc [start angle=180,end angle=360,x radius=2,y radius=0.5] --++(0,1) arc [start angle=360,end angle=180,x radius=2,y radius=0.5] --++(0,-1);
\draw (A) ++(0,1) arc [start angle=180,end angle=0,x radius=2,y radius=0.5];
\draw[densely dashed] (A) arc [start angle=180,end angle=0,x radius=2,y radius=0.25];
%X
\draw (A) ++ (2,1) circle[x radius =0.8,y radius = 0.25];
%Y \subset \partial X

\draw (A) ++(1.2,1) ++ (0,1) arc [start angle=180,end angle=360,x radius=0.8,y radius=0.25] --++(0,0.5) arc [start angle=360,end angle=180,x radius=0.8,y radius=0.25] --++(0,-0.5);
\draw (A) ++ (1.2,1)++(0,1)++(0,0.5) arc [start angle=180,end angle=0,x radius=0.8,y radius=0.25];
\draw[densely dashed] (A) ++ (1.2,1)++(0,1) arc [start angle=180,end angle=0,x radius=0.8,y radius=0.1];
%Y \times [0,1]

\draw (A) ++ (0.5,0.3) to[out=200,in=20] ++(-0.9,-0.5) node[left]{$X$};
\draw (A) ++ (2.3,1) to[out=60,in=200] ++(1.3,0.8)node[right]{$Y$};
\draw (A) ++ (1.5,2.1) to[out=180,in=0]++(-0.7,0.3) node[left]{$Y \times [0,1]$};
%%%ここから右側

\draw (B) arc [start angle=180,end angle=360,x radius=2,y radius=0.5] --++(0,1) arc [start angle=360,end angle=180,x radius=2,y radius=0.5] --++(0,-1);
\draw (B) ++(0,1) arc [start angle=180,end angle=0,x radius=2,y radius=0.5];
\draw[densely dashed] (B) arc [start angle=180,end angle=0,x radius=2,y radius=0.25];
%X
\draw (B) ++ (1.2,1) arc [start angle=180,end angle=360,x radius=0.8,y radius=0.25];
%Y \subset \partial X

\fill[white] (B) ++ (1.2,1) ++(0,0.25)++(0.8,0.25) circle [x radius=0.8,y radius=0.25];
\draw (B) ++(1.2,1)  arc [start angle=180,end angle=360,x radius=0.8,y radius=0.25] --++(0,0.5) arc [start angle=360,end angle=180,x radius=0.8,y radius=0.25] --++(0,-0.5);
\draw (B) ++ (1.2,1)++(0,0.5) arc [start angle=180,end angle=0,x radius=0.8,y radius=0.25];
%\draw[densely dashed] (B) ++ (1.2,1)++(0,1) arc [start angle=180,end angle=0,x radius=0.8,y radius=0.1];
% Yの接着
\draw (B) ++ (3.2,1.2) to[out=60,in=270]++(1,0.7)node[above]{$X \cup_{Y} (Y \times [0,1])$};

\draw (C) arc [start angle=180,end angle=360,x radius=2,y radius=0.5] --++(0,1) arc [start angle=360,end angle=180,x radius=2,y radius=0.5] --++(0,-1);
\draw (C) ++(0,1) arc [start angle=180,end angle=0,x radius=2,y radius=0.5];
\draw[densely dashed] (C) arc [start angle=180,end angle=0,x radius=2,y radius=0.25];
%X
\draw (C) ++ (2,1) circle[x radius =0.8,y radius = 0.25];
\draw[densely dashed] (C) ++ (1.2,1) --++(0,-0.4);
\draw[densely dashed] (C) ++ (2.8,1) --++(0,-0.4);
\draw[densely dashed] (C) ++(2,0.6) circle [x radius=0.8, y radius =0.25]; 
\draw[preaction={line width=1.5 pt, draw = white}] (C) ++(4,1) arc [start angle=0,end angle=-180,x radius=2,y radius=0.5];
\draw (C) ++ (3,1.2) to[out=60,in=270]++(1,0.7)node[above]{$X \setminus (\Int(Y) \times (0,1))$};

\end{tikzpicture}
\caption{Attachment and removal of $Y \times [0,1]$. In this example, $X \cong D^3$ and $Y \cong D^2$.}

\label{fig:attachinglemma}
\end{figure}
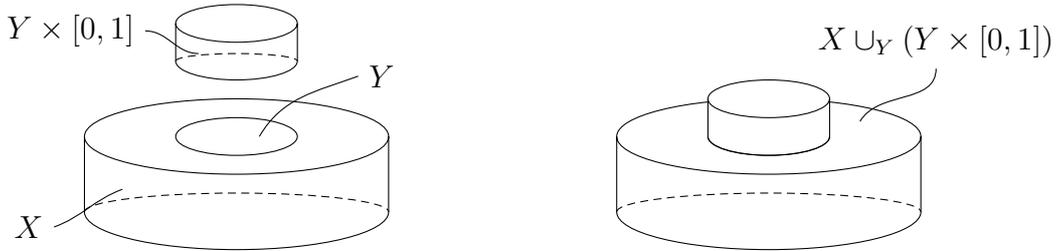

\begin{lemma} \label{lem:handle}
Let $X$ and $Y$ be compact manifolds with boundary, such that $Y \subset \partial X$ and $\dim Y = \dim \partial X$. 
Then, we have the following (see Figure \ref{fig:attachinglemma}).
\begin{enumerate}[(1)]
\item The manifold $X$ is diffeomorphic to the manifold $X\cup_{Y} (Y \times [0,1])$ which is obtained from $X$ by attaching $Y\times [0,1]$ along $Y \times \{0\} \cong Y$.
In other words, attaching $Y \times [0,1]$ to $X$ does not change the diffeomorphism type of $X$. 

\item Suppose that $Y \times [0,1]$ is embedded in $X$ with $(Y \times [0,1]) \cap \partial X = Y \times \{1\}$. Then, $X$ is diffeomorphic to the manifold $X \setminus (\Int(Y) \times (0,1])$. In other words, removing $Y \times [0,1]$ from $X$ does not change the diffeomorphism type of $X$. 
\end{enumerate}
\end{lemma}

\begin{figure}[htbp]
\centering
\begin{tikzpicture}
\coordinate (P) at (0,0);
\coordinate (Q) at  (0,-3.5);

\draw (Q) ++ (0.5,1.2) --++ (1.2,0.8) --++ (0,4) --++ (-1.2,-0.8) --cycle;
%Qのファイバー

\draw [line width=5pt, white](P) to[out=5,in=210] (4,2);
\draw [line width=5pt, white](P) to[out=-5,in=150] (4,-2);
\draw [very thick](P) to[out=5,in=210] (4,2);
\draw [very thick](P) to[out=-5,in=150] (4,-2);
%曲線

\draw [line width=5pt, white](P)++(-0.3,0)to[out=90,in=180] ++(0.3,0.3)[out=0, in=210] to ++(4,2.2);
\draw [line width=5pt, white](P)++(-0.3,0) to[out=270,in=180] ++(0.3,-0.3)[out=0,in=150] to ++(4,-2.2);
\draw (P)++(-0.3,0)to[out=90,in=180] ++(0.3,0.3)[out=0, in=210] to ++(4,2.2);
\draw (P)++(-0.3,0) to[out=270,in=180] ++(0.3,-0.3)[out=0,in=150] to ++(4,-2.2);
%近傍の外側

\draw [line width=5pt, white](P)++(1.3,0) to[out=90,in=200] ++(0.1,0.1) [out=20,in=210]to ++(2.6,1.4);
\draw [line width=5pt, white](P)++(1.3,0) to[out=270,in=160] ++(0.1,-0.1)[out=-20,in=140]to ++(2.6,-1.4);
\draw (P)++(1.3,0) to[out=90,in=200] ++(0.1,0.1) [out=20,in=210]to ++(2.6,1.4);
\draw (P)++(1.3,0) to[out=270,in=160] ++(0.1,-0.1)[out=-20,in=140]to ++(2.6,-1.4);
%近傍の内側

%\fill (Q) ++(2,0)++(0,4.5) circle (0.06);
%\fill (Q) ++(2,0)++(0,3.85) circle (0.06);
%\draw [thick](Q) ++(1,0)++ (0,3.5) circle(0.15cm and 0.47cm);
%Qのファイバーと曲線の交差

\draw[densely dashed,thick] (Q) ++ (1,0) ++ (0,3.04) [out=0,in=-90] to++(0.17,0.5) [out=90, in=0] to ++(-0.17,0.41);
%Qのファイバーと曲線の交差右

\draw [line width =4pt, white](Q) ++ (1,0) ++ (0,3.04) [out=180,in=-90] to++(-0.17,0.5) [out=90, in=180] to ++(0.17,0.41);
\draw [thick](Q) ++ (1,0) ++ (0,3.04) [out=180,in=-90] to++(-0.17,0.5) [out=90, in=180] to ++(0.17,0.41);
%Qのファイバーと曲線の交差左

\draw[line width=3pt, white] (Q) ++ (0.5,1) ++ (0,0.3) -- ++(0,3.5);
\draw(Q) ++ (0.5,1) ++ (0,0.3) -- ++(0,3.5);
%Qのファイバー塗り直し

\draw[line width = 1.5pt,white] (Q) ++ (0.65,2.5) --++ (0,2);
%\draw[line width = 1.1pt,red] (Q) ++ (0.7,2.5) --++ (0,2);
%背後の点線づけ

\draw [densely dashed](P) ++ (0,-1) --++ (0,-2.5);
\fill[black] (Q) circle (0.06);
\draw (Q) node [left] {$p_i$};
%Pまわり

\draw[densely dotted](Q)circle(2cm and 1cm);
%U_i

\draw(Q)circle(1cm and 0.5cm);

\fill (Q) ++(1,0) circle (0.06);
\draw  (Q) ++(1,0) node[right] {$q$};
\draw [densely dashed](Q) ++(1,0) --++ (0,1);
%Q

\end{tikzpicture}
\caption{The neighborhood $V_i$}
\label{fig:nbd_of_V}
\end{figure}
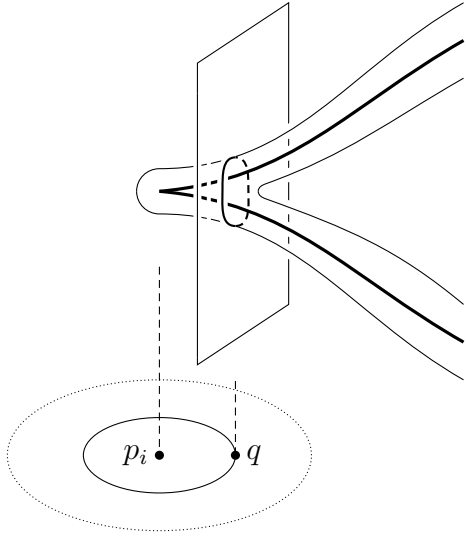

For each $p_i$, let $V_i \subset U_i$ be a small disk around $p_i$ such that $\pi_{0}^{-1} (a) \cong D_{n-(m_{i}-1)}$ for all $a \in V_i$ (see Figure \ref{fig:nbd_of_V}).
For $a \in \partial U_{i}$, we denote by $\ell_{a}$ the intersection of the segment connecting $p_i$ with $a$ and $\overline{U_{i} \setminus V_{i}}$. 
We will study the stable map $\pi_0 |_{\partial M_0}$ restricted to $\overline{U_i \setminus V_i}$.

\begin{proposition}\label{prop:stable-1}
If necessary, by shrinking the neighbourhood $\nu(C)$ in $\pi^{-1} (\overline{U_i \setminus V_i})$, we may assume that the stable map given by the restriction $\pi_{0} |_{\partial M_0}$ has only indefinite folds as critical points on $\pi^{-1}_{0} (\overline{U_i \setminus V_i})$.
\end{proposition}

\proof
Recall that the critical points of a stable map from a $3$-manifold to a $2$-manifold are classified into three types: definite folds, indefinite folds, and cusps \cite{lev}.
The image of a cusp appears as an intersection of the images of a definite fold and an indefinite fold (see Figure 2 in \cite{sae}, for example).
Thus, to show that there are only indefinite folds, it suffices to eliminate the definite folds.
Assume that there is a definite fold singularity $q \in (\pi_0 |_{\partial M_0})^{-1} (\overline{U_i \setminus V_i})$. 
Let $p = \pi_{0} (q)$ and assume $p$ lies in a segment $\ell_a$. We can generically assume that $\ell_a$ intersects the set of critical values transversally.
When crossing a definite fold singularity, the fiber of $\pi_{0}$ changes by the collapse of a circle to a point or the reverse process (emergence of a circle) (Figure \ref{fig:definite}).
Let $\delta_\varepsilon \subset (\pi_{0} |_{\partial M_0})^{-1} (p + \varepsilon)$ be the collapsing (or emerging) circle component, where $\varepsilon >0$ is a sufficiently small value. 
The circle $\delta_{\varepsilon}$ bounds a disk $\Delta_\varepsilon$ in $\pi_0^{-1} ([p, p+\varepsilon]) \cap \partial \nu (C)$. 
%Note that the disk $\Delta_\varepsilon$ does not intersect the curve $C$, since the curve $C$ defines a multi-section of the projection $\pi$ over $\ell_a$. 
The restriction of $\pi_{0} |_{\partial M_{0}}$ to $(\pi_0 |_{\partial M_0})^{-1} (\ell_a)$ is a Morse function on a surface, and the definite fold points correspond to the critical points of indices $0$ or $2$. 
As the value $\varepsilon$ increases, the circle $\delta_\varepsilon$ meets a circle component of some other fiber of $\pi_0 |_{\partial M_0}$ (in Figure \ref{fig:definite} (a), the fiber $\pi^{-1}_0 (p')$).
Otherwise, the circle $\delta_\varepsilon$ would either coincide with a circle component of a fiber of $\partial U_i$ (Figure \ref{fig:definite} (b)), or be capped off at a paired definite fold point and appear as a surface component belonging to a different connected component (Figure \ref{fig:definite} (c)). 
Consider the former case. Let $\delta$ be the resulting circle. 
The circle $\delta_\varepsilon$ bounds a disk $\Delta_\varepsilon$ in $\pi_0^{-1} ([p,p+\varepsilon]) \cap \partial \nu (C)$, and this holds even if the value $\varepsilon$ is increased.
So, in this case, the resulting circle $\delta$, which is a component of $\pi^{-1}_0 (a) \cap \partial \nu (C)$, would bound a disk in $\pi_0^{-1} (\ell_a) \cap \partial \nu (C)$.
Besides, since $\nu (C)$ deformation retracts to $C$ over each radius of $U_i$, each connected component of $\pi^{-1} (\ell_a) \cap \nu (C)$ is homeomorphic to $D^3$.
Each connected component of $\pi_0^{-1} (a) \cap \partial \nu (C)$ represents a non-trivial element in $H_1(\pi^{-1} (\ell_a) \cap (\nu (C) \setminus C); \bZ)$, since the curve $C$ defines a multi-section of the projection $\pi$ over $\ell_a$ (like the ``meridian'' of each component of $C$).
Thus it does not bound any surface in $\pi^{-1} (\ell_a) \cap (\nu (C) \setminus C)$.
In particular, the resulting circle $\delta$ does not bound any disk in $\pi_0^{-1}(\ell_a) \cap \partial \nu (C)$.
This yields a contradiction.
In the latter case, a new connected component in $\pi_{0}^{-1}(\ell_a) \cap \partial \nu(C)$ appears; however, this does not happen since we took the neighbourhood $\nu(C)$ such that this deformation retracts to $C$ in the preimage of each radius of $U_i$.

%Therefore, when a definite fold appears, the  
The intersection point of the two circle components and the original definite fold critical point forms a canceling pair of a Morse function. By shrinking 
$\nu(C)$, we can eliminate the definite fold critical point.
Thus, the proposition follows.
%Otherwise, then by increasing 
%$\varepsilon$, the circle $\delta$ would either coincide with a circle component of a fiber of $\partial U_i$, or be capped off at a paired definite fold point and appear as a surface component belonging to a different connected component.
\endproof

\begin{figure}[htbp]
\centering
\begin{tikzpicture}[scale=0.8]
\coordinate (Q1) at (5,0);
\coordinate (Q2) at (5,0);
\coordinate (R1) at (0,-8);
\coordinate (R2) at (8,-8);
\coordinate (S1) at (10,-8);
\coordinate (S2) at (18,-8);

%%a%%

\draw (Q1) --++ (8,0);
\draw (Q1) ++ (-0.5,-1) --++(8,0);
\filldraw (Q1) ++ (-0.2,-0.4) --++(8,0) circle (0.04) node[right]{$a$};
\draw (Q1) ++ (-0.2,-0.4)  node[left]{$\ell_a$};
%底空間

\draw (Q1) ++ (0,3) circle [x radius=0.2,y radius=0.5];
%left fiber

\draw (Q1) ++ (8,3) ++(0,2) circle [x radius=0.2,y radius=0.5];
\draw (Q1) ++ (8,3) ++(0,-0.5) circle [x radius=0.2,y radius=0.5];
%right fiber

\draw (Q1) ++ (0,3) ++ (0,0.5) to[out=0,in=180] ++(8,2);
\draw (Q1) ++ (0,3)++(0,-0.5) to[out=0,in=90] ++ (4,0);
\draw (Q1) ++ (0,3) ++ (2.5,-1.5)to[out=90,in=-90] ++ (1.5,1);
\draw (Q1) ++ (0,3) ++ (2.5,-1.5)to[out=-90,in=190] ++ (3.5,0.5) to[out=10,in=180] ++(2,0);
\draw (Q1) ++ (8,3)++(0,1.5) to[out=180, in=90] ++ (-2,-1) to[out=-90, in=180] ++ (2,-0.5);
%neighbourhood

\draw (Q1) ++ (0,3) ++(-0.3,0) node[left]{$C$};
\draw (Q1) ++ (0,3) ++(-0.3,0.1) to[out=0, in=180] ++ (8.6,2);
\draw (Q1) ++ (0,3) ++(-0.3,-0.1) to[out=0, in=180] ++ (3.6,0.5)to[out=0, in=180] ++ (5,-0.9);
%Curve

\filldraw (Q1) ++ (2.5,1.5) circle (0.05) node[left]{$q$};
\draw[densely dotted, thick] (Q1) ++ (3.2,1.5) circle[x radius =0.2, y radius=0.5];
\draw (Q1) ++ (3.5,0.8) node{$\delta_{\varepsilon}$};
\draw (Q1) ++ (2.8,1.6) to[out=120, in=180] ++(-0.3,0.9) node[right]{$\Delta_\varepsilon$};
 
\filldraw[densely dashed] (Q1) ++(2.5,1.5) --++(0,-1.9) circle (0.05) node[below]{$p$};
\filldraw[densely dashed] (Q1) ++(4,3) --++(0,-3.4) circle (0.05) node[below]{$p'$};

%臨界点の点線

\draw (Q1) ++ (-1,1) node[left]{$\pi_0$};
\draw[->] (Q1) ++ (-1,1.5) --++ (0,-1);

\draw (Q1) ++ (4,-2) node{(a): A removable critical point};

%%aおわり%%

%%b%%%

\draw (R1) --++ (8,0);
\draw (R1) ++ (-0.5,-1) --++(8,0);
\filldraw (R1) ++ (-0.2,-0.4) --++(8,0) circle (0.04) node[right]{$a$};
\draw (R1) ++ (-0.2,-0.4)  node[left]{$\ell_a$};
%底空間

\draw (R1) ++ (0,3) circle [x radius=0.2,y radius=0.5];
%left fiber

\draw (R1) ++ (8,3) ++(0,2) circle [x radius=0.2,y radius=0.5];
\draw (R1) ++ (8,3) ++(0,0) circle [x radius=0.2,y radius=0.5];
\draw (R1) ++ (8,3) ++(0,-1.5) circle [x radius=0.2,y radius=0.5];
%right fiber

\draw (R1) ++ (0,3) ++ (0,0.5) to[out=0,in=180] ++(8,2);
\draw (R1) ++ (0,3)++(0,-0.5) to[out=0,in=180] ++ (4,0.2) to[out=0, in =180] ++(4,-0.2);
\draw (R1) ++ (8,3)++(0,1.5) to[out=180, in=90] ++ (-1.5,-0.6) to[out=-90, in=180] ++ (1.5,-0.4);
\draw (R1) ++ (8,2) to[out=180, in =90] ++ (-2,-0.5) to[out=-90, in =180] ++ (2,-0.5);
%neighbourhood 

\draw (R1) ++ (0,3) ++ (-0.3,0) node[left]{$C$};
\draw (R1) ++ (0,3) ++ (-0.3,0.1) to[out=0, in=180] ++(8.6,2);
\draw (R1) ++ (0,3) ++ (-0.3,-0.1) to[out=0, in=180] ++(4,0.3)to[out=0, in=180] ++(4.6,-0.3);
\draw[dashed] (R1) ++ (8.3,1.5) --++(-0.8,0);
%Curve

\filldraw (R1) ++ (6,1.5) circle (0.05) node[left]{$q$};
%\draw[densely dotted, thick] (R1) ++ (3.2,1.5) circle[x radius =0.2, y radius=0.5];
%\draw (R1) ++ (3.5,0.8) node{$\delta_{\varepsilon}$};

\filldraw[densely dashed] (R1) ++(6,1.5) --++(0,-1.9) circle (0.05) node[below]{$p$};
%臨界点の点線

\draw (R1) ++ (-0.5,1) node[left]{$\pi_0$};
\draw[->] (R1) ++ (-0.5,1.5) --++ (0,-1);
\draw (R1) ++ (8.2,1.2) --++(0.1,-0.3) node[below]{$\delta$};

\draw (R1) ++ (4,-2) node{(b): NG};

%%bおわり%%

%%c%%%

\draw (S1) --++ (8,0);
\draw (S1) ++ (-0.5,-1) --++(8,0);
\filldraw (S1) ++ (-0.2,-0.4) --++(8,0) circle (0.04) node[right]{$a$};
\draw (S1) ++ (-0.2,-0.4)  node[left]{$\ell_a$};
%底空間

\draw (S1) ++ (0,3) circle [x radius=0.2,y radius=0.5];
%left fiber

\draw (S1) ++ (8,3) ++(0,2) circle [x radius=0.2,y radius=0.5];
\draw (S1) ++ (8,3) ++(0,0) circle [x radius=0.2,y radius=0.5];
%right fiber

\draw (S1) ++ (0,3) ++ (0,0.5) to[out=0,in=180] ++(8,2);
\draw (S1) ++ (0,3)++(0,-0.5) to[out=0,in=180] ++ (4,0.2) to[out=0, in =180] ++(4,-0.2);
\draw (S1) ++ (8,3)++(0,1.5) to[out=180, in=90] ++ (-1.5,-0.6) to[out=-90, in=180] ++ (1.5,-0.4);
\draw (S1) ++ (5,2) to[out=180, in =90] ++ (-2,-0.5) to[out=-90, in =180] ++ (2,-0.5);
\draw (S1) ++ (5,2) to[out=0, in =90] ++ (2,-0.5) to[out=-90, in =0] ++ (-2,-0.5);
%neighbourhood 

\draw (S1) ++ (0,3) ++ (-0.3,0) node[left]{$C$};
\draw (S1) ++ (0,3) ++ (-0.3,0.1) to[out=0, in=180] ++(8.6,2);
\draw (S1) ++ (0,3) ++ (-0.3,-0.1) to[out=0, in=180] ++(4,0.3)to[out=0, in=180] ++(4.6,-0.3);
%Curve

\filldraw (S1) ++ (3,1.5) circle (0.05) node[left]{$q$};
%\draw[densely dotted, thick] (S1) ++ (3.2,1.5) circle[x radius =0.2, y radius=0.5];
%\draw (S1) ++ (3.5,0.8) node{$\delta_{\varepsilon}$};

\filldraw[densely dashed] (S1) ++(3,1.5) --++(0,-1.9) circle (0.05) node[below]{$p$};
%臨界点の点線

\draw (S1) ++ (-0.5,1) node[left]{$\pi_0$};
\draw[->] (S1) ++ (-0.5,1.5) --++ (0,-1);

\draw (S1) ++ (4,-2) node{(c): NG};

%%cおわり%%

\end{tikzpicture}
\caption{Critical points in $\pi^{-1}_0 (\ell_a)$}
\label{fig:definite}
\end{figure}

From now on, we consider the neighborhood $\nu(C)$ satisfying this proposition. Since there are only indefinite fold critical points, the critical value set in $\overline{U_{i} \setminus V_i}$ is an immersed $1$-dimensional manifold, and we denote it by $S$ (Figure \ref{fig:criticalvalue}).

%Also, from now on, we are interested in only the connected components of $\nu(C)$ and $\partial \nu(C)$ that are related to the singular points of $C$ contained in $\pi^{-1} (p_i)$, locally. 
%By a slight abuse of notation, we shall henceforth use the same notation $\nu(C) \cap \pi_0^{-1}(\ell_0)$ to denote the connected components of $\nu(C) \cap \pi_0^{-1}(\ell_0)$ corresponding to the singular points.

\begin{figure}[htbp]
\centering
\begin{tikzpicture}
\coordinate (P) at (0,0);

\fill [black] (P) circle (0.06);
\draw (P) node[below] {$p_i$};
%p_i

\draw[thick] (P)++(1.2,0) to[out=90,in=0] ++(-1.2,1.5)to[out=180,in=90] ++(-1.8,-1.5)to[out=-90,in=180] ++(1.8,-1.8) to[out=0,in=-90] ++ (1.5,1.8);
\draw[thick] (P)++(1.5,0) to[out=90,in=0] ++(-1.5,1.2)to[out=180,in=90] ++(-1.5,-1.2)to[out=-90,in=180] ++(1.5,-1.5)to[out=0,in=-90] ++ (1.2,1.5);;
\draw[thick] (P)++(1.8,0) to[out=90,in=0] ++(-1.8,1.8)to[out=180,in=90] ++(-1.2,-1.8)to[out=-90,in=180] ++(1.2,-1.2)to[out=0,in=-90] ++ (1.8,1.2);
\draw (P) ++ (1.5,1) to[out=60,in=180]++(0.5,0.3)node[right]{$S$};
%critical values

\draw (P) ++ (-0.86,-0.5) --++(-0.86,-0.5) node[left]{$a$};
\fill (P) ++ (-0.86,-0.5) ++(-0.86,-0.5) circle(0.04);
\draw[preaction={draw = white, line width=3 pt}] (P) ++ (-1,-0.58) ++(-0.43,-0.25) to[out=120,in=0] ++(-0.7,0.3) node[left]{$\ell_{a}$};
%\ell_aとa

\draw (P) circle (1);
\draw (P)++(-0.5,0) node[above]{$V_i$};
%V_i
\draw (P) circle (2);
%U_i

%\draw (P) --++ (1.9,0.6) --++ (4,0);
%\draw (P) --++ (1.9,-0.6) --++ (4,0);
\draw [densely dashed](P) --++ (2,0) node[right]{$p'_i$};
%\draw (P) ++(1.9,0) ++(4,0) node[right]{$s_i$};
%s_iの近傍

\end{tikzpicture}
\caption{The critical value set}
\label{fig:criticalvalue}
\end{figure}

\begin{proposition}\label{prop:stable-2}
For each $a \in \partial U_i$, there are $(m_i -1)$ critical points of $\pi_{0} |_{\partial M_0}$ on $\pi_{0}^{-1} (\ell_a)$.
\end{proposition}

\proof
For generic $a \in \partial U_i$, $\ell_a$ intersects the critical value set transversally. For such $a \in \partial U_i$, the restriction of $\pi_{0} |_{\partial M_0}$ to $\pi_0^{-1}(\ell_a)$ is a Morse function.
Let us consider such $a \in \partial U_i$ at first.
Since $\nu(C)$ deformation retracts to $C$ in the preimage of $\pi$ over each radius of $U_i$, $\nu(C) \cap \pi^{-1} (\ell_a)$ is homeomorphic to a disjoint union of $(n-m_i +1)$ $3$-balls, and each boundary is homeomorphic to $S^2$. 
Now, we have 
\[
\partial (\nu(C) \cap \pi^{-1} (\ell_a)) = (\pi^{-1} (\partial \ell_a) \cap \nu(C)) \cup (\pi^{-1}_{0} (\ell_a) \cap \partial \nu(C)),
\]
and $\pi^{-1} (\partial \ell_a) \cap \nu(C)$ consists of $(n - m_i + 1) + n$ disks (Figure \ref{fig:indexone}). Therefore, $\pi^{-1}_{0} (\ell_a) \cap \partial \nu(C)$ is a disjoint union of $(n-m_i + 1)$ $2$-spheres with $(2n-m_{i} +1)$ boundary components in total.
Thus, the relative Euler characteristic satisfies 
\[
\chi (\pi^{-1}_0 (\ell_a) \cap \partial \nu (C), \pi^{-1}_0 (\partial \ell_a \cap V_i) \cap \partial \nu (C) ) = 1 -m_i.
\]
By Proposition \ref{prop:stable-1}, the Morse function obtained as the restriction of $\pi_{0}$ to $\pi_0^{-1} (\ell_a)$ has only critical points of index $1$. Therefore, by the Euler-Poincar\'e formula in Morse theory, there are $(m_i-1)$ critical points in $\pi_0^{-1}(\ell_a)$.

We initially assumed that the statement holds for a generic $a \in \partial U_i$, but in fact it holds for every $a \in \partial U_i$.
Suppose that there exists $a \in \partial U_i$ for which it does not hold. 
Then, the corresponding segment $\ell_a$ would be tangent to $S$. We may assume that this tangency is quadratic (after a suitable coordinate change in the base space, if necessary).
In that case, for the values close to $a$, the number of the intersection points between the segment $\ell$ and $S$ would differ by two.
However, the argument given above on the number of critical points for generic $a$ does not depend on the choice of $a$. This is a contradiction. Therefore, for every $a \in \partial U_i$, the segment $\ell_a$ intersects the critical values set transversally, and the proposition follows. 
\endproof

\begin{figure}[htbp]
\centering
\begin{tikzpicture}[scale=0.8]
\coordinate (Q1) at (0,0);
\coordinate (Q2) at (0,-3.5);

%%a%%

\draw (Q2) --++ (8,0);
\draw (Q2) ++ (-0.5,-1) --++(8,0);
\draw (Q2) ++ (-0.2,-0.4) --++(8,0) node[right]{$\partial U_i$};
\draw (Q2) ++ (-0.2,-0.4)  node[left]{$V_i$};
\draw (Q2) ++(4,-0.3)node[below]{$\ell_a$};
%底空間

\draw (Q1)++(0,0.5) circle [x radius=0.2,y radius=0.5];
\draw (Q1) ++ (8,0)++(0,0.5) circle [x radius=0.2,y radius=0.5];
\draw (Q1) ++(0,0.5)++(0,0.5) --++(8,0);
\draw (Q1) ++(0,0.5)++(0,-0.5) --++(8,0);

\draw (Q1) ++ (4,-0.5) node{$\vdots$};

\draw (Q1)++(0,-1.5) circle [x radius=0.2,y radius=0.5];
\draw (Q1) ++ (8,0)++(0,-1.5) circle [x radius=0.2,y radius=0.5];
\draw (Q1) ++(0,-1.5)++(0,0.5) --++(8,0);
\draw (Q1) ++(0,-1.5)++(0,-0.5) --++(8,0);
%irrelevant fibers

\draw (Q1) ++ (0,3) circle [x radius=0.2,y radius=0.5];
\draw (Q1) ++ (-0.5,1) node[left]{$\pi^{-1}_0 (\partial \ell_a \cap V_i) \cap \nu(C)$};
\draw (Q1) ++ (-0.5,0.3) node[left]{($(n-m_i +1)$ disks)};
%left fiber

\draw (Q1) ++ (8,3) ++(0,2) circle [x radius=0.2,y radius=0.5];
\draw (Q1) ++ (8,3) ++(0,-1) circle [x radius=0.2,y radius=0.5];
\draw (Q1) ++ (8,0) ++ (3,1) node{$\pi^{-1}_0 (\partial \ell_a \cap \partial U_i) \cap \nu(C)$};
\draw (Q1) ++ (8,0) ++ (3,0.3) node{($n$ disks)};
%right fiber

\draw (Q1) ++ (0,3) ++ (0,0.5) to[out=0,in=180] ++(8,2);
\draw (Q1) ++ (0,3) ++(0,-0.5) to[out=0,in=180] ++ (8,-1);
\draw (Q1) ++ (0,3) ++(8,1.5) to[out=180,in=30] ++ (-0.7,-0.3);
\draw (Q1) ++ (0,3) ++(8,-0.5) to[out=180,in=-30] ++ (-0.7,0.3);
\draw (Q1) ++ (0,3) ++(8,0.5) node{$\vdots$}; 
%neighbourhood

\draw (Q1) ++ (0,-1.5) ++(4,-0.7) node[below]{$\pi^{-1}(\ell_a) \cap \partial \nu(C)$};

\draw (Q1) ++ (-1,-3) node[left]{$\pi_0$};
\draw[->] (Q1) ++ (-1,-2.5) --++ (0,-1);
%%aおわり%%

\end{tikzpicture}
\caption{The topology of $\partial (\nu(C) \cap \pi^{-1} (\ell_a))$}
\label{fig:indexone}
\end{figure}

The critical value set $S$ is the union of immersed circles in $\overline{U_i \setminus V_i}$. Thus, for generic $a \in \partial U_i$, there are no self-intersections of $S$ on $\ell_a$.
By performing a coordinate change if necessary, we may assume from the outset that there are no self-intersections of $S$ in $\nu(s_i)$.

\begin{definition}\label{def:a&b}
We set $c_i = \partial U_i \setminus \nu(s_i)$, $A = \bigcup_{a \in c_i} \ell_a$, and $B = \overline{U_i \setminus (A \cup V_i)}$ (Figure \ref{fig:AandB}).
\end{definition}

\begin{figure}[htbp]
\centering
\begin{tikzpicture}
\coordinate (P) at (0,0);
\fill [red!20!white] (P) --++(1.9,0.6) arc (17.5:342.5:2cm) --cycle;
\draw (P) ++ (0,1.5) node[right] {$A$};
%領域A塗りつぶし

\fill [blue!20!white] (P) --++ (1.9,-0.6) arc (-17.5:17.5:2cm) --cycle;
\draw (P) ++(1.5,0) node {$B$};
%領域B塗りつぶし

%\fill [green!20!white] (P) --++ (0.95,0.3) arc (17.5:342.5:1cm) --cycle;
\fill [green!20!white] (P) circle (1);
%領域V_i塗りつぶし

\fill [black] (P) circle (0.06);
\draw (P) node[below] {$p_i$};
%p_i

\draw (P) circle (1);
\draw (P)++(-0.5,0) node[above]{$V_i$};
%V_i
\draw (P) circle (2);
%U_i

\draw (P) --++ (1.9,0.6) --++ (4,0);
\draw (P) --++ (1.9,-0.6) --++ (4,0);
\draw [densely dashed](P) ++ (2,0) --++(3.9,0);
\draw (P) ++(1.9,0) ++(4,0) node[right]{$s_i$};
%s_iの近傍

\draw (P) ++ (4.5,-0.6) [out=0,in=0] to ++(0,1.2);
\draw (P) ++ (4.82,0.2) [out=30,in=240] to ++(0.4,1) node [right] {$\nu (s_i)$};
%\nu (s_i)

%\draw [<->, <-> = stealth ](P) ++ (2.8,0) --++ (0,0.6);
%\draw (P) ++ (2.8,0) ++ (0,0.3) node[right] {$\delta$};

%\draw (P) --++(1.96,0.62) arc (17.5:342.5:2.05cm) --cycle;
\draw (P) ++ (1.5,1.3) [out=30,in=-120] to ++(0.6,0.3) node[above]{$c_i$} ;
%c_i

\draw (P) ++ (-0.86,-0.5) --++(-0.86,-0.5) node[left]{$a$};
\fill (P) ++ (-0.86,-0.5) ++(-0.86,-0.5) circle(0.04);
\draw (P) ++ (-0.86,-0.5) ++(-0.43,-0.25) to[out=120,in=-90] ++(-0.1,0.3) node[above]{$\ell_{a}$};
%\ell_aとa

\fill (P) ++ (2,0) circle (0.06);
\draw (P) ++ (2,-0.3) node[right]{$p'_{i}$};
%p'_{i}

\end{tikzpicture}
\caption{The regions $A$ and $B$}
\label{fig:AandB}
\end{figure}

We now examine the local behavior around each critical point.
Since each critical point is an indefinite fold, by using appropriate local coordinates $(X, Y,v)$ of $\partial M_{0}$ and $(u,v)$ of $U_i \setminus V_i$, $\pi_0 |_{\partial M_0}$ is described as follows:
\[
\pi_0 |_{\partial M_0}: (X,Y, v) \mapsto (X^2-Y^2, v).
\]
Now, the critical value set $\{u=0\}$ intersects $\ell_a$ transversally. We take the coordinate $u$ so that it is positive in the direction of increasing radius.
On the other hand, if we denote by $(s,t)$ the complex line coordinates of the fiber over each point $(u,v)$, then the projection can be written as
$
\pi(s,t,u,v)=(u,v).
$
Since the local coordinates of $\partial M_0$ were expressed using $(X,Y,v)$, the embedding $F: \partial M_0 \rightarrow \bC^2 \cong \bR^4$ is expressed as:
\[
F: (X,Y,v) \mapsto (f(X,Y,v), g(X,Y,v), X^2-Y^2, v),
\]
where $f$ and $g$ are smooth functions. 
Let us fix $v=0$ and define a smooth map $F_0: \bR^2 \rightarrow \bR^2$ as $F_0: (X,Y) \mapsto (s=f(X,Y,0), t=g(X,Y,0))$.
Since $F$ induces the embedding $(X,Y) \mapsto (f(X,Y,0), g(X,Y,0), X^2 - Y^2)$ and $(X,Y)=(0,0)$ is a critial point of $X^2- Y^2$, $dF_0$ is a rank $2$-matrix at $(0,0)$. Hence, by the inverse function theorem, $F_0$ is locally a diffeomorphism around a neighbourhood $U_0$ of $(X, Y) = (0,0)$.
Similarly, it follows that  $(X,Y) \mapsto (f(X,Y,v_0), g(X,Y,v_0))$ is a diffeormorphism on $U_0$ for fixed small $v_0$.
Therefore, there exists a small neighbourhood $U_1 \subset \bC^2$ of the critical point such that the map 
\[
\varphi: U_1 \rightarrow U_0 \times U_2 ; (s,t,u,v) \mapsto (X,Y,u,v)
\] 
is a diffeomorphism, where $s = f(X,Y,v)$ and $t=g(X,Y,v)$ and $U_2 \subset U_i \setminus V_i$ is a neighbourhood of $(u,v)=(0,0)$. Thus, after composing $\varphi$, the projection $\pi$ is expressed as $\pi (X,Y,u,v)=(u,v)$ and $\partial M_0$ is locally described as 
\[
\{(X,Y,u,v) \mid u= X^2-Y^2\}.
\]
Moreover, since the positive direction of $u$ points to the outward of $\partial U_i$, the interior $M_0$ is locally expressed as 
\[
\{(X,Y,u,v) \mid u \geq X^2- Y^2\}
\]

\begin{proposition}\label{prop:attachingA}
Attaching $\pi_0^{-1} (A)$ to $M_1$ does not change the topology.
\end{proposition}

\proof
Consider the radius function $\mathbf{r}: A \rightarrow \bR$.
If necessary, by an appropriate coordinate change on $A$, we may assume that each component of the critical value set $S \cap A$ of $\pi_0 |_{\partial M_0}$ is quadratically tangent to the level set of $\mathbf{r}$ at exactly one point, from inside on $A$.
The critical point of the composition $\mathbf{r} \circ \pi_{0} |_{\partial M_0}$ over $A$ is the critical point of $\pi_0 |_{\partial M_0}$ whose image is the tangent point with the level set of the radius function. 
Since the tangent is quadratic, it is a non-degenerate critical point. 
We will use $- \mathbf{r}$ as the Morse function.
Since the tangent is from inside, $\mathbf{r}$ and locally described as $-\mathbf{r}
 :(u,v) \mapsto -u+v^2$ around the critical point (Figure \ref{fig:radius}, recall that $u$ points to the outward of $\partial U_i$). 
Thus the composition is described as $(-\mathbf{r} \circ \pi_0)(X,Y,u,v)= - u + v^{2}=-X^2+Y^2+v^2$.
Since the interior $M_0$ is locally described as $\{(X,Y,u,v) \mid u \geq X^2- Y^2\}$, the differential $d(-\mathbf{r} \circ \pi_0)$ evaluated on the outward normal vector (for example, $(0,0,-1,0)\, ^{t}(\partial_X, \partial_Y, \partial_u, \partial_v)$) is positive. Thus, by Morse theory for manifolds with boundary (see Satz 4.1 in \cite{bra}, Lemma 2.19 and Theorem 2.27 in \cite{BNR}, Lemma 3.3 in \cite{jr}, for example), passing through this critical point does not change the topology.
Thus, attaching $\pi^{-1}_0  (A)$ do not change the topology.
\endproof

\begin{figure}[htbp]
\centering
\begin{tikzpicture}
\draw[->] (-2.3,0) --(2,0)node[right]{$u$};
\draw (2.5,0) node[right]{$(\partial U_i)$};
\draw (-2.5,0) node[left]{$(V_i)$};
\draw[thick] (-1,-2) -- (-1,2)node[above]{$v$};
\draw[->] (-1,2);
\draw[domain=-1.4:1.4, smooth, variable=\x]
    plot ({\x*\x-1}, {\x});
\draw[domain=-1.4:1.4, smooth, variable=\x]
    plot ({\x*\x-2}, {\x});
\draw[domain=-1.4:1.4, smooth, variable=\x]
    plot ({\x*\x-0}, {\x});
\draw (-1,-1.5) to[out=180,in=00] ++ (-0.3,-0.3)node[left]{$S$};
\draw (1,-1.35) node[below]{level sets of $-\mathbf{r}$};
\draw[very thick,->] (1.25,1.6) --++ (-1,0);
\end{tikzpicture}
\caption{Local behavior of $-\mathbf{r}$}
\label{fig:radius}
\end{figure}

Since $\pi_0^{-1} (q) \cong D_{n-(m_i -1)}$ for all $q \in V_i$, we have
\[
\pi^{-1}_0 (V_i) \cong D_{n-(m_i -1) } \times D^2 \cong  (D_{n-(m_i -1) } \times [0,1] ) \times [0,1].
\]
The piece $\pi^{-1}_0 (V_i)$ is attached to $\pi^{-1}_0 (V_i \cap A) \cong D_{n-(m_i-1)} \times [0,1]$. Thus by Lemma \ref{lem:handle} (1), attaching $\pi^{-1}_0 (V_i)$ to $\pi^{-1}_0 (A)$ does not change the topology.

Let us denote $M_2 = M_1 \cup \pi^{-1}_0 (A \cup V_i)$, which is diffeomorphic to $M_1$ by the above discussion.

\subsection{The attachment of $\pi^{-1}_0 (B)$}

Next, we will study the attachment of $\pi^{-1}_0 (B)$ to $M_2$. For the attachment of $\pi^{-1}_0 (B)$, we have the following. 

\begin{figure}[htbp]
\centering
\begin{tikzpicture}
\coordinate (P1) at (4.5,0) ;
%ファイバー1
\coordinate (P2) at (0,0) ;
%ファイバー2
\coordinate (P3) at (-4.5,0) ;
%ファイバー3
\coordinate (P4) at (-9,0) ;
%ファイバー4
\coordinate (P5) at (-7.5,-5.5) ;
%底空間

%\fill (P1) circle (0.06);
\fill (P2) circle (0.06);
%\fill (P3) circle (0.06);
%\fill (P4) circle (0.06);
%\fill (P5) circle (0.06);

\draw (P1) ++ (-1.5,2) --++ (3,0);
\draw (P1) ++ (-1.5,1) --++ (3,0);
%%上筒
\draw (P1) ++ (-1.5,1) [out =0,in=-90] to++ (0.3,0.5) [out=90,in=0] to++ (-0.3,0.5);
\draw (P1) ++ (-1.5,1) [out =180,in=-90] to++ (-0.3,0.5) [out=90,in=180] to++ (0.3,0.5);
%%左
\draw (P1) ++ (-1.5,1) ++ (3,0) [out =0,in=-90] to++ (0.3,0.5) [out=90,in=0] to++ (-0.3,0.5);
\draw [densely dashed](P1) ++ (-1.5,1)++ (3,0) [out =180,in=-90] to++ (-0.3,0.5) [out=90,in=180] to++ (0.3,0.5);
%%右
%P1ファイバーの上側

\draw (P1) ++ (-1.5,-1) --++ (3,0);
\draw (P1) ++ (-1.5,-2) --++ (3,0);
%%下筒
\draw (P1) ++ (-1.5,-2) [out =0,in=-90] to++ (0.3,0.5) [out=90,in=0] to++ (-0.3,0.5);
\draw (P1) ++ (-1.5,-2) [out =180,in=-90] to++ (-0.3,0.5) [out=90,in=180] to++ (0.3,0.5);
%%左
\draw (P1) ++ (-1.5,-2) ++ (3,0) [out =0,in=-90] to++ (0.3,0.5) [out=90,in=0] to++ (-0.3,0.5);
\draw [densely dashed](P1) ++ (-1.5,-2)++ (3,0) [out =180,in=-90] to++ (-0.3,0.5) [out=90,in=180] to++ (0.3,0.5);
%%右
%P1ファイバーの下側

\draw (P2) ++ (-1.5,2) --++ (3,0);
\draw (P2) ++ (-1.5,1) [out=0,in=140] to++ (1,-0.5) --++ (0.5,-0.5);
\draw (P2) --++ (0.5,0.5)[out=40,in=180] to ++(1,0.5);
%%上筒
\draw (P2) ++ (-1.5,1) [out =0,in=-90] to ++ (0.3,0.5) [out=90,in=0] to++ (-0.3,0.5);
\draw (P2) ++ (-1.5,1) [out =180,in=-90] to++ (-0.3,0.5) [out=90,in=180] to++ (0.3,0.5);
%%左
\draw (P2) ++ (-1.5,1) ++ (3,0) [out =0,in=-90] to++ (0.3,0.5) [out=90,in=0] to++ (-0.3,0.5);
\draw [densely dashed](P2) ++ (-1.5,1)++ (3,0) [out =180,in=-90] to++ (-0.3,0.5) [out=90,in=180] to++ (0.3,0.5);
%%右
%P2ファイバーの上側

\draw (P2) ++ (-1.5,-1) [out=0,in=220] to++ (1,0.5) --++ (0.5,0.5);
\draw (P2) --++ (0.5,-0.5)[out=-40,in=180] to ++(1,-0.5);
\draw (P2) ++ (-1.5,-2) --++ (3,0);
%%下筒
\draw (P2) ++ (-1.5,-2) [out =0,in=-90] to++ (0.3,0.5) [out=90,in=0] to++ (-0.3,0.5);
\draw (P2) ++ (-1.5,-2) [out =180,in=-90] to++ (-0.3,0.5) [out=90,in=180] to++ (0.3,0.5);
%%左
\draw (P2) ++ (-1.5,-2) ++ (3,0) [out =0,in=-90] to++ (0.3,0.5) [out=90,in=0] to++ (-0.3,0.5);
\draw [densely dashed](P2) ++ (-1.5,-2)++ (3,0) [out =180,in=-90] to++ (-0.3,0.5) [out=90,in=180] to++ (0.3,0.5);
%%右
%P2ファイバーの下側

\draw (P3) ++ (-1.5,2) --++ (3,0);
%\draw (P3) ++ (-1.5,1) --++ (3,0);
\draw (P3) ++ (-1.5,1) [out=0,in=90] to++(1,-1);
\draw (P3) ++ (1.5,1) [out=180,in=90] to++(-1,-1);
%%上筒
\draw (P3) ++ (-1.5,1) [out =0,in=-90] to++ (0.3,0.5) [out=90,in=0] to++ (-0.3,0.5);
\draw (P3) ++ (-1.5,1) [out =180,in=-90] to++ (-0.3,0.5) [out=90,in=180] to++ (0.3,0.5);
%%左
\draw (P3) ++ (-1.5,1) ++ (3,0) [out =0,in=-90] to++ (0.3,0.5) [out=90,in=0] to++ (-0.3,0.5);
\draw [densely dashed](P3) ++ (-1.5,1)++ (3,0) [out =180,in=-90] to++ (-0.3,0.5) [out=90,in=180] to++ (0.3,0.5);
%%右
%P3ファイバーの上側

%\draw (P3) ++ (-1.5,-1) --++ (3,0);
\draw (P3) ++ (-1.5,-2) --++ (3,0);
\draw (P3) ++ (-1.5,-1) [out=0,in=-90] to++(1,1);
\draw (P3) ++ (1.5,-1) [out=180,in=-90] to++(-1,1);
%%下筒
\draw (P3) ++ (-1.5,-2) [out =0,in=-90] to++ (0.3,0.5) [out=90,in=0] to++ (-0.3,0.5);
\draw (P3) ++ (-1.5,-2) [out =180,in=-90] to++ (-0.3,0.5) [out=90,in=180] to++ (0.3,0.5);
%%左
\draw (P3) ++ (-1.5,-2) ++ (3,0) [out =0,in=-90] to++ (0.3,0.5) [out=90,in=0] to++ (-0.3,0.5);
\draw [densely dashed](P3) ++ (-1.5,-2)++ (3,0) [out =180,in=-90] to++ (-0.3,0.5) [out=90,in=180] to++ (0.3,0.5);
%%右
%P3ファイバーの下側

\draw (P4) ++ (-1.5,2) --++ (3,0);
%\draw (P4) ++ (-1.5,1) --++ (3,0);
%%上筒
%\draw (P4) ++ (-1.5,1) [out =0,in=-90] to++ (0.3,0.5) [out=90,in=0] to++ (-0.3,0.5);
%\draw (P4) ++ (-1.5,1) [out =180,in=-90] to++ (-0.3,0.5) [out=90,in=180] to++ (0.3,0.5);
%%左
%\draw (P4) ++ (-1.5,1) ++ (3,0) [out =0,in=-90] to++ (0.3,0.5) [out=90,in=0] to++ (-0.3,0.5);
%\draw [densely dashed](P4) ++ (-1.5,1)++ (3,0) [out =180,in=-90] to++ (-0.3,0.5) [out=90,in=180] to++ (0.3,0.5);
%%右
\draw (P4) ++ (-1.5,-2)  [out =0,in=-90] to++ (0.5,2) [out=90,in=0] to++ (-0.5,2);
\draw (P4) ++ (-1.5,-2)  [out =180,in=-90] to++ (-0.5,2) [out=90,in=180] to++ (0.5,2);
%P4ファイバーの左側

%\draw (P4) ++ (-1.5,-1) --++ (3,0);
\draw (P4) ++ (-1.5,-2) --++ (3,0);
%%下筒
%\draw (P4) ++ (-1.5,-2) [out =0,in=-90] to++ (0.3,0.5) [out=90,in=0] to++ (-0.3,0.5);
%\draw (P4) ++ (-1.5,-2) [out =180,in=-90] to++ (-0.3,0.5) [out=90,in=180] to++ (0.3,0.5);
%%左
%\draw (P4) ++ (-1.5,-2) ++ (3,0) [out =0,in=-90] to++ (0.3,0.5) [out=90,in=0] to++ (-0.3,0.5);
%\draw [densely dashed](P4) ++ (-1.5,-2)++ (3,0) [out =180,in=-90] to++ (-0.3,0.5) [out=90,in=180] to++ (0.3,0.5);
%%右
\draw (P4) ++ (-1.5,-2) ++(3,0)  [out =0,in=-90] to++ (0.5,2) [out=90,in=0] to++ (-0.5,2);
\draw [densely dashed](P4) ++ (-1.5,-2) ++(3,0) [out =180,in=-90] to++ (-0.5,2) [out=90,in=180] to++ (0.5,2);
%P4ファイバーの右側

\draw (P5) ++ (0,1.5) --++(10,0);
\draw (P5) ++ (0,0) --++(10,0);
%\draw (P5) ++ (2.5,0) --++ (0,1.5);
\draw (P5) ++ (3.5,0) node[below] {$\alpha_{r_1}$};
%\draw (P5) ++ (5,0) --++ (0,1.5);
\draw (P5) ++ (6,0) node[below] {$\alpha_{r_2}$};
%\draw (P5) ++ (6.5,0) --++ (0,1.5);
\draw (P5) ++ (6.7,0) node[below] {$\alpha_{0}$};
%\draw (P5) ++ (8,0) --++ (0,1.5);
\draw (P5) ++ (8.5,0) node[below] {$\alpha_{r_3}$};
\draw [very thick](P5) ++ (5.5,0) to++ (0,1.5);
%\draw [very thick](P5) ++ (4.5,1.5) [out=-30, in=90]to++ (2,-0.75);
\draw (P5) ++ (5.5,0.4) [out=180,in=90]to++ (-0.5,-0.5)node[below]{$S$};
\draw[thick] plot[domain=-0.75:0.75] ({2*pow(\x,2)-5},\x-4.75);
\draw[thick] plot[domain=-0.75:0.75] ({2*pow(\x,2)-2.6},\x-4.75);
\draw[thick] plot[domain=-0.75:0.75] ({2*pow(\x,2)-2},\x-4.75);
\draw[thick] plot[domain=-0.75:0.75] ({2*pow(\x,2)-0},\x-4.75);
%底空間

\draw [->,->=stealth](P4) ++ (1.5,-2.5) --++(3,-1); 
\draw [->,->=stealth](P3) ++ (1.3,-2.5) --++(1.5,-1);
\draw [->,->=stealth](P2) ++ (-0.3,-2.5) --++(-0.5,-1);  
\draw [->,->=stealth](P1) ++ (-1.5,-2.5) --++(-1.8,-1); 
%矢印たち

\draw (P4) ++ (0,2) node[above]{$(\pi_0 |_{\partial M_{0}})^{-1} (\alpha_{r_1})$};
\draw (P3) ++ (0,2) node[above]{$(\pi_0 |_{\partial M_{0}})^{-1} (\alpha_{r_2})$};
\draw (P2) ++ (0,2) node[above]{$(\pi_0 |_{\partial M_{0}})^{-1} (\alpha_{0})$};
\draw (P1) ++ (0,2) node[above]{$(\pi_0 |_{\partial M_{0}})^{-1} (\alpha_{r_3})$};
%ファイバー名前

\draw (P2) ++ (0.2,0)node[right] {$q$};
%q名前

\draw[densely dashed,blue] (P3) ++ (-0.5,0) [out=90,in=90] to++(1,0);
\draw[blue] (P3) ++ (-0.5,0) [out=-90,in=-90] to++(1,0);
%消滅サイクル

\end{tikzpicture}
\caption{The fibers of $\pi_{0} |_{\partial M}$ over $B$}
\label{fig:fibers}
\end{figure}
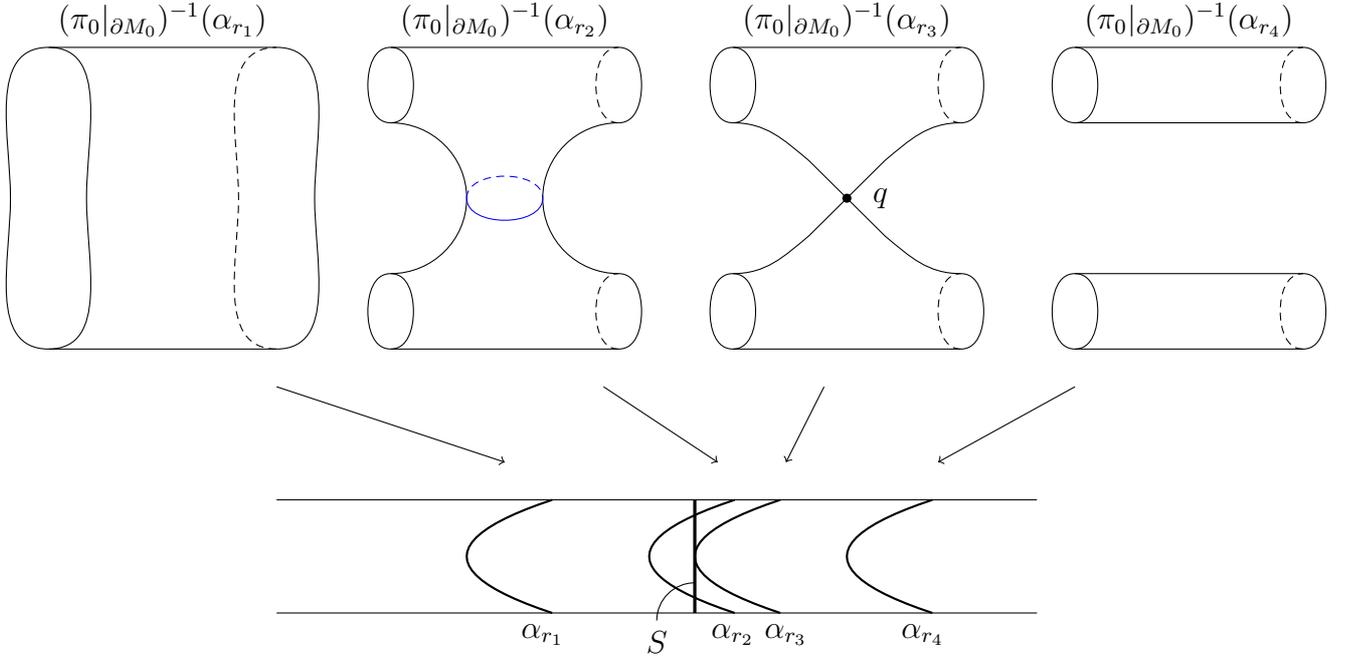

\begin{theorem}\label{thm:2-handle}
Attaching $\pi^{-1}_0 (B)$ to $M_2$ is equivalent to attaching $(m_i-1)$ $2$-handles. 
\end{theorem}
From now on, we will give the proof of this theorem.
The intersection of $S$ and $B$ is a disjoint union of $(m_i -1)$ segments (see the arguments just before Definition \ref{def:a&b}).
The topological change occurs when the fiber of $\pi_0$ passes a critical point. 
For simplicity, we restrict to the case $m_i=2$.
Again, we will use local coordinates $(X,Y,u,v)$, which were used in the previous discussion.
Recall that $\pi_0 |_{\partial M_0}$ is locally described as $\pi_0 |_{\partial M_0}: (X,Y,v) \mapsto (X^2-Y^2, v)$, and the critical value set $S = \{(u,v) \mid u=0\}$.
Let $\alpha_{r}$ be the intersection of the curve $\{(u,v) \mid u=v^2 + r\}$ and $B$.
Then there exist $\delta_1 <0$ and $\delta_2 >0$ such that $B = \bigcup_{\delta_1 \leq r \leq \delta_2}\alpha_r$ and $\pi^{-1}_0(B) = \bigcup_{\delta_1 \leq r \leq \delta_2}\pi^{-1}_0 (\alpha_r)$.
We will look the attachment of $\pi^{-1}_0 (\alpha_r)$ for each $\delta_1 \leq r \leq \delta_2$ (Figure \ref{fig:fibers}).

Let $\eta$ be a sufficiently small positive number. At first, attaching $\bigcup_{\delta_1 \leq r \leq -\eta} \pi^{-1}_0 (\alpha_r)$ to $M_2$ does not change the diffeomorphism type by Lemma \ref{lem:handle} (2).

Next, we consider the part $-\eta \leq r \leq \eta$.
Let $q$ be the intersection point of the critical point set and $\pi^{-1}_0(\alpha_0)$. Around the critical point in $\bigcup_{-\eta \leq r \leq \eta}\pi^{-1}_0 (\alpha_r)$, $\partial M_0$ is described as $\{(X,Y,u,v) \mid u =X^2 - Y^2\}$. 
Thus, by using $r = u-v^2$ instead of $u$, we have the following description:
\[
\{(X,Y,r,v) \mid X^2-Y^2=v^2+r, \, -\eta \leq r \leq \eta\}.
\]

\begin{figure}[htbp]
\centering
\begin{tikzpicture}

%%%%%接着前！%%%%%%%%

\coordinate (P1) at (1.5,0) ;
%ファイバー1
\coordinate (P2) at (-3,0) ;
%ファイバー2
\coordinate (P3) at (-7.5,0) ;
%ファイバー3
\coordinate (P4) at (-9,0) ;
%ファイバー4
\coordinate (P5) at (-7.5,-5.5) ;
%底空間

%\fill (P1) circle (0.06);
%\fill (P2) circle (0.06);
%\fill (P3) circle (0.06);
%\fill (P4) circle (0.06);
%\fill (P5) circle (0.06);

\draw (P1) ++ (-1.5,2) --++ (3,0);
%\draw (P1) ++ (-1.5,1) --++ (3,0);
\draw (P1) ++ (-1.5,1) [out=0,in=90] to++(1,-1);
\draw (P1) ++ (1.5,1) [out=180,in=90] to++(-1,-1);
%%上筒
\draw (P1) ++ (-1.5,1) [out =0,in=-90] to++ (0.3,0.5) [out=90,in=0] to++ (-0.3,0.5);
\draw (P1) ++ (-1.5,1) [out =180,in=-90] to++ (-0.3,0.5) [out=90,in=180] to++ (0.3,0.5);
%%左
\draw (P1) ++ (-1.5,1) ++ (3,0) [out =0,in=-90] to++ (0.3,0.5) [out=90,in=0] to++ (-0.3,0.5);
\draw [densely dashed](P1) ++ (-1.5,1)++ (3,0) [out =180,in=-90] to++ (-0.3,0.5) [out=90,in=180] to++ (0.3,0.5);
%%右
%P1ファイバーの上側

%\draw (P1) ++ (-1.5,-1) --++ (3,0);
\draw (P1) ++ (-1.5,-2) --++ (3,0);
\draw (P1) ++ (-1.5,-1) [out=0,in=-90] to++(1,1);
\draw (P1) ++ (1.5,-1) [out=180,in=-90] to++(-1,1);
%%下筒
\draw (P1) ++ (-1.5,-2) [out =0,in=-90] to++ (0.3,0.5) [out=90,in=0] to++ (-0.3,0.5);
\draw (P1) ++ (-1.5,-2) [out =180,in=-90] to++ (-0.3,0.5) [out=90,in=180] to++ (0.3,0.5);
%%左
\draw (P1) ++ (-1.5,-2) ++ (3,0) [out =0,in=-90] to++ (0.3,0.5) [out=90,in=0] to++ (-0.3,0.5);
\draw [densely dashed](P1) ++ (-1.5,-2)++ (3,0) [out =180,in=-90] to++ (-0.3,0.5) [out=90,in=180] to++ (0.3,0.5);
%%右
%P1ファイバーの下側

\draw (P2) ++ (-1.5,2) --++ (3,0);
%\draw (P2) ++ (-1.5,1) --++ (3,0);
\draw (P2) ++ (-1.5,1) [out=0,in=90] to++(1,-1);
\draw (P2) ++ (1.5,1) [out=180,in=90] to++(-1,-1);
%%上筒
\draw (P2) ++ (-1.5,1) [out =0,in=-90] to++ (0.3,0.5) [out=90,in=0] to++ (-0.3,0.5);
\draw (P2) ++ (-1.5,1) [out =180,in=-90] to++ (-0.3,0.5) [out=90,in=180] to++ (0.3,0.5);
%%左
\draw (P2) ++ (-1.5,1) ++ (3,0) [out =0,in=-90] to++ (0.3,0.5) [out=90,in=0] to++ (-0.3,0.5);
\draw [densely dashed](P2) ++ (-1.5,1)++ (3,0) [out =180,in=-90] to++ (-0.3,0.5) [out=90,in=180] to++ (0.3,0.5);
%%右
%P2ファイバーの上側

%\draw (P2) ++ (-1.5,-1) --++ (3,0);
\draw (P2) ++ (-1.5,-2) --++ (3,0);
\draw (P2) ++ (-1.5,-1) [out=0,in=-90] to++(1,1);
\draw (P2) ++ (1.5,-1) [out=180,in=-90] to++(-1,1);
%%下筒
\draw (P2) ++ (-1.5,-2) [out =0,in=-90] to++ (0.3,0.5) [out=90,in=0] to++ (-0.3,0.5);
\draw (P2) ++ (-1.5,-2) [out =180,in=-90] to++ (-0.3,0.5) [out=90,in=180] to++ (0.3,0.5);
%%左
\draw (P2) ++ (-1.5,-2) ++ (3,0) [out =0,in=-90] to++ (0.3,0.5) [out=90,in=0] to++ (-0.3,0.5);
\draw [densely dashed](P2) ++ (-1.5,-2)++ (3,0) [out =180,in=-90] to++ (-0.3,0.5) [out=90,in=180] to++ (0.3,0.5);
%%右
%P2ファイバーの下側

\draw (P3) ++ (-1.5,2) --++ (3,0);
%\draw (P3) ++ (-1.5,1) --++ (3,0);
\draw (P3) ++ (-1.5,1) [out=0,in=90] to++(1,-1);
\draw (P3) ++ (1.5,1) [out=180,in=90] to++(-1,-1);
%%上筒
\draw (P3) ++ (-1.5,1) [out =0,in=-90] to++ (0.3,0.5) [out=90,in=0] to++ (-0.3,0.5);
\draw (P3) ++ (-1.5,1) [out =180,in=-90] to++ (-0.3,0.5) [out=90,in=180] to++ (0.3,0.5);
%%左
\draw (P3) ++ (-1.5,1) ++ (3,0) [out =0,in=-90] to++ (0.3,0.5) [out=90,in=0] to++ (-0.3,0.5);
\draw [densely dashed](P3) ++ (-1.5,1)++ (3,0) [out =180,in=-90] to++ (-0.3,0.5) [out=90,in=180] to++ (0.3,0.5);
%%右
%P3ファイバーの上側

%\draw (P3) ++ (-1.5,-1) --++ (3,0);
\draw (P3) ++ (-1.5,-2) --++ (3,0);
\draw (P3) ++ (-1.5,-1) [out=0,in=-90] to++(1,1);
\draw (P3) ++ (1.5,-1) [out=180,in=-90] to++(-1,1);
%%下筒
\draw (P3) ++ (-1.5,-2) [out =0,in=-90] to++ (0.3,0.5) [out=90,in=0] to++ (-0.3,0.5);
\draw (P3) ++ (-1.5,-2) [out =180,in=-90] to++ (-0.3,0.5) [out=90,in=180] to++ (0.3,0.5);
%%左
\draw (P3) ++ (-1.5,-2) ++ (3,0) [out =0,in=-90] to++ (0.3,0.5) [out=90,in=0] to++ (-0.3,0.5);
\draw [densely dashed](P3) ++ (-1.5,-2)++ (3,0) [out =180,in=-90] to++ (-0.3,0.5) [out=90,in=180] to++ (0.3,0.5);
%%右
%P3ファイバーの下側

\draw (P5) ++ (0,1.5) --++(8,0);
\draw (P5) ++ (0,0) --++(8,0);
%\draw (P5) ++ (2.5,0) --++ (0,1.5);
\draw (P5) ++ (3.9,0) node[below] {$\alpha_{-\eta}$};
%\draw (P5) ++ (5,0) --++ (0,1.5);
\draw (P5) ++ (4.7,0) node[below] {$\alpha_{0}$};
%\draw (P5) ++ (6.5,0) --++ (0,1.5);
\draw (P5) ++ (5.5,0) node[below] {$\alpha_{\eta}$};
%\draw (P5) ++ (8,0) --++ (0,1.5);
%\draw (P5) ++ (8.5,0) node[below] {$\alpha_{r_4}$};
\draw [very thick](P5) ++ (3.5,0) to++ (0,1.5);
%\draw [very thick](P5) ++ (4.5,1.5) [out=-30, in=90]to++ (2,-0.75);
\draw (P5) ++ (3.5,0.4) [out=180,in=90]to++ (-0.5,-0.5)node[below]{$S$};
\draw[thick] plot[domain=-0.75:0.75] ({2*pow(\x,2)-4.8},\x-4.75);
\draw[thick] plot[domain=-0.75:0.75] ({2*pow(\x,2)-4},\x-4.75);
\draw[thick] plot[domain=-0.75:0.75] ({2*pow(\x,2)-3.2},\x-4.75);
%底空間

%\draw [->,->=stealth](P4) ++ (1.5,-2.5) --++(3,-1); 
\draw [->,->=stealth](P3) ++ (1.3,-2.5) --++(1.5,-1);
\draw [->,->=stealth](P2) ++ (0,-2.5) --++(0,-1);  
\draw [->,->=stealth](P1) ++ (-1.5,-2.5) --++(-1.8,-1);
%矢印たち

%\draw (P4) ++ (0,2) node[above]{$(\pi_0 |_{\partial M})^{-1} (\alpha_{r_1})$};
\draw (P3) ++ (0,2) node[above]{$(\pi_0 |_{\partial M})^{-1} (\alpha_{-\eta})$};
\draw (P2) ++ (0,2) node[above]{$(\pi_0 |_{\partial M})^{-1} (\alpha_{0})$};
\draw (P1) ++ (0,2) node[above]{$(\pi_0 |_{\partial M})^{-1} (\alpha_{\eta})$};
%ファイバー名前

%\draw (P2) ++ (0.2,0)node[right] {$q$};
%q名前

%\draw[densely dashed,blue] (P3) ++ (-0.5,0) [out=90,in=90] to++(1,0);
%\draw[blue] (P3) ++ (-0.5,0) [out=-90,in=-90] to++(1,0);
%消滅サイクル

\draw (P2) ++ (0,-7) node{Before attaching $N_2$};
%After attaching $N_2$

\draw[->,line width=5pt] (P2) ++ (0,-7.7)--++(0,-1);

%%%%%%%ここから接着後！%%%%%%%%%%%

\coordinate (Q1) at (1.5,-12) ;
%ファイバー1
\coordinate (Q2) at (-3,-12) ;
%ファイバー2
\coordinate (Q3) at (-7.5,-12) ;
%ファイバー3
\coordinate (Q4) at (-9,-12) ;
%ファイバー4
\coordinate (Q5) at (-7.5,-17.5) ;
%底空間

%\fill (Q1) circle (0.06);
%\fill (Q2) circle (0.06);
%\fill (Q3) circle (0.06);
%\fill (P4) circle (0.06);
%\fill (Q5) circle (0.06);

\draw (Q1) ++ (-1.5,2) --++ (3,0);
\draw (Q1) ++ (-1.5,1) --++ (3,0);
%%上筒
\draw (Q1) ++ (-1.5,1) [out =0,in=-90] to++ (0.3,0.5) [out=90,in=0] to++ (-0.3,0.5);
\draw (Q1) ++ (-1.5,1) [out =180,in=-90] to++ (-0.3,0.5) [out=90,in=180] to++ (0.3,0.5);
%%左
\draw (Q1) ++ (-1.5,1) ++ (3,0) [out =0,in=-90] to++ (0.3,0.5) [out=90,in=0] to++ (-0.3,0.5);
\draw [densely dashed](Q1) ++ (-1.5,1)++ (3,0) [out =180,in=-90] to++ (-0.3,0.5) [out=90,in=180] to++ (0.3,0.5);
%%右
%Q1ファイバーの上側

\draw (Q1) ++ (-1.5,-1) --++ (3,0);
\draw (Q1) ++ (-1.5,-2) --++ (3,0);
%%下筒
\draw (Q1) ++ (-1.5,-2) [out =0,in=-90] to++ (0.3,0.5) [out=90,in=0] to++ (-0.3,0.5);
\draw (Q1) ++ (-1.5,-2) [out =180,in=-90] to++ (-0.3,0.5) [out=90,in=180] to++ (0.3,0.5);
%%左
\draw (Q1) ++ (-1.5,-2) ++ (3,0) [out =0,in=-90] to++ (0.3,0.5) [out=90,in=0] to++ (-0.3,0.5);
\draw [densely dashed](Q1) ++ (-1.5,-2)++ (3,0) [out =180,in=-90] to++ (-0.3,0.5) [out=90,in=180] to++ (0.3,0.5);
%%右
%Q1ファイバーの下側

\draw[white,thick] (Q1) ++ (-0.79,1) --++ (1.58,0);
\fill[white!80!blue] (Q1) ++(-0.8,1) to [out=-20,in=200] ++(1.6,0) to[out=-120,in=120] ++(0,-2)to [out=200,in=-20] ++(-1.6,0) to[out=60,in=-60]++(0,2);
\fill[white!80!blue](Q1)++(0.8,1)arc[start angle=0,end angle=180,x radius=0.8,y radius=0.2] to[out=-20,in=200] ++(1.6,0);
\draw (Q1) ++(-0.8,1) to [out=-20,in=200] ++(1.6,0) to[out=-120,in=120] ++(0,-2)to [out=200,in=-20] ++(-1.6,0) to[out=60,in=-60]++(0,2);

\draw[densely dashed](Q1)++(0.8,1)arc[start angle=0,end angle=180,x radius=0.8,y radius=0.2];
\draw[densely dashed](Q1)++(0.8,-1)arc[start angle=0,end angle=180,x radius=0.8,y radius=0.2];
%Q1ファイバー内のN_2

\fill[white!80!blue] (Q2) ++ (1,-0.88) to[out=120,in=-120]++(0,1.76) to[out=200,in=-20] ++(-2,0) to[out=-60,in=60]++(0,-1.76) to[out=-20,in=200] ++(2,0);
\draw (Q2) ++ (1,-0.88) to[out=120,in=-120]++(0,1.76) to[out=200,in=-20] ++(-2,0) to[out=-60,in=60]++(0,-1.76) to[out=-20,in=200] ++(2,0);
%\draw(Q2)++(1,0.88)arc[start angle=0,end angle=-180,x radius=1,y radius=0.2];
%\draw(Q2)++(1,-0.88)arc[start angle=0,end angle=-180,x radius=1,y radius=0.2];
\draw[densely dashed](Q2)++(1,0.88)arc[start angle=0,end angle=180,x radius=1,y radius=0.2];
\draw[densely dashed](Q2)++(1,-0.88)arc[start angle=0,end angle=180,x radius=1,y radius=0.2];
%Q2ファイバー内のN_2

\draw (Q2) ++ (-1.5,2) --++ (3,0);
\draw (Q2) ++ (-1.5,1) [out=0,in=140] to++ (1,-0.5) --++ (0.5,-0.5);
\draw (Q2) --++ (0.5,0.5)[out=40,in=180] to ++(1,0.5);
%%上筒
\draw (Q2) ++ (-1.5,1) [out =0,in=-90] to ++ (0.3,0.5) [out=90,in=0] to++ (-0.3,0.5);
\draw (Q2) ++ (-1.5,1) [out =180,in=-90] to++ (-0.3,0.5) [out=90,in=180] to++ (0.3,0.5);
%%左
\draw (Q2) ++ (-1.5,1) ++ (3,0) [out =0,in=-90] to++ (0.3,0.5) [out=90,in=0] to++ (-0.3,0.5);
\draw [densely dashed](Q2) ++ (-1.5,1)++ (3,0) [out =180,in=-90] to++ (-0.3,0.5) [out=90,in=180] to++ (0.3,0.5);
%%右
%Q2ファイバーの上側

\draw (Q2) ++ (-1.5,-1) [out=0,in=220] to++ (1,0.5) --++ (0.5,0.5);
\draw (Q2) --++ (0.5,-0.5)[out=-40,in=180] to ++(1,-0.5);
\draw (Q2) ++ (-1.5,-2) --++ (3,0);
%%下筒
\draw (Q2) ++ (-1.5,-2) [out =0,in=-90] to++ (0.3,0.5) [out=90,in=0] to++ (-0.3,0.5);
\draw (Q2) ++ (-1.5,-2) [out =180,in=-90] to++ (-0.3,0.5) [out=90,in=180] to++ (0.3,0.5);
%%左
\draw (Q2) ++ (-1.5,-2) ++ (3,0) [out =0,in=-90] to++ (0.3,0.5) [out=90,in=0] to++ (-0.3,0.5);
\draw [densely dashed](Q2) ++ (-1.5,-2)++ (3,0) [out =180,in=-90] to++ (-0.3,0.5) [out=90,in=180] to++ (0.3,0.5);
%%右
%Q2ファイバーの下側

%\fill[white!80!blue] (Q3) ++ (-0.77,-0.7) to[out=-20,in=200] ++(1.54,0) to[out=130,in=-130]++(0,1.4) to[out=200,in=-20]++(-1.54,0)to[out=-50,in=50] ++(0,-1.4);
%\draw (Q3) ++ (-0.77,-0.7) to[out=-20,in=200] ++(1.54,0);
%\draw (Q3) ++ (-0.77,0.7) to[out=-20,in=200] ++(1.54,0);
%\draw[densely dashed] (Q3) ++ (-0.77,0.7) to[out=10,in=-190] ++(1.54,0);
%\draw[densely dashed] (Q3) ++ (-0.77,-0.7) to[out=10,in=-190] ++(1.54,0);

\draw (Q3) ++ (-1.5,2) --++ (3,0);
%\draw (Q3) ++ (-1.5,1) --++ (3,0);
\draw (Q3) ++ (-1.5,1) [out=0,in=90] to++(1,-1);
\draw (Q3) ++ (1.5,1) [out=180,in=90] to++(-1,-1);
%%上筒
\draw (Q3) ++ (-1.5,1) [out =0,in=-90] to++ (0.3,0.5) [out=90,in=0] to++ (-0.3,0.5);
\draw (Q3) ++ (-1.5,1) [out =180,in=-90] to++ (-0.3,0.5) [out=90,in=180] to++ (0.3,0.5);
%%左
\draw (Q3) ++ (-1.5,1) ++ (3,0) [out =0,in=-90] to++ (0.3,0.5) [out=90,in=0] to++ (-0.3,0.5);
\draw [densely dashed](Q3) ++ (-1.5,1)++ (3,0) [out =180,in=-90] to++ (-0.3,0.5) [out=90,in=180] to++ (0.3,0.5);
%%右
%Q3ファイバーの上側

%\draw (Q3) ++ (-1.5,-1) --++ (3,0);
\draw (Q3) ++ (-1.5,-2) --++ (3,0);
\draw (Q3) ++ (-1.5,-1) [out=0,in=-90] to++(1,1);
\draw (Q3) ++ (1.5,-1) [out=180,in=-90] to++(-1,1);
%%下筒
\draw (Q3) ++ (-1.5,-2) [out =0,in=-90] to++ (0.3,0.5) [out=90,in=0] to++ (-0.3,0.5);
\draw (Q3) ++ (-1.5,-2) [out =180,in=-90] to++ (-0.3,0.5) [out=90,in=180] to++ (0.3,0.5);
%%左
\draw (Q3) ++ (-1.5,-2) ++ (3,0) [out =0,in=-90] to++ (0.3,0.5) [out=90,in=0] to++ (-0.3,0.5);
\draw [densely dashed](Q3) ++ (-1.5,-2)++ (3,0) [out =180,in=-90] to++ (-0.3,0.5) [out=90,in=180] to++ (0.3,0.5);
%%右
%Q3ファイバーの下側

\draw (Q5) ++ (0,1.5) --++(8,0);
\draw (Q5) ++ (0,0) --++(8,0);
%\draw (Q5) ++ (2.5,0) --++ (0,1.5);
\draw (Q5) ++ (3.9,0) node[below] {$\alpha_{-\eta}$};
%\draw (Q5) ++ (5,0) --++ (0,1.5);
\draw (Q5) ++ (4.7,0) node[below] {$\alpha_{0}$};
%\draw (Q5) ++ (6.5,0) --++ (0,1.5);
\draw (Q5) ++ (5.5,0) node[below] {$\alpha_{\eta}$};
%\draw (Q5) ++ (8,0) --++ (0,1.5);
%\draw (Q5) ++ (8.5,0) node[below] {$\alpha_{r_4}$};
\draw [very thick](Q5) ++ (3.5,0) to++ (0,1.5);
%\draw [very thick](Q5) ++ (4.5,1.5) [out=-30, in=90]to++ (2,-0.75);
\draw (Q5) ++ (3.5,0.4) [out=180,in=90]to++ (-0.5,-0.5)node[below]{$S$};
\draw[thick] plot[domain=-0.75:0.75] ({2*pow(\x,2)-4.8},\x-16.75);
\draw[thick] plot[domain=-0.75:0.75] ({2*pow(\x,2)-4},\x-16.75);
\draw[thick] plot[domain=-0.75:0.75] ({2*pow(\x,2)-3.2},\x-16.75);
%底空間

%\draw [->,->=stealth](P4) ++ (1.5,-2.5) --++(3,-1); 
\draw [->,->=stealth](Q3) ++ (1.3,-2.5) --++(1.5,-1);
\draw [->,->=stealth](Q2) ++ (0,-2.5) --++(0,-1);  
\draw [->,->=stealth](Q1) ++ (-1.5,-2.5) --++(-1.8,-1);
%矢印たち

%\draw (P4) ++ (0,2) node[above]{$(\pi_0 |_{\partial M})^{-1} (\alpha_{r_1})$};
\draw (Q3) ++ (0,2) node[above]{$(\pi_0 |_{\partial M})^{-1} (\alpha_{-\eta})$};
\draw (Q2) ++ (0,2) node[above]{$(\pi_0 |_{\partial M})^{-1} (\alpha_{0})$};
\draw (Q1) ++ (0,2) node[above]{$(\pi_0 |_{\partial M})^{-1} (\alpha_{\eta})$};
%ファイバー名前

%\draw (Q2) ++ (0.2,0)node[right] {$q$};
%q名前

\draw[densely dashed,blue] (Q3) ++ (-0.5,0) [out=90,in=90] to++(1,0);
\draw[blue] (Q3) ++ (-0.5,0) [out=-90,in=-90] to++(1,0);
%消滅サイクル

\draw (Q2) ++ (0,-7) node{After attaching $N_2$};
%After attaching $N_2$

%\draw (Q3) ++(-0.5,0)node[left]{attaching circle $\rightarrow$};
\draw (Q1) ++(1.3,-3) node[below]{$N_2$};
\draw[->] (Q1) ++(1.3,-3) to++(-1,2);
\draw[->] (Q1) ++(1.3,-3) to++(-5.5,2.2);
%attaching circle, N_2名前

\end{tikzpicture}
\caption{Attachment of $N_2$}
\label{fig:attaching_N2}
\end{figure}
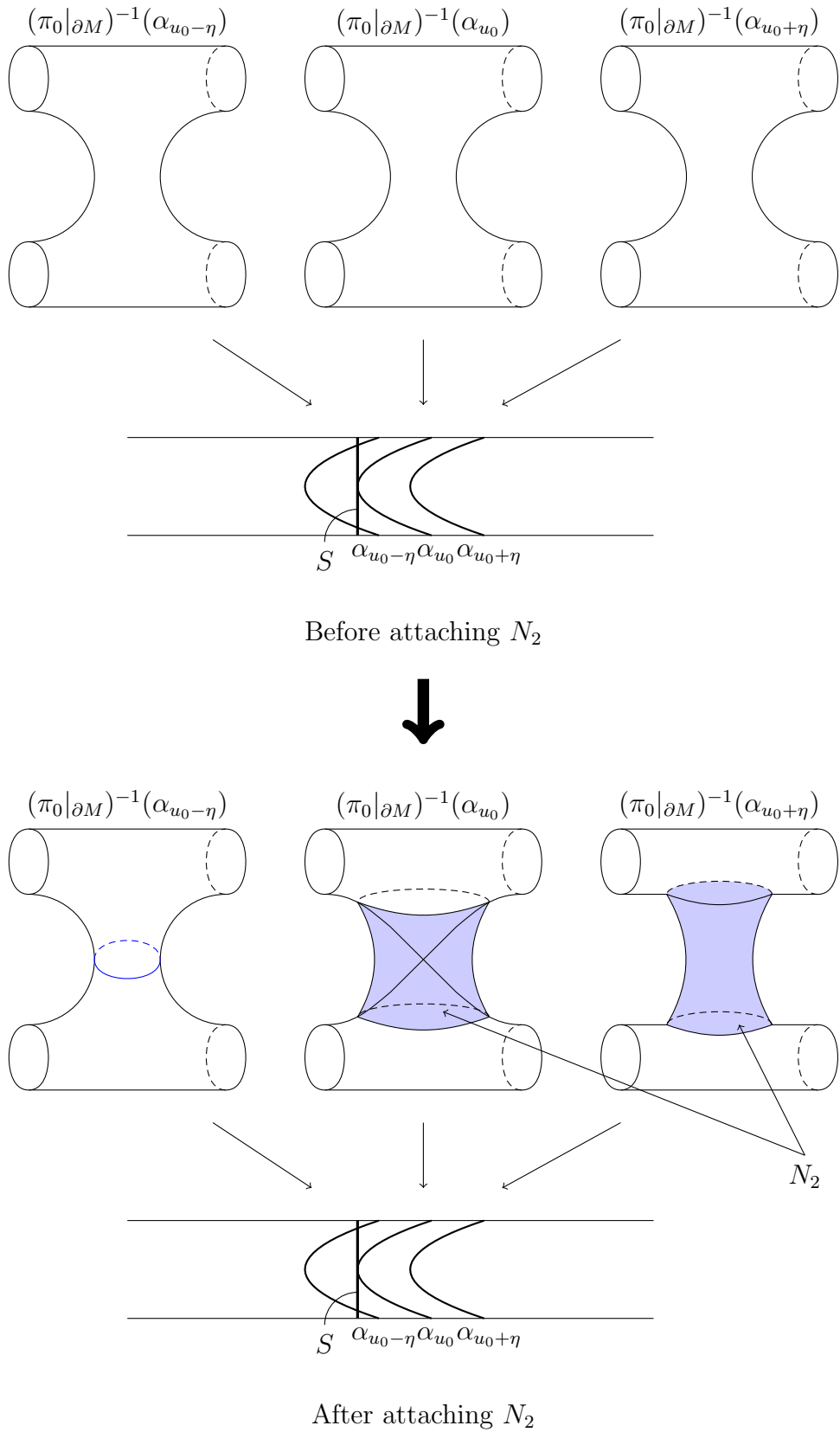

Therefore, around the critical point, the part $N$ which is obtained from $M_2$ by attaching $\bigcup_{-\eta \leq r \leq \eta} \pi_{0}^{-1} (\alpha_{r})$ is described as follows:
\[
N = \{(X,Y,r,v) \mid X^2-Y^2 \leq v^2+r, -\eta \leq r \leq \eta\} \cap B^{4}_{R},
\]
where $B^4_{R}$ is a $4$-ball with radius $R \, (> \eta)$ centered at the origin.
We decompose $N$ into the following two pieces $N_1$ and $N_2$ (see Figure \ref{fig:attaching_N2}):
\begin{eqnarray*}
N_1 &=& \{(X,Y,r,v) \mid X^2 - Y^2 - v^2 \leq -\eta , \, -\eta \leq r \leq \eta \} \cap B^{4}_{R}, \\
N_2 &=& \{(X,Y,r,v) \mid -\eta \leq X^2 - Y^2 - v^2 \leq r \leq \eta \} \cap B^{4}_{R} .
\end{eqnarray*}

The piece $N_1$ is first attached, and next $N_2$ is attached.

\begin{proposition}
The attachment of $N_1$ does not change the topology.
\end{proposition}
\proof
Define a projection $p: N_1 \rightarrow N_1 \cap \partial B^4_R$ as 
\[
p(X,Y,r,v) = \left ( \frac{R^2-r^2}{X^2+Y^2+v^2} X,\frac{R^2-r^2}{X^2+Y^2+v^2} Y, r, \frac{R^2-r^2}{X^2+Y^2+v^2} v \right ).
\]
This is a deformation retraction and a trivial $D^1$-bundle on the interior (the fiber is given by the internally dividing points between $(X,Y,r,v)$ and $p(X,Y,r,v)$). Therefore, by smoothing corners, we have $N_1 \cong (N_1 \cap \partial B^4_R) \times D^1$. The part $N_1$ is attached to $M_2$ along $N_1 \cap \partial B^4_R$. Thus, by Lemma \ref{lem:handle} (1), the attachment of $N_1$ does not change the topology. 
\endproof

On attaching $N_2$, we have the following proposition. 
\begin{proposition}
$N_2$ is diffeomorphic to a $4$-ball, and the attachment of $N_2$ to $M_2 \cup N_1$ is a $2$-handle attachment. 
\end{proposition}

\proof 
Define the region $N'_2$ in $\bR^3$ as follows:
\[
N'_2 = \{(X,Y,v) \mid -\eta \leq X^2 -Y^2 -v^2 \leq \eta\} \cap B^{3}_{R} .
\]
Then, $N'_2$ is diffeomorphic to a $3$-ball $D^3$ after smoothing corners ($B^3_{R}$ is a $3$-ball with radius $R$ centered at the origin). Let $p : N_2 \rightarrow N'_2$ be the projection $p(X,Y,r,v) = (X,Y,v)$. Then, it is a surjective map and each fiber $\{r \mid X^2 -Y^2 - v^2 \leq r \leq \eta \}$ is diffeomorphic to $D^1$, over the interior of $N'_2$. Therefore, by smoothing corners, we obtain that $N_2 \cong N'_2 \times D^1 \cong D^3 \times D^1 \cong D^4$.

The pieces $N_1$ and $N_2$ are attached along $N_{3}$ which is defined as
\[
N_{3} = \{(X,Y,r,v) \mid X^2-Y^2 -v^2 = -\eta, \, -\eta \leq r \leq \eta\} \cap B_{R}^{4}.
\]
This attaching area $N_{3}$ is diffeomorphic to $S^1 \times D^2$ and its core is unknotted in $\partial N_{2}$. 
Also, $N_2$ does not intersect with $M_2$.
Therefore, attaching $N_2$ to $M_2 \cup N_1$ is equivalent to attaching a $2$-handle.
\endproof

%\vspace{5cm}

For the remaining part, attaching $\bigcup_{ \eta \leq r \leq \delta_{2}} \pi^{-1}_{0} (\alpha_{r})$ does not change the diffeomorphism type for the same reason as in the part $\delta_{1} \leq r \leq - \eta$.
Therefore, attaching $\pi_{0}^{-1} (B)$ is equivalent to attaching a $2$-handle (see Figure \ref{fig:attaching_N2} for the attaching of $N_2$). 
%When $m_{i} \geq 3$, we can see that attaching $\pi_{0}^{-1} (B)$ is attaching $(m_{i}-1)$ $2$-handles by repeating this procedure.

Let us consider the case where $m_{i} \geq 3$. In this case, the critical value set $S$ is the disjoint union of $(m_{i}-1)$ segments.
The above discussion tells us that a $2$-handle attaching appears each time $\alpha_{r}$ crosses the critical value segment.
Therefore, attaching $\pi^{-1}_{0} (B)$ is equivalent to attaching $(m_{i}-1)$ $2$-handles and this proves Theorem \ref{thm:2-handle}.

Finally, only the part $\nu(\Gamma) \setminus U_i$ is left. 
Since $\nu (\Gamma) \setminus U_{i}$ is diffeomorphic to $D^2$, $\pi^{-1}_{0} (\nu(\Gamma) \setminus U_{i})$ is diffeomorphic to $\natural_{n} (S^1 \times D^3) \cong \natural_{n} (S^1 \times D^2) \times [0,1]$,
and $\pi^{-1}_{0} (\partial (\nu (\Gamma) \setminus U_{i}) ) \cong \natural_{n} (S^1 \times D^1) \times [0,1] \cong \natural_{n} S^1 \times D^2$. 
Thus,  by Lemma \ref{lem:handle} (2), attaching $\pi^{-1} (\nu (\Gamma) \setminus U_{i} )$ does not change the topology.
This completes the construction of a handle decomposition of $M_0$.

\begin{remark}
Similar arguments as Proposition \ref{prop:attachingA} using Morse theory for manifolds with boundary can show that the attachment of $\pi^{-1}_0 (B)$ is the attachment of $(m_i -1)$ $2$-handles; however, we need a detailed analysis of the attaching circle to describe Kirby diagrams. So, we took this specific method.
\end{remark}

\subsection{Remarks on attaching circles}

\begin{remark}\label{rmk:attachingcircle}
The attaching circle of the $2$-handle is written as $\{(X_0, Y, r_{0}, v) \mid Y^2 + v^2 = X_{0}^{2} + \eta \}$ by using fixed $X_0$ and $r_0$ (see the definition of $N_3$).
This can be moved into $\pi_{0}^{-1} (\alpha_{3}) \subset \partial M_{1}$ by the following procedure (see Figure \ref{fig:attachingcircle}) 
and let $C'_{0}$ be the resulting attaching circle.
\begin{enumerate}[(i)]
\item Move the attaching circle $C_{0}$ in the fiber $(\pi_{0} |_{\partial M_{2}})^{-1} (\alpha_{1})$.

\item Corresponding to moving the arc from $\alpha_{1}$ to $\alpha_{2}$ in the base space, move the attaching circle into $(\pi_{0} |_{\partial (M_{1} \cup \pi^{-1}_{0}(A))})^{-1} (\alpha_{2})$. 

\item Move the attaching circle in $\pi^{-1}_{0} (\alpha_{2})$ so that the fiber of the endpoints of $\alpha_{2}$ link to the attaching circle and except for the endpoints of $\alpha_{2}$, the attaching circle runs parallel to the outside cylinder of $\pi^{-1}_{0} (\alpha_{2})$. Observe that $\pi^{-1}_{0} (\alpha_{2})$ is contained in $\partial (M_1 \cup \pi^{-1}_{0} (A))$.

\item 
Corresponding to moving the arc from $\alpha_{2}$ to $\alpha_{3}$ in the base space, move the attaching circle into $(\pi_{0} |_{\partial M_{1}})^{-1} (\alpha_{3})$. Here, $\beta$ denotes the braid representing the local monodromy. 
Note that $\alpha_{3} = c_{i} \subset \partial U_{i}$. 
The braid representing the local monodromy appears in $\pi^{-1} (c_{i})$ (see Remark \ref{rmk:braid}).
For the precise motion with a specific example, see the following Example \ref{exam:specific}.

\end{enumerate}

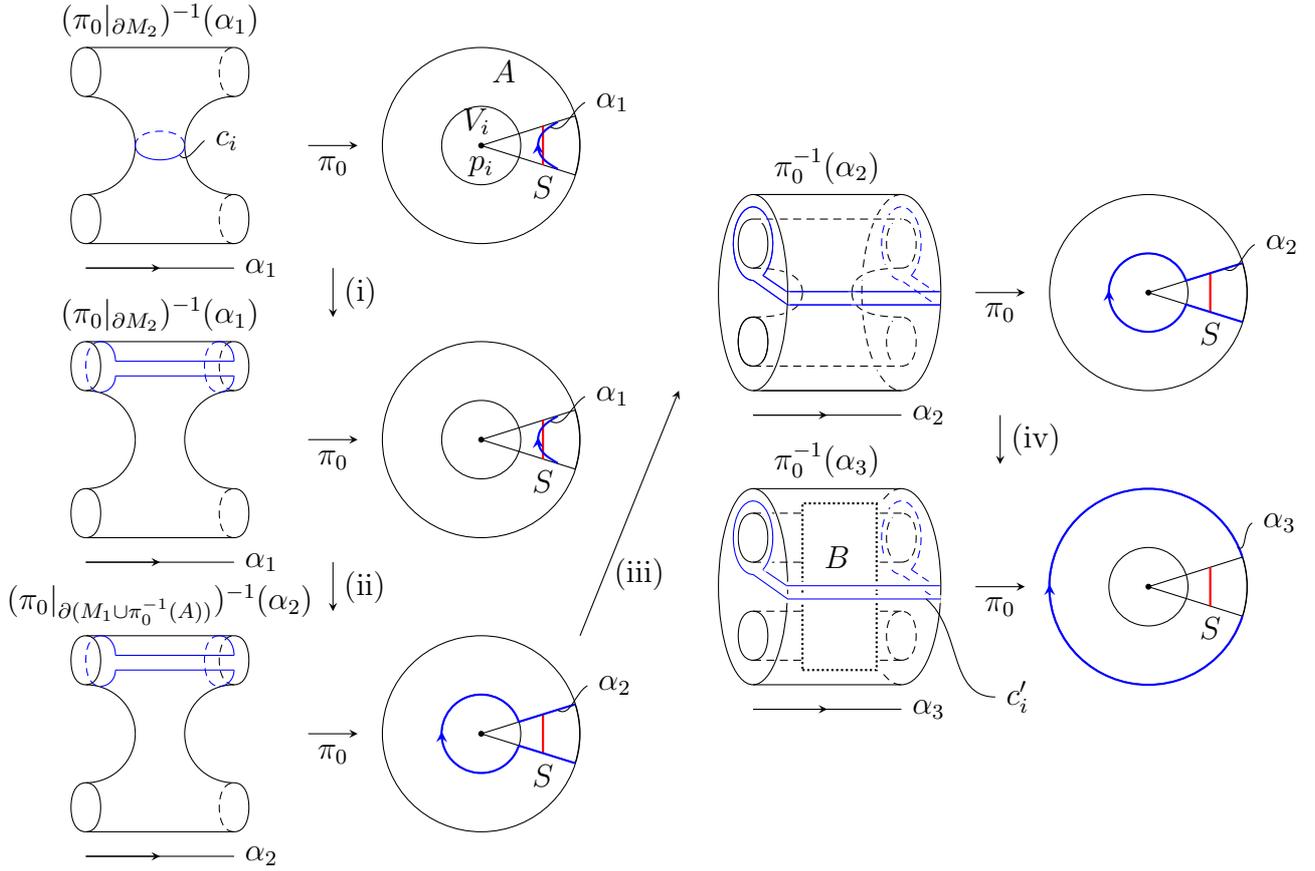
\begin{figure}
\centering
\begin{tikzpicture}[scale =0.65]
\coordinate (P1) at (-3.5,0);
\coordinate (Q1) at (3,0);
\coordinate (P2) at (-3.5,-6);
\coordinate (Q2) at (3,-6);
\coordinate (P3) at (-3.5,-12);
\coordinate (Q3) at (3,-12);
\coordinate (P4) at (10,-3);
\coordinate (Q4) at (16.5,-3);
\coordinate (P5) at (10,-9);
\coordinate (Q5) at (16.5,-9);
%\fill (P1) circle (0.06);

\draw (P1) ++ (-1.5,2) --++ (3,0);
%\draw (P1) ++ (-1.5,1) --++ (3,0);
\draw (P1) ++ (-1.5,1) [out=0,in=90] to++(1,-1);
\draw (P1) ++ (1.5,1) [out=180,in=90] to++(-1,-1);
%%上筒
\draw (P1) ++ (-1.5,1) [out =0,in=-90] to++ (0.3,0.5) [out=90,in=0] to++ (-0.3,0.5);
\draw (P1) ++ (-1.5,1) [out =180,in=-90] to++ (-0.3,0.5) [out=90,in=180] to++ (0.3,0.5);
%%左
\draw (P1) ++ (-1.5,1) ++ (3,0) [out =0,in=-90] to++ (0.3,0.5) [out=90,in=0] to++ (-0.3,0.5);
\draw [densely dashed](P1) ++ (-1.5,1)++ (3,0) [out =180,in=-90] to++ (-0.3,0.5) [out=90,in=180] to++ (0.3,0.5);
%%右
%P1ファイバーの上側

%\draw (P1) ++ (-1.5,-1) --++ (3,0);
\draw (P1) ++ (-1.5,-2) --++ (3,0);
\draw (P1) ++ (-1.5,-1) [out=0,in=-90] to++(1,1);
\draw (P1) ++ (1.5,-1) [out=180,in=-90] to++(-1,1);
%%下筒
\draw (P1) ++ (-1.5,-2) [out =0,in=-90] to++ (0.3,0.5) [out=90,in=0] to++ (-0.3,0.5);
\draw (P1) ++ (-1.5,-2) [out =180,in=-90] to++ (-0.3,0.5) [out=90,in=180] to++ (0.3,0.5);
%%左
\draw (P1) ++ (-1.5,-2) ++ (3,0) [out =0,in=-90] to++ (0.3,0.5) [out=90,in=0] to++ (-0.3,0.5);
\draw [densely dashed](P1) ++ (-1.5,-2)++ (3,0) [out =180,in=-90] to++ (-0.3,0.5) [out=90,in=180] to++ (0.3,0.5);
%%右
%P1ファイバーの下側

\draw[blue] (P1) ++(-0.5,0) to[out=-90,in=-90] ++(1,0);
\draw[densely dashed, blue] (P1) ++(-0.5,0) to[out=90,in=90] ++(1,0);
\draw (P1) ++ (0.4,-0.2) to[out=-50,in=-140] ++(0.5,0.3) node[right]{$C_{0}$};
%P1attaching circle

\draw (P1) ++(0,2)node[above]{$(\pi_{0} |_{\partial M_{2}})^{-1} (\alpha_{1})$};

\draw[->,>=stealth] (P1) ++ (-1.5,-2.5) --++(1.5,0);
\draw (P1) ++ (-1.5,-2.5) --++(3,0) node[right]{$\alpha_{1}$};

\draw[->,>=stealth] (P1) ++ (3,0) --++(1,0);
\draw (P1) ++(3.5,0) node[below]{$\pi_{0}$};

%\fill [red!20!white] (P) --++(1.9,0.6) arc (17.5:342.5:2cm) --cycle;
\draw (Q1) ++ (0,1.5) node[right] {$A$};
%領域A塗りつぶし
\draw (Q1) --++ (1.9,-0.6) arc (-17.5:17.5:2cm) --cycle;
%\draw (P) ++(1.5,0) node[above] {$B$};
%領域B塗りつぶし
%\fill [green!20!white] (P) --++ (0.95,0.3) arc (17.5:342.5:1cm) --cycle;
%領域V_i塗りつぶし

\fill [black] (Q1) circle (0.06);
\draw (Q1) node[below] {$p_i$};
%p_i

\draw (Q1) circle (0.8);
\draw (Q1)++(-0.1,0) node[above]{$V_i$};
%V_i
\draw (Q1) circle (2);
%U_i
\draw[red,thick] (Q1) ++(0.95+0.3,-0.4) --++ (0,0.8);
\draw (Q1) ++(0.95+0.3,-0.4) node[below]{$S$};
%S
\draw[blue,thick] (Q1) ++(0.95+0.6,-0.48) to[out=150,in=-90]++(-0.4,0.48) to[out=90,in=-150] ++(0.4,0.48); 
\draw [->,>=stealth,blue,thick] (Q1) ++ (1.15,0.05);
\draw (Q1) ++(0.95+0.6,0.48) ++(-0.15,-0.1) to[out=-30,in=-150]++(0.7,0.5)node[right]{$\alpha_{1}$};
%\alpha_{1}

%%%%↑ここまで(P1)、(Q1)%%%%

\draw[->,>=stealth] (P1) ++(3.5,-2.5) --++(0,-1);
\draw (P1) ++ (3.5,-3) node[right]{(i)};

%%%%↓ここから(P2)、(Q2)%%%%
\draw (P2) ++ (-1.5,2) --++ (3,0);
%\draw (P2) ++ (-1.5,1) --++ (3,0);
\draw (P2) ++ (-1.5,1) [out=0,in=90] to++(1,-1);
\draw (P2) ++ (1.5,1) [out=180,in=90] to++(-1,-1);
%%上筒
\draw (P2) ++ (-1.5,1) [out =0,in=-90] to++ (0.3,0.5) [out=90,in=0] to++ (-0.3,0.5);
\draw (P2) ++ (-1.5,1) [out =180,in=-90] to++ (-0.3,0.5) [out=90,in=180] to++ (0.3,0.5);
%%左
\draw (P2) ++ (-1.5,1) ++ (3,0) [out =0,in=-90] to++ (0.3,0.5) [out=90,in=0] to++ (-0.3,0.5);
\draw [densely dashed](P2) ++ (-1.5,1)++ (3,0) [out =180,in=-90] to++ (-0.3,0.5) [out=90,in=180] to++ (0.3,0.5);
%%右
%P2ファイバーの上側

%\draw (P2) ++ (-1.5,-1) --++ (3,0);
\draw (P2) ++ (-1.5,-2) --++ (3,0);
\draw (P2) ++ (-1.5,-1) [out=0,in=-90] to++(1,1);
\draw (P2) ++ (1.5,-1) [out=180,in=-90] to++(-1,1);
%%下筒
\draw (P2) ++ (-1.5,-2) [out =0,in=-90] to++ (0.3,0.5) [out=90,in=0] to++ (-0.3,0.5);
\draw (P2) ++ (-1.5,-2) [out =180,in=-90] to++ (-0.3,0.5) [out=90,in=180] to++ (0.3,0.5);
%%左
\draw (P2) ++ (-1.5,-2) ++ (3,0) [out =0,in=-90] to++ (0.3,0.5) [out=90,in=0] to++ (-0.3,0.5);
\draw [densely dashed](P2) ++ (-1.5,-2)++ (3,0) [out =180,in=-90] to++ (-0.3,0.5) [out=90,in=180] to++ (0.3,0.5);
%%右
%P2ファイバーの下側

\draw [blue] (P2) ++ (-0.9,1.6) to[out=90,in=0]++(-0.3,0.4);
\draw [blue] (P2) ++ (-0.9,1.3) to[out=-90,in=0]++(-0.3,-0.33);
\draw [blue, densely dashed] (P2) ++ (-1.2,2) to[out=180,in=90] ++(-0.3,-0.5) to[out=-90,in=180] ++(0.3,-0.53);
%%左側のからみ

\draw [blue] (P2) ++(-0.9,1.6) --++(2.4,0);
\draw [blue] (P2) ++(-0.9,1.3) --++(2.4,0);
%%中央

\draw [blue] (P2) ++ (1.5,1.6) to[out=90,in=0]++(-0.3,0.4);
\draw [blue] (P2) ++ (1.5,1.3) to[out=-90,in=0]++(-0.3,-0.33);
\draw [blue, densely dashed] (P2) ++ (1.2,2) to[out=180,in=90] ++(-0.3,-0.5) to[out=-90,in=180] ++(0.3,-0.53);
%%左側のからみ

%P2attaching circle

\draw (P2) ++(0,2)node[above]{$(\pi_{0} |_{\partial M_{2}})^{-1} (\alpha_{1})$};

\draw[->,>=stealth] (P2) ++ (-1.5,-2.5) --++(1.5,0);
\draw (P2) ++ (-1.5,-2.5) --++(3,0) node[right]{$\alpha_{1}$};

\draw[->,>=stealth] (P2) ++ (3,0) --++(1,0);
\draw (P2) ++(3.5,0) node[below]{$\pi_{0}$};

%\fill [red!20!white] (P) --++(1.9,0.6) arc (17.5:342.5:2cm) --cycle;
%\draw (Q2) ++ (0,1.5) node[right] {$A$};
%領域A塗りつぶし
\draw (Q2) --++ (1.9,-0.6) arc (-17.5:17.5:2cm) --cycle;
%\draw (P) ++(1.5,0) node[above] {$B$};
%領域B塗りつぶし
%\fill [green!20!white] (P) --++ (0.95,0.3) arc (17.5:342.5:1cm) --cycle;
%領域V_i塗りつぶし

\fill [black] (Q2) circle (0.06);
%\draw (Q2) node[below] {$p_i$};
%p_i

\draw (Q2) circle (0.8);
%\draw (Q2)++(-0.3,0) node[above]{$V_i$};
%V_i
\draw (Q2) circle (2);
%U_i
\draw[red,thick] (Q2) ++(0.95+0.3,-0.4) --++ (0,0.8);
\draw (Q2) ++(0.95+0.3,-0.4) node[below]{$S$};
%S
\draw[blue,thick] (Q2) ++(0.95+0.6,-0.48) to[out=150,in=-90]++(-0.4,0.48) to[out=90,in=-150] ++(0.4,0.48); 
\draw [->,>=stealth,blue,thick] (Q2) ++ (1.15,0.05);
\draw (Q2) ++(0.95+0.6,0.48) ++(-0.15,-0.1) to[out=-30,in=-150]++(0.7,0.5)node[right]{$\alpha_{1}$};
%\alpha_{1}

%%%%↑ここまで(P2)、(Q2)%%%%

\draw[->,>=stealth] (P2) ++(3.5,-2.5) --++(0,-1);
\draw (P2) ++ (3.5,-3) node[right]{(ii)};

%%%%↓ここから(P3)、(Q3)%%%%
\draw (P3) ++ (-1.5,2) --++ (3,0);
%\draw (P3) ++ (-1.5,1) --++ (3,0);
\draw (P3) ++ (-1.5,1) [out=0,in=90] to++(1,-1);
\draw (P3) ++ (1.5,1) [out=180,in=90] to++(-1,-1);
%%上筒
\draw (P3) ++ (-1.5,1) [out =0,in=-90] to++ (0.3,0.5) [out=90,in=0] to++ (-0.3,0.5);
\draw (P3) ++ (-1.5,1) [out =180,in=-90] to++ (-0.3,0.5) [out=90,in=180] to++ (0.3,0.5);
%%左
\draw (P3) ++ (-1.5,1) ++ (3,0) [out =0,in=-90] to++ (0.3,0.5) [out=90,in=0] to++ (-0.3,0.5);
\draw [densely dashed](P3) ++ (-1.5,1)++ (3,0) [out =180,in=-90] to++ (-0.3,0.5) [out=90,in=180] to++ (0.3,0.5);
%%右
%P2ファイバーの上側

%\draw (P3) ++ (-1.5,-1) --++ (3,0);
\draw (P3) ++ (-1.5,-2) --++ (3,0);
\draw (P3) ++ (-1.5,-1) [out=0,in=-90] to++(1,1);
\draw (P3) ++ (1.5,-1) [out=180,in=-90] to++(-1,1);
%%下筒
\draw (P3) ++ (-1.5,-2) [out =0,in=-90] to++ (0.3,0.5) [out=90,in=0] to++ (-0.3,0.5);
\draw (P3) ++ (-1.5,-2) [out =180,in=-90] to++ (-0.3,0.5) [out=90,in=180] to++ (0.3,0.5);
%%左
\draw (P3) ++ (-1.5,-2) ++ (3,0) [out =0,in=-90] to++ (0.3,0.5) [out=90,in=0] to++ (-0.3,0.5);
\draw [densely dashed](P3) ++ (-1.5,-2)++ (3,0) [out =180,in=-90] to++ (-0.3,0.5) [out=90,in=180] to++ (0.3,0.5);
%%右
%P2ファイバーの下側

\draw [blue] (P3) ++ (-0.9,1.6) to[out=90,in=0]++(-0.3,0.4);
\draw [blue] (P3) ++ (-0.9,1.3) to[out=-90,in=0]++(-0.3,-0.33);
\draw [blue, densely dashed] (P3) ++ (-1.2,2) to[out=180,in=90] ++(-0.3,-0.5) to[out=-90,in=180] ++(0.3,-0.53);
%%左側のからみ

\draw [blue] (P3) ++(-0.9,1.6) --++(2.4,0);
\draw [blue] (P3) ++(-0.9,1.3) --++(2.4,0);
%%中央

\draw [blue] (P3) ++ (1.5,1.6) to[out=90,in=0]++(-0.3,0.4);
\draw [blue] (P3) ++ (1.5,1.3) to[out=-90,in=0]++(-0.3,-0.33);
\draw [blue, densely dashed] (P3) ++ (1.2,2) to[out=180,in=90] ++(-0.3,-0.5) to[out=-90,in=180] ++(0.3,-0.53);
%%左側のからみ

%P2attaching circle

\draw (P3) ++(0,2)node[above]{$(\pi_{0} |_{\partial (M_{1} \cup \pi^{-1}_{0} (A) )})^{-1} (\alpha_{2})$};

\draw[->,>=stealth] (P3) ++ (-1.5,-2.5) --++(1.5,0);
\draw (P3) ++ (-1.5,-2.5) --++(3,0) node[right]{$\alpha_{2}$};

\draw[->,>=stealth] (P3) ++ (3,0) --++(1,0);
\draw (P3) ++(3.5,0) node[below]{$\pi_{0}$};

%\fill [red!20!white] (P) --++(1.9,0.6) arc (17.5:342.5:2cm) --cycle;
%\draw (Q3) ++ (0,1.5) node[right] {$A$};
%領域A塗りつぶし
\draw (Q3) --++ (1.9,-0.6) arc (-17.5:17.5:2cm) --cycle;
%\draw (P) ++(1.5,0) node[above] {$B$};
%領域B塗りつぶし
%\fill [green!20!white] (P) --++ (0.95,0.3) arc (17.5:342.5:1cm) --cycle;
%領域V_i塗りつぶし

\fill [black] (Q3) circle (0.06);
%\draw (Q3) node[below] {$p_i$};
%p_i

\draw (Q3) circle (0.8);
%\draw (Q3)++(-0.3,0) node[above]{$V_i$};
%V_i
\draw (Q3) circle (2);
%U_i
\draw[red,thick] (Q3) ++(0.95+0.3,-0.4) --++ (0,0.8);
\draw (Q3) ++(0.95+0.3,-0.4) node[below]{$S$};
%S
\draw[blue, thick] (Q3)++(0.76,0.24) arc (17.5:342.5:0.8cm);
\draw [->,>=stealth,blue,thick] (Q3) ++ (-0.8,0.05);
\draw[blue, thick] (Q3) ++(1.9,-0.6) --++(-1.15,0.36);
\draw[blue, thick] (Q3) ++(1.9,0.6) --++(-1.15,-0.36);
\draw (Q3) ++(0.95+0.6,0.48)  to[out=-30,in=-150]++(0.6,0.5)node[right]{$\alpha_{2}$};

%\alpha_{2}

%%%%↑ここまで(P3)、(Q3)%%%%

\draw[->,>=stealth] (P3) ++(8.5,2) --++(2,5);
\draw (P3) ++ (9.7,3.8) node[below]{(iii)};

%%%%↓ここから(P4)、(Q4)%%%%

\draw (P4) ++ (-1.5,0) circle [x radius =0.7, y radius = 2];
\draw (P4) ++ (-1.5,2) --++ (3,0);
\draw (P4) ++ (-1.5,-2) --++ (3,0);

\draw (P4) ++ (1.5,-2) arc [x radius =0.8, y radius = 2, start angle =-90, end angle =90];
\draw [densely dashed](P4) ++ (1.5,2) arc [x radius =0.8, y radius = 2, start angle =90, end angle =270];
%%外円周

\draw [densely dashed](P4)++ (-1.5,1.5) --++ (3,0);
\draw[densely dashed] (P4) ++ (-1.5,0.5) [out=0,in=90] to++(1,-0.5);
\draw[densely dashed] (P4) ++ (1.5,0.5) [out=180,in=90] to++(-1,-0.5);
%%上筒
\draw(P4) ++ (-1.5,0.5) [out =0,in=-90] to++ (0.3,0.5) [out=90,in=0] to++ (-0.3,0.5);
\draw (P4) ++ (-1.5,0.5) [out =180,in=-90] to++ (-0.3,0.5) [out=90,in=180] to++ (0.3,0.5);
%%左
\draw[densely dashed] (P4) ++ (-1.5,0.5) ++ (3,0) [out =0,in=-90] to++ (0.3,0.5) [out=90,in=0] to++ (-0.3,0.5);
\draw [densely dashed](P4) ++ (-1.5,0.5)++ (3,0) [out =180,in=-90] to++ (-0.3,0.5) [out=90,in=180] to++ (0.3,0.5);
%%右
%P2ファイバーの上側

%\draw (P4) ++ (-1.5,-1) --++ (3,0);
\draw [densely dashed](P4) ++ (-1.5,-1.5) --++ (3,0);
\draw [densely dashed](P4) ++ (-1.5,-0.5) [out=0,in=-90] to++(1,0.5);
\draw [densely dashed](P4) ++ (1.5,-0.5) [out=180,in=-90] to++(-1,0.5);
%%下筒
\draw (P4) ++ (-1.5,-1.5) [out =0,in=-90] to++ (0.3,0.5) [out=90,in=0] to++ (-0.3,0.5);
\draw(P4) ++ (-1.5,-1.5) [out =180,in=-90] to++ (-0.3,0.5) [out=90,in=180] to++ (0.3,0.5);
%%左
\draw [densely dashed](P4) ++ (-1.5,-1.5) ++ (3,0) [out =0,in=-90] to++ (0.3,0.5) [out=90,in=0] to++ (-0.3,0.5);
\draw [densely dashed](P4) ++ (-1.5,-1.5)++ (3,0) [out =180,in=-90] to++ (-0.3,0.5) [out=90,in=180] to++ (0.3,0.5);
%%右
%P2ファイバーの下側

%\draw [blue] (P4) ++ (-0.9,1.6) to[out=90,in=0]++(-0.3,0.4);
%\draw [blue] (P4) ++ (-0.9,1.3) to[out=-90,in=0]++(-0.3,-0.33);
%\draw [blue, densely dashed] (P4) ++ (-1.2,2) to[out=180,in=90] ++(-0.3,-0.5) to[out=-90,in=180] ++(0.3,-0.53);
%%左側のからみ

%\draw [blue] (P4) ++(-0.9,1.6) --++(2.4,0);
%\draw [blue] (P4) ++(-0.9,1.3) --++(2.4,0);
%%中央

%\draw [blue] (P4) ++ (1.5,1.6) to[out=90,in=0]++(-0.3,0.4);
%\draw [blue] (P4) ++ (1.5,1.3) to[out=-90,in=0]++(-0.3,-0.33);
%\draw [blue, densely dashed] (P4) ++ (1.2,2) to[out=180,in=90] ++(-0.3,-0.5) to[out=-90,in=180] ++(0.3,-0.53);
%%左側のからみ

\draw [blue] (P4) ++(-1.3,0.35) --++(0.48,-0.33) --++(3.1,0);
\draw [blue] (P4) ++ (-1.3,0.35) arc [x radius =0.4, y radius =0.75, start angle = -60, end angle =270] --++(0.7,-0.5) --++(3.1,0);
\draw [blue, densely dashed] (P4) ++(-1.3,0.35)++(3,0) --++(0.48,-0.33) ;
\draw [blue, densely dashed] (P4) ++ (-1.3,0.35)++(3,0) arc [x radius =0.4, y radius =0.75, start angle = -60, end angle =270] --++(0.7,-0.5);
%P2attaching circle

%ぬり直しコーナー↓
\draw [preaction={draw = white, line width=3 pt}, densely dashed](P4)++ (-1.5,1.5) --++ (3,0);
\draw[preaction={draw = white, line width=3 pt}, densely dashed] (P4) ++ (1.5,0.5) [out=180,in=90] to++(-1,-0.5);
\draw [preaction={draw = white, line width=3 pt}, densely dashed](P4) ++ (-1.5,-1.5) --++ (3,0);
\draw [preaction={draw = white, line width=3 pt}, densely dashed](P4) ++ (1.5,-0.5) [out=180,in=-90] to++(-1,0.5);
\draw[preaction={draw = white, line width=2 pt}] (P4) ++ (-1.5,0) circle [x radius =0.7, y radius = 2];
\draw (P4) ++ (-1.5,-1.5) [out =0,in=-90] to++ (0.3,0.5) [out=90,in=0] to++ (-0.3,0.5);
\draw(P4) ++ (-1.5,-1.5) [out =180,in=-90] to++ (-0.3,0.5) [out=90,in=180] to++ (0.3,0.5);
\draw (P4) ++ (-1.5,2) --++ (3,0);
\draw (P4) ++ (-1.5,-2) --++ (3,0);
\draw [blue] (P4) ++(-1.3,0.35) --++(0.48,-0.33) --++(3.1,0);
\draw [blue] (P4) ++ (-1.3,0.35) arc [x radius =0.4, y radius =0.75, start angle = -60, end angle =270] --++(0.7,-0.5) --++(3.1,0);
%ぬり直しコーナー↑

\draw[preaction={draw = white, line width=2 pt}] (P4) ++(0,2)node[above]{$\pi_{0}^{-1} (\alpha_{2})$};

\draw[->,>=stealth] (P4) ++ (-1.5,-2.5) --++(1.5,0);
\draw (P4) ++ (-1.5,-2.5) --++(3,0) node[right]{$\alpha_{2}$};

\draw[->,>=stealth] (P4) ++ (3,0) --++(1,0);
\draw (P4) ++(3.5,0) node[below]{$\pi_{0}$};

%\fill [red!20!white] (P) --++(1.9,0.6) arc (17.5:342.5:2cm) --cycle;
%\draw (Q4) ++ (0,1.5) node[right] {$A$};
%領域A塗りつぶし
\draw (Q4) --++ (1.9,-0.6) arc (-17.5:17.5:2cm) --cycle;
%\draw (P) ++(1.5,0) node[above] {$B$};
%領域B塗りつぶし
%\fill [green!20!white] (P) --++ (0.95,0.3) arc (17.5:342.5:1cm) --cycle;
%領域V_i塗りつぶし

\fill [black] (Q4) circle (0.06);
%\draw (Q4) node[below] {$p_i$};
%p_i

\draw (Q4) circle (0.8);
%\draw (Q4)++(-0.3,0) node[above]{$V_i$};
%V_i
\draw (Q4) circle (2);
%U_i
\draw[red,thick] (Q4) ++(0.95+0.3,-0.4) --++ (0,0.8);
\draw (Q4) ++(0.95+0.3,-0.4) node[below]{$S$};
%S
\draw[blue, thick] (Q4)++(0.76,0.24) arc (17.5:342.5:0.8cm);
\draw [->,>=stealth,blue,thick] (Q4) ++ (-0.8,0.05);
\draw[blue, thick] (Q4) ++(1.9,-0.6) --++(-1.15,0.36);
\draw[blue, thick] (Q4) ++(1.9,0.6) --++(-1.15,-0.36);
\draw (Q4) ++(0.95+0.6,0.48)  to[out=-30,in=-150]++(0.6,0.5)node[right]{$\alpha_{2}$};

%\alpha_{2}

%%%%↑ここまで(P4)、(Q4)%%%%

\draw[->,>=stealth] (P4) ++(3.5,-2.5) --++(0,-1);
\draw (P4) ++ (3.5,-3) node[right]{(iv)};

%%%%↓ここから(P5)、(Q5)%%%%
\draw (P5) ++ (-1.5,0) circle [x radius =0.7, y radius = 2];
\draw (P5) ++ (-1.5,2) --++ (3,0);
\draw (P5) ++ (-1.5,-2) --++ (3,0);

\draw (P5) ++ (1.5,-2) arc [x radius =0.8, y radius = 2, start angle =-90, end angle =90];
\draw [densely dashed](P5) ++ (1.5,2) arc [x radius =0.8, y radius = 2, start angle =90, end angle =120];
\draw [densely dashed](P5) ++ (1.5,-2) arc [x radius =0.8, y radius = 2, start angle =-90, end angle =-120];
%%外円周

\draw [densely dashed](P5)++ (-1.5,1.5) --++ (1,0);
\draw[densely dashed] (P5) ++ (-1.5,0.5) --++(1,0);
\draw[densely dashed] (P5) ++ (1.5,0.5) --++(-0.5,0);
\draw[densely dashed] (P5) ++ (1.5,1.5) --++(-0.5,0);
%%上筒
\draw(P5) ++ (-1.5,0.5) [out =0,in=-90] to++ (0.3,0.5) [out=90,in=0] to++ (-0.3,0.5);
\draw (P5) ++ (-1.5,0.5) [out =180,in=-90] to++ (-0.3,0.5) [out=90,in=180] to++ (0.3,0.5);
%%左
\draw[densely dashed] (P5) ++ (-1.5,0.5) ++ (3,0) [out =0,in=-90] to++ (0.3,0.5) [out=90,in=0] to++ (-0.3,0.5);
\draw [densely dashed](P5) ++ (-1.5,0.5)++ (3,0) [out =180,in=-90] to++ (-0.3,0.5) [out=90,in=180] to++ (0.3,0.5);
%%右
%P2ファイバーの上側

%\draw (P5) ++ (-1.5,-1) --++ (3,0);
\draw [densely dashed](P5) ++ (-1.5,-1.5) --++ (1,0);
\draw [densely dashed](P5) ++ (-1.5,-0.5) --++(1,0);
\draw [densely dashed](P5) ++ (1.5,-0.5) --++(-0.5,0);
\draw [densely dashed](P5) ++ (1.5,-1.5) --++(-0.5,0);
%下筒
\draw (P5) ++ (-1.5,-1.5) [out =0,in=-90] to++ (0.3,0.5) [out=90,in=0] to++ (-0.3,0.5);
\draw(P5) ++ (-1.5,-1.5) [out =180,in=-90] to++ (-0.3,0.5) [out=90,in=180] to++ (0.3,0.5);
%%左
\draw [densely dashed](P5) ++ (-1.5,-1.5) ++ (3,0) [out =0,in=-90] to++ (0.3,0.5) [out=90,in=0] to++ (-0.3,0.5);
\draw [densely dashed](P5) ++ (-1.5,-1.5)++ (3,0) [out =180,in=-90] to++ (-0.3,0.5) [out=90,in=180] to++ (0.3,0.5);
%%右
%P2ファイバーの下側

\draw[densely dotted, thick] (P5) ++(-0.5,1.7) --++(0,-3.4) --++ (1.5,0) --++(0,3.4) --cycle;
\draw (P5) ++(0.2,0.6) node{$\beta$};
%ブレイドの箱

%\draw [blue] (P5) ++ (-0.9,1.6) to[out=90,in=0]++(-0.3,0.4);
%\draw [blue] (P5) ++ (-0.9,1.3) to[out=-90,in=0]++(-0.3,-0.33);
%\draw [blue, densely dashed] (P5) ++ (-1.2,2) to[out=180,in=90] ++(-0.3,-0.5) to[out=-90,in=180] ++(0.3,-0.53);
%%左側のからみ

%\draw [blue] (P5) ++(-0.9,1.6) --++(2.4,0);
%\draw [blue] (P5) ++(-0.9,1.3) --++(2.4,0);
%%中央

%\draw [blue] (P5) ++ (1.5,1.6) to[out=90,in=0]++(-0.3,0.4);
%\draw [blue] (P5) ++ (1.5,1.3) to[out=-90,in=0]++(-0.3,-0.33);
%\draw [blue, densely dashed] (P5) ++ (1.2,2) to[out=180,in=90] ++(-0.3,-0.5) to[out=-90,in=180] ++(0.3,-0.53);
%%左側のからみ

\draw [blue] (P5) ++(-1.3,0.35) --++(0.48,-0.33) --++(3.1,0);
\draw [blue] (P5) ++ (-1.3,0.35) arc [x radius =0.4, y radius =0.75, start angle = -60, end angle =270] --++(0.7,-0.5) --++(3.1,0);
\draw [blue, densely dashed] (P5) ++(-1.3,0.35)++(3,0) --++(0.48,-0.33) ;
\draw [blue, densely dashed] (P5) ++ (-1.3,0.35)++(3,0) arc [x radius =0.4, y radius =0.75, start angle = -60, end angle =270] --++(0.7,-0.5);
%P2attaching circle

%ぬり直しコーナー
\draw[preaction={draw = white, line width = 3 pt}, densely dashed] (P5) ++ (1.5,0.5) --++(-0.4,0);
\draw[preaction={draw = white, line width = 3 pt}, densely dashed] (P5) ++ (1.5,1.5) --++(-0.4,0);
\draw [preaction={draw = white, line width = 3 pt},blue] (P5) ++(-1.3,0.35) ++(0.48,-0.33) --++(3.1,0);
\draw [preaction={draw = white, line width = 3 pt},blue] (P5) ++ (-1.3,0.35) arc [x radius =0.4, y radius =0.75, start angle = -60, end angle =270] --++(0.7,-0.5) --++(3.1,0);
%ぬり直しコーナー

\draw (P5) ++(0,2)node[above]{$\pi_{0}^{-1} (\alpha_{3})$};

\draw[->,>=stealth] (P5) ++ (-1.5,-2.5) --++(1.5,0);
\draw (P5) ++ (-1.5,-2.5) --++(3,0) node[right]{$\alpha_{3}$};

\draw[->,>=stealth] (P5) ++ (3,0) --++(1,0);
\draw (P5) ++(3.5,0) node[below]{$\pi_{0}$};

\draw (P5) ++(0,-3.2) node{($\beta$ is the braid of};
\draw (P5) ++(0,-3.9) node{local monodromy.)};

\draw (P5) ++(2,-0.25) to[out=-40, in=180] ++(1.4,-2)node[right] {$C'_{0}$};
%c'_i名前

%\fill [red!20!white] (P) --++(1.9,0.6) arc (17.5:342.5:2cm) --cycle;
%\draw (Q4) ++ (0,1.5) node[right] {$A$};
%領域A塗りつぶし
\draw (Q5) --++ (1.9,-0.6) arc (-17.5:17.5:2cm) --cycle;
%\draw (P) ++(1.5,0) node[above] {$B$};
%領域B塗りつぶし
%\fill [green!20!white] (P) --++ (0.95,0.3) arc (17.5:342.5:1cm) --cycle;
%領域V_i塗りつぶし

\fill [black] (Q5) circle (0.06);
%\draw (Q5) node[below] {$p_i$};
%p_i

\draw (Q5) circle (0.8);
%\draw (Q5)++(-0.3,0) node[above]{$V_i$};
%V_i
\draw (Q5) circle (2);
%U_i
\draw[red,thick] (Q5) ++(0.95+0.3,-0.4) --++ (0,0.8);
\draw (Q5) ++(0.95+0.3,-0.4) node[below]{$S$};
%S
\draw [->,>=stealth,blue,thick] (Q5) ++ (-2,0.05);
\draw[blue, thick] (Q5) ++(1.9,0.6) arc (17.5:342.5:2cm);
\draw (Q5) ++(1.8,0.8)  to[out=-30,in=-150]++(0.3,0.5)node[right]{$\alpha_{3}$};
%\alpha_{3}

\end{tikzpicture}
\caption{Moving the attaching circle}
\label{fig:attachingcircle}
\end{figure}

\end{remark}

\begin{figure}[htbp]
\centering
\begin{tikzpicture}[scale=0.75]
%%%1つ目

\coordinate (R5) at (0,8.5);
\coordinate (R6) at (3,3.5);

\draw (R5) ++ (-1.5,0) circle [x radius =0.7, y radius = 2];
\draw (R5) ++ (-1.5,2) --++ (9,0);
\draw (R5) ++ (-1.5,-2) --++ (9,0);

\draw (R5) ++ (7.5,-2) arc [x radius =0.8, y radius = 2, start angle =-90, end angle =90];
\draw [densely dashed](R5) ++ (7.5,2) arc [x radius =0.8, y radius = 2, start angle =90, end angle =270];
%%外円周

\draw [densely dashed](R5)++ (-1.5,1.5) --++ (9,0);
%\draw[densely dashed] (R5) ++ (-1.5,0.5) --++(1,0);
%\draw[densely dashed] (R5) ++ (7.5,0.5) --++(-1,0);
%\draw[densely dashed] (R5) ++ (7.5,1.5) --++(-1,0);
%%上筒
\draw(R5) ++ (-1.5,0.5) [out =0,in=-90] to++ (0.3,0.5) [out=90,in=0] to++ (-0.3,0.5);
\draw (R5) ++ (-1.5,0.5) [out =180,in=-90] to++ (-0.3,0.5) [out=90,in=180] to++ (0.3,0.5);
%%左
\draw[densely dashed] (R5) ++ (-1.5,0.5) ++ (9,0) [out =0,in=-90] to++ (0.3,0.5) [out=90,in=0] to++ (-0.3,0.5);
\draw [densely dashed](R5) ++ (-1.5,0.5)++ (9,0) [out =180,in=-90] to++ (-0.3,0.5) [out=90,in=180] to++ (0.3,0.5);
%%右
%P2ファイバーの上側

%\draw (R5) ++ (-1.5,-1) --++ (3,0);
\draw [densely dashed](R5) ++ (-1.5,-1.5) --++ (9,0);
%\draw [densely dashed](R5) ++ (-1.5,-0.5) --++(1,0);
%\draw [densely dashed](R5) ++ (7.5,-0.5) --++(-1,0);
%\draw [densely dashed](R5) ++ (7.5,-1.5) --++(-1,0);
%下筒
\draw (R5) ++ (-1.5,-1.5) [out =0,in=-90] to++ (0.3,0.5) [out=90,in=0] to++ (-0.3,0.5);
\draw(R5) ++ (-1.5,-1.5) [out =180,in=-90] to++ (-0.3,0.5) [out=90,in=180] to++ (0.3,0.5);
%%左
\draw [densely dashed](R5) ++ (-1.5,-1.5) ++ (9,0) [out =0,in=-90] to++ (0.3,0.5) [out=90,in=0] to++ (-0.3,0.5);
\draw [densely dashed](R5) ++ (-1.5,-1.5)++ (9,0) [out =180,in=-90] to++ (-0.3,0.5) [out=90,in=180] to++ (0.3,0.5);
%%右
%P2ファイバーの下側

\draw[densely dashed] (R5) ++ (-1.5,0.5)  to[out=0,in=90] ++(1.2,-0.5) to[out=-90,in=0] ++ (-1.2,-0.5);

\draw [densely dashed] (R5) ++ (7.5,0.5) to[out=180,in=90] ++(-1,-0.5) to[out=-90,in=180] ++ (1,-0.5);

\draw [blue] (R5) ++(-1.3,0.35) --++(0.48,-0.33) --++(9.1,0);
\draw [blue] (R5) ++ (-1.3,0.35) arc [x radius =0.4, y radius =0.75, start angle = -60, end angle =270] --++(0.7,-0.5) --++(9.1,0);
\draw [blue, densely dashed] (R5) ++(-1.3,0.35)++(9,0) --++(0.48,-0.33) ;
\draw [blue, densely dashed] (R5) ++ (-1.3,0.35)++(9,0) arc [x radius =0.4, y radius =0.75, start angle = -60, end angle =270] --++(0.7,-0.5);
%P2attaching circle

%ぬり直しコーナー
\draw[preaction={draw = white, line width = 3 pt}, densely dashed] (R5) ++ (7.5,0.5) --++(-0.2,0);
\draw[preaction={draw = white, line width = 3 pt}, densely dashed] (R5) ++ (7.5,1.5) --++(-0.5,0);
\draw [preaction={draw = white, line width = 3 pt},blue] (R5) ++(-1.3,0.35) ++(0.48,-0.33) --++(3.1,0);
\draw [preaction={draw = white, line width = 3 pt},blue] (R5) ++ (-1.3,0.35) arc [x radius =0.4, y radius =0.75, start angle = -60, end angle =270] --++(0.7,-0.5) --++(3.1,0);
%ぬり直しコーナー

\draw (R6) --++ (1.425,-0.45) arc (-17.5:17.5:1.5cm) --cycle;
\fill [black] (R6) circle (0.06);
%\draw (R6) node[below] {$p_i$};
%p_i

\draw (R6) circle (0.5);
%\draw (R6)++(-0.3,0) node[above]{$V_i$};
%V_i
\draw (R6) circle (1.5);
%U_i
\draw[red,thick] (R6) ++(0.95,-0.3) --++ (0,0.6);
\draw (R6) ++(0.895,-0) node[right]{$S$};
\draw (R6) ++ (0.95,0) --++(0.2,0);
%S
\draw[blue, thick] (R6)++(0.475,0.15) arc (17.5:342.5:0.5cm);
\draw [->,>=stealth,blue,thick] (R6) ++ (-0.5,0.05);
\draw[blue, thick] (R6) ++(1.425,-0.45) --++(-0.95,0.3);
\draw[blue, thick] (R6) ++(1.425,0.45) --++(-0.95,-0.3);
\draw (R6) ++(0.95++-0.313 ,0.3++-0.1)  to[out=60,in=-150]++(0.8,1)node[right]{$\alpha_{2}$};
\draw (R6) ++(0.95++0.313 ,0.3++0.1)  to[out=-60,in=150]++(0.5,-0)node[right]{$\ell_{\theta_0, \varepsilon /2}$};
\draw (R6) ++(0.95++0.313 ,-0.3++-0.1)  to[out=60,in=-150]++(0.5,-0)node[right]{$\ell_{2 \pi -\theta_0, \varepsilon /2}$};
\draw (R6) ++ (-0.5,0) to[out=150, in=30] ++(-1.2,0) node [left]{$T_{\varepsilon /2}$};

\draw[red,densely dotted] (R6)++ (0.95,0.3) arc (17.5:342.5:1cm);

\draw (R5) ++ (7,2) node[above]{$\pi_{0}^{-1} (\alpha_{2})$};

\draw (R6) ++ (0,-1.9) node{$(a)$};

\draw (R5) ++ (-0.3,-2.7) node{$\ell_{2\pi- \theta_0, \varepsilon /2}$};
\draw (R5) ++ (-1.5, -2) [out=-90,in=150] to++ (0.2,-0.4);
\draw (R5) ++ (0.3, -2) [out=-90,in=30] to++ (-0.2,-0.4);
\draw (R5) ++(8,0)++ (-1,-2.7) node{$\ell_{\theta_0, \varepsilon /2}$};
\draw (R5) ++(8,0)++ (-2.1, -2) [out=-90,in=150] to++ (0.2,-0.4);
\draw (R5) ++(8,0)++ (-0.3, -2) [out=-90,in=30] to++ (-0.2,-0.4);
\draw (R5) ++ (0.3,-2) [out=-90,in=150] to++ (0.4,-0.3);
\draw (R5) ++(8,0) ++ (-2.1, -2) [out=-90,in=30] to++ (-0.4,-0.3);
\draw (R5) ++ (3,-2.5) node{$T_{\varepsilon /2}$};

\draw[->, line width=3pt] (R6) ++ (0,-2.4) --++(0,-0.7);

%%%二つ目

\coordinate (Q5) at (0,-2);
\coordinate (Q6) at (3,-7);

\draw (Q5) ++ (-1.5,0) circle [x radius =0.7, y radius = 2];
\draw (Q5) ++ (-1.5,2) --++ (9,0);
\draw (Q5) ++ (-1.5,-2) --++ (9,0);

\draw (Q5) ++ (7.5,-2) arc [x radius =0.8, y radius = 2, start angle =-90, end angle =90];
\draw [densely dashed](Q5) ++ (7.5,2) arc [x radius =0.8, y radius = 2, start angle =90, end angle =270];
%%外円周

\draw [densely dashed](Q5)++ (-1.5,1.5) --++ (1.5,0);
%\draw[densely dashed] (Q5) ++ (-1.5,0.5) --++(1,0);
%\draw[densely dashed] (Q5) ++ (7.5,0.5) --++(-1,0);
\draw[densely dashed] (Q5) ++ (7.5,1.5) --++(-1.5,0);
%%上筒
\draw(Q5) ++ (-1.5,0.5) [out =0,in=-90] to++ (0.3,0.5) [out=90,in=0] to++ (-0.3,0.5);
\draw (Q5) ++ (-1.5,0.5) [out =180,in=-90] to++ (-0.3,0.5) [out=90,in=180] to++ (0.3,0.5);
%%左
\draw[densely dashed] (Q5) ++ (-1.5,0.5) ++ (9,0) [out =0,in=-90] to++ (0.3,0.5) [out=90,in=0] to++ (-0.3,0.5);
\draw [densely dashed](Q5) ++ (-1.5,0.5)++ (9,0) [out =180,in=-90] to++ (-0.3,0.5) [out=90,in=180] to++ (0.3,0.5);
%%右
%P2ファイバーの上側

%\draw (Q5) ++ (-1.5,-1) --++ (3,0);
\draw [densely dashed](Q5) ++ (-1.5,-1.5) --++ (1.5,0);
%\draw [densely dashed](Q5) ++ (-1.5,-0.5) --++(-1,0);
%\draw [densely dashed](Q5) ++ (7.5,-0.5) --++(-1,0);
\draw [densely dashed](Q5) ++ (7.5,-1.5) --++(-1.5,0);
%下筒
\draw (Q5) ++ (-1.5,-1.5) [out =0,in=-90] to++ (0.3,0.5) [out=90,in=0] to++ (-0.3,0.5);
\draw(Q5) ++ (-1.5,-1.5) [out =180,in=-90] to++ (-0.3,0.5) [out=90,in=180] to++ (0.3,0.5);
%%左
\draw [densely dashed](Q5) ++ (-1.5,-1.5) ++ (9,0) [out =0,in=-90] to++ (0.3,0.5) [out=90,in=0] to++ (-0.3,0.5);
\draw [densely dashed](Q5) ++ (-1.5,-1.5)++ (9,0) [out =180,in=-90] to++ (-0.3,0.5) [out=90,in=180] to++ (0.3,0.5);
%%右
%P2ファイバーの下側

%\draw [densely dashed](Q5)++ (-1.5,1.5) ++ (1,0) --++(0.5,0) to[out=0,in=135] ++(0.7,-0.5);
%\draw[densely dashed] (Q5) ++ (-1.5,0.5) ++(1,0) to[out=0,in=135] ++(0.7,-0.5);
\draw [densely dashed] (Q5) ++ (-1.5,0.5) to [out=0, in=90] ++(1.5,-0.5) to[out=-90,in=0] ++ (-1.5,-0.5);
\draw [densely dashed ](Q5) ++ (7.5,0.5)  to [out=180, in=90] ++(-1.5,-0.5) to[out=-90,in=180] ++ (1.5,-0.5) ;

%\draw [densely dashed, thick] (Q5) ++ (-0,-0.1) --++(6,0);
%\draw [densely dashed, semithick] (Q5) ++ (-0,-0.1) --++(1.2,0);
%\draw [densely dashed, semithick] (Q5) ++ (3,-0.1) --++(-1.1,0);
%\draw [densely dashed] (Q5) ++ (0,-1.5) --++(0.5,0) to[out=0,in=180] ++ (2,3) --++(0.5,0);
%\draw [densely dashed] (Q5) ++ (0,1.5) --++(0.5,0) to[out=0,in=110] ++ (0.9,-1.1);
%\draw [densely dashed] (Q5) ++ (3,-1.5) --++(-0.5,0) to[out=180,in=-70] ++ (-0.9,1.1);
\draw[densely dashed] (Q5) ++(0,-0.1) to[out=0,in=180] ++ (3,1.6);
\draw[densely dashed] (Q5) ++(0,-1.5) to[out=0,in=180] ++ (2.9,1.4) --++(0.1,0);
\draw[densely dashed] (Q5) ++ (0,1.5) to[out=0,in=135] ++ (1.1,-0.7);
\draw[densely dashed] (Q5) ++ (0,1.5)++ (1.3,-0.8) ++ (0.2,-0.3) to[out=-45,in=180] ++ (1.4,-0.5) --++(0.1,0);
\draw[densely dashed] (Q5) ++ (0,-0.1) to[out=0,in=150] ++ (1.1,-0.7);
\draw[densely dashed] (Q5) ++ (0,-0.1)++ (1.3,-0.8) ++ (0.2,-0.2) to[out=-45,in=180] ++ (1.4,-0.4) --++(0.1,0);

\draw[densely dashed] (Q5)++(3,0) ++(0,-0.1) to[out=0,in=180] ++ (3,1.6);
\draw[densely dashed] (Q5)++(3,0) ++(0,-1.5) to[out=0,in=180] ++ (2.9,1.4) --++(0.1,0);
\draw[densely dashed] (Q5)++(3,0) ++ (0,1.5) to[out=0,in=135] ++ (1.1,-0.7);
\draw[densely dashed] (Q5)++(3,0) ++ (0,1.5)++ (1.3,-0.8) ++ (0.2,-0.3) to[out=-45,in=180] ++ (1.4,-0.5) --++(0.1,0);
\draw[densely dashed] (Q5)++(3,0) ++ (0,-0.1) to[out=0,in=150] ++ (1.1,-0.7);
\draw[densely dashed] (Q5)++(3,0) ++ (0,-0.1)++ (1.3,-0.8) ++ (0.2,-0.2) to[out=-45,in=180] ++ (1.4,-0.4) --++(0.1,0);

\fill (Q5)++(0.05,-0.1) circle (0.06);
\fill (Q5)++(3,0)++(0,-0.1) circle (0.06);
\fill (Q5)++(5.95,0)++(0,-0.1) circle (0.06);
%ブレイド

\draw [densely dashed] (Q5) ++ (0.05,-0.1) to[out=70, in=-70] ++ (0,1.6);
\draw [densely dashed] (Q5) ++ (0.05,-0.1) to[out=110, in=-110] ++ (0,1.6);
\draw [densely dashed] (Q5) ++ (0.05,-0.1) to[out=-70, in=70] ++ (0,-1.4);
\draw [densely dashed] (Q5) ++ (0.05,-0.1) to[out=-110, in=110] ++ (0,-1.4);

\draw [densely dashed] (Q5) ++ (3,-0.1) to[out=70, in=-70] ++ (0,1.6);
\draw [densely dashed] (Q5) ++ (3,-0.1) to[out=110, in=-110] ++ (0,1.6);
\draw [densely dashed] (Q5) ++ (3,-0.1) to[out=-70, in=70] ++ (0,-1.4);
\draw [densely dashed] (Q5) ++ (3,-0.1) to[out=-110, in=110] ++ (0,-1.4);

\draw [densely dashed] (Q5) ++ (5.95,-0.1) to[out=70, in=-70] ++ (0,1.6);
\draw [densely dashed] (Q5) ++ (5.95,-0.1) to[out=110, in=-110] ++ (0,1.6);
\draw [densely dashed] (Q5) ++ (5.95,-0.1) to[out=-70, in=70] ++ (0,-1.4);
\draw [densely dashed] (Q5) ++ (5.95,-0.1) to[out=-110, in=110] ++ (0,-1.4);
%ファイバーの八の字

\draw [blue] (Q5) ++(-1.3,0.35) --++(0.48,-0.33) --++(9.1,0);
\draw [blue] (Q5) ++ (-1.3,0.35) arc [x radius =0.4, y radius =0.75, start angle = -60, end angle =270] --++(0.7,-0.5) --++(9.1,0);
\draw [blue, densely dashed] (Q5) ++(-1.3,0.35)++(9,0) --++(0.48,-0.33) ;
\draw [blue, densely dashed] (Q5) ++ (-1.3,0.35)++(9,0) arc [x radius =0.4, y radius =0.75, start angle = -60, end angle =270] --++(0.7,-0.5);
%P2attaching circle

%ぬり直しコーナー
\draw[preaction={draw = white, line width = 2.5 pt}, densely dashed] (Q5) ++ (7.5,0.5) --++(-0.5,0);
\draw[preaction={draw = white, line width = 2.5 pt}, densely dashed] (Q5) ++ (7.5,1.5) --++(-0.5,0);
\draw [preaction={draw = white, line width = 2.5 pt},blue] (Q5) ++(-1.3,0.35) ++(0.48,-0.33) --++(7,0);
\draw [preaction={draw = white, line width = 2.5 pt},blue] (Q5) ++ (-1.3,0.35) arc [x radius =0.4, y radius =0.75, start angle = -60, end angle =270] --++(0.7,-0.5) --++(7,0);
%ぬり直しコーナー

\draw (Q6) --++ (1.425,-0.45) arc (-17.5:17.5:1.5cm) --cycle;
\fill [black] (Q6) circle (0.06);
%\draw (Q6) node[below] {$p_i$};
%p_i

\draw (Q6) circle (0.5);
%\draw (Q6)++(-0.3,0) node[above]{$V_i$};
%V_i
\draw (Q6) circle (1.5);
%U_i
\draw[red,thick] (Q6) ++(0.95,-0.3) --++ (0,0.6);
\draw (Q6) ++(0.895,-0) node[right]{$S$};
%S
\draw[blue, thick] (Q6)++(0.95,0.3) arc (17.5:342.5:1cm);
\draw [->,>=stealth,blue,thick] (Q6) ++ (-1,0.05);
\draw[blue, thick] (Q6) ++(1.425,-0.45) --++(-0.495,0.15);
\draw[blue, thick] (Q6) ++(1.425,0.45) --++(-0.495,-0.15);
\draw (Q6) ++(0.95++0.2,0.3++0.05)  to[out=60,in=-150]++(0.6,0.5)node[right]{$\alpha'_{2}$};
\draw (Q6) ++ (-1,0) to[out=150,in=30] ++(-0.7,0) node[left]{$T_{\varepsilon}$};

\draw[red,densely dotted] (Q6)++ (0.92,0.29) arc (17.5:342.5:0.96cm);

\draw (Q5) ++ (7,2) node[above]{$\pi_{0}^{-1} (\alpha'_{2})$};

\draw (Q6) ++(0,-2) node {$(b)$};

\draw (Q5) ++ (-0.3,-2.7) node{$\ell_{2\pi- \theta_0, \varepsilon}$};
\draw (Q5) ++ (-1.5, -2) [out=-90,in=150] to++ (0.2,-0.4);
\draw (Q5) ++ (0, -2) [out=-90,in=30] to++ (-0.2,-0.4);
\draw (Q5) ++(8,0)++ (-1,-2.7) node{$\ell_{\theta_0, \varepsilon}$};
\draw (Q5) ++(8,0)++ (-2, -2) [out=-90,in=150] to++ (0.2,-0.4);
\draw (Q5) ++(8,0)++ (-0.3, -2) [out=-90,in=30] to++ (-0.2,-0.4);
\draw (Q5) ++ (0,-2) [out=-90,in=150] to++ (0.4,-0.3);
\draw (Q5) ++(8,0) ++ (-2, -2) [out=-90,in=30] to++ (-0.4,-0.3);
\draw (Q5) ++ (3,-2.5) node{$T_{\varepsilon}$};

\draw[->, line width =3pt] (Q5) ++ (8.7,2) --++(0.7,0.7);

%%3つ目
\coordinate (S5) at (12,4);
\coordinate (S6) at (15,-0);

\draw (S5) ++ (-1.5,0) circle [x radius =0.7, y radius = 2];
\draw (S5) ++ (-1.5,2) --++ (9,0);
\draw (S5) ++ (-1.5,-2) --++ (9,0);

\draw (S5) ++ (7.5,-2) arc [x radius =0.8, y radius = 2, start angle =-90, end angle =90];
\draw [densely dashed](S5) ++ (7.5,2) arc [x radius =0.8, y radius = 2, start angle =90, end angle =270];
%%外円周

\draw [densely dashed](S5)++ (-1.5,1.5) --++ (1,0);
\draw[densely dashed] (S5) ++ (-1.5,0.5) --++(1,0);
\draw[densely dashed] (S5) ++ (7.5,0.5) --++(-1,0);
\draw[densely dashed] (S5) ++ (7.5,1.5) --++(-1,0);
%%上筒
\draw(S5) ++ (-1.5,0.5) [out =0,in=-90] to++ (0.3,0.5) [out=90,in=0] to++ (-0.3,0.5);
\draw (S5) ++ (-1.5,0.5) [out =180,in=-90] to++ (-0.3,0.5) [out=90,in=180] to++ (0.3,0.5);
%%左
\draw[densely dashed] (S5) ++ (-1.5,0.5) ++ (9,0) [out =0,in=-90] to++ (0.3,0.5) [out=90,in=0] to++ (-0.3,0.5);
\draw [densely dashed](S5) ++ (-1.5,0.5)++ (9,0) [out =180,in=-90] to++ (-0.3,0.5) [out=90,in=180] to++ (0.3,0.5);
%%右
%P2ファイバーの上側

%\draw (S5) ++ (-1.5,-1) --++ (3,0);
\draw [densely dashed](S5) ++ (-1.5,-1.5) --++ (1,0);
\draw [densely dashed](S5) ++ (-1.5,-0.5) --++(1,0);
\draw [densely dashed](S5) ++ (7.5,-0.5) --++(-1,0);
\draw [densely dashed](S5) ++ (7.5,-1.5) --++(-1,0);
%下筒
\draw (S5) ++ (-1.5,-1.5) [out =0,in=-90] to++ (0.3,0.5) [out=90,in=0] to++ (-0.3,0.5);
\draw(S5) ++ (-1.5,-1.5) [out =180,in=-90] to++ (-0.3,0.5) [out=90,in=180] to++ (0.3,0.5);
%%左
\draw [densely dashed](S5) ++ (-1.5,-1.5) ++ (9,0) [out =0,in=-90] to++ (0.3,0.5) [out=90,in=0] to++ (-0.3,0.5);
\draw [densely dashed](S5) ++ (-1.5,-1.5)++ (9,0) [out =180,in=-90] to++ (-0.3,0.5) [out=90,in=180] to++ (0.3,0.5);
%%右
%P2ファイバーの下側

\draw [densely dashed](S5)++ (-1.5,1.5) ++ (1,0) --++(0.5,0) to[out=0,in=135] ++(0.7,-0.5);
\draw[densely dashed] (S5) ++ (-1.5,0.5) ++(1,0) to[out=0,in=135] ++(0.7,-0.5);

\draw[densely dashed] (S5) ++ (-1.5,-1.5) ++(5,0) --++(-1.7,0) to[out=180,in=-45] ++(-0.8,0.5);
\draw[densely dashed] (S5) ++ (-1.5,-0.5) ++(5,0) --++(-1.4,0) to[out=180,in=-45] ++(-0.7,0.5);

\draw [densely dashed](S5) ++ (-1.5,-1.5) ++ (1,0) --++(0.5,0) to[out=0,in=180] ++(2,2) --++(1.5,0);
\draw [densely dashed](S5) ++ (-1.5,-0.5) ++(1,0) --++(0,0) to[out=0,in=180] ++(2,2) --++(2,0);
%%前半

\draw [densely dashed](S5) ++(4,0)++ (-1.5,1.5) ++ (1,0) --++(0.5,0) to[out=0,in=135] ++(0.7,-0.5);
\draw[densely dashed] (S5)++(4,0) ++ (-1.5,0.5) ++(1,0) to[out=0,in=135] ++(0.7,-0.5);

\draw[densely dashed] (S5)++(4,0) ++ (-1.5,-1.5) ++(5,0) --++(-1.7,0) to[out=180,in=-45] ++(-0.8,0.5);
\draw[densely dashed] (S5)++(4,0) ++ (-1.5,-0.5) ++(5,0) --++(-1.4,0) to[out=180,in=-45] ++(-0.7,0.5);

\draw [densely dashed](S5) ++(4,0)++ (-1.5,-1.5) ++ (1,0) --++(0.5,0) to[out=0,in=180] ++(2,2) --++(1.5,0);
\draw [densely dashed](S5) ++(4,0)++ (-1.5,-0.5) ++(1,0) --++(0,0) to[out=0,in=180] ++(2,2) --++(2,0);

%%後半
%中央のブレイド

\draw [blue] (S5) ++(-1.3,0.35) --++(0.48,-0.33) --++(9.1,0);
\draw [blue] (S5) ++ (-1.3,0.35) arc [x radius =0.4, y radius =0.75, start angle = -60, end angle =270] --++(0.7,-0.5) --++(9.1,0);
\draw [blue, densely dashed] (S5) ++(-1.3,0.35)++(9,0) --++(0.48,-0.33) ;
\draw [blue, densely dashed] (S5) ++ (-1.3,0.35)++(9,0) arc [x radius =0.4, y radius =0.75, start angle = -60, end angle =270] --++(0.7,-0.5);
%P2attaching circle

%ぬり直しコーナー
\draw[preaction={draw = white, line width = 3 pt}, densely dashed] (S5) ++ (7.5,0.5) --++(-0.5,0);
\draw[preaction={draw = white, line width = 3 pt}, densely dashed] (S5) ++ (7.5,1.5) --++(-0.5,0);
\draw [preaction={draw = white, line width = 3 pt},blue] (S5) ++(-1.3,0.35) ++(0.48,-0.33) --++(3.1,0);
\draw [preaction={draw = white, line width = 3 pt},blue] (S5) ++ (-1.3,0.35) arc [x radius =0.4, y radius =0.75, start angle = -60, end angle =270] --++(0.7,-0.5) --++(3.1,0);
%ぬり直しコーナー

\draw (S6) --++ (1.425,-0.45) arc (-17.5:17.5:1.5cm) --cycle;
\fill [black] (S6) circle (0.06);
%\draw (S6) node[below] {$p_i$};
%p_i

\draw (S6) circle (0.5);
%\draw (S6)++(-0.3,0) node[above]{$V_i$};
%V_i
\draw (S6) circle (1.5);
%U_i
\draw[red,thick] (S6) ++(0.95,-0.3) --++ (0,0.6);
\draw (S6) ++(0.895,-0) node[right]{$S$};
%S
\draw[blue, thick] (S6)++(1.425,0.45) arc (17.5:342.5:1.5cm);
\draw [->,>=stealth,blue,thick] (S6) ++ (-1.5,0.05);
\draw (S6) ++(1.4,0.6)  to[out=60,in=-150]++(0.6,0.5)node[right]{$\alpha_{3}$};

\draw[red,densely dotted] (S6)++ (0.95,0.3) arc (17.5:342.5:1cm);

\draw (S5) ++ (7,2) node[above]{$\pi_{0}^{-1} (\alpha_{3})$};

\draw (S6) ++ (0,-2) node {$(c)$};

\end{tikzpicture}
\caption{Moving the attaching circle from $\pi_{0}^{-1}(\alpha_{2})$ to $\pi_{0}^{-1} (\alpha_{3})$}
\label{fig:precisemove}
\end{figure}

\begin{example}
\label{exam:specific}
Let us illustrate the motion of the attaching circle in the last step (iv) in Remark \ref{rmk:attachingcircle} with a specific example.
Suppose that the curve is defined by $f(x,y) = x^2 - y^2$ and define a tubular neighborhood $\nu(C)$ of the curve by $\nu(C)=\{(x,y) \in \bC^2 \mid |x^2 - y^2| \leq \varepsilon^2 \}$ where $\varepsilon$ is a small positive number. 
Consider the topology around the only singularity $(0,0)$ of this curve.
Define the neighborhoods of $0 \in \bC$ as $U=\{x \in \bC \mid |x|\leq 2 \varepsilon\}$ and $V=\{x \in \bC \mid |x| \leq \varepsilon /2\}$.
Since we have that  
\[
\pi^{-1}(x_{0}) \cap \partial \nu (C)
=\partial \{(x_0, y) \in \bC^2 \mid |y-x_0||y+x_{0}| \leq \varepsilon^2 \},
\]
for $x_{0} \in \bC$ (so-called Cassinian oval), the number of connected components is two when $|x_{0}| > \varepsilon $ and one otherwise.
Therefore, the neighborhoods $U$ and $V$ match the situations we have been considering. 
When $|x_{0}| > \varepsilon$,  $\pi^{-1}_0 (x_0) \cap \partial \nu (C)$ consists of two simply closed curves. Each of the points $\pm x_{0} = \pm |x_0| e^{i \arg x_{0}}$ lies in the interior region bounded by one of these curves, and we refer to these points as the ``centers'' of the corresponding lobes. Also, when, $|x_{0}| = \varepsilon$, $\pi^{-1}_0 (x_0) \cap \partial \nu (C)$ is the union of two lobes centered at $\pm \varepsilon e^{i \arg x_{0}}$, intersecting at the origin.
The critical value curve $S$ is described as $S=\{x \in \bC \mid |x| = \varepsilon \}$.
For real numbers $r$ and $\theta_0$, define the segment 
$\ell_{\theta_0, r} = \{t e^{i \theta_{0}} \mid r \leq t \leq 2 \varepsilon\}$ 
and the arc
$T_{r} = \{r e^{i \theta} \mid \theta_{0} \leq \theta \leq 2 \pi - \theta_{0}\}$.
Let $\alpha_{2} = \ell_{\theta_{0}, \varepsilon /2} \cup T_{\varepsilon /2} \cup \ell_{2\pi - \theta_0, \varepsilon/2} $ which is considered in Remark \ref{rmk:attachingcircle}.
The fiber of $\pi_{0} |_{\partial M_0}$ over the arc $\alpha_2$ is described as in ($a$) in Figure \ref{fig:precisemove}, since near the endpoints of $\alpha_2$, the absolute value of $x \in \alpha_2$ becomes larger than $\varepsilon$.

Next, move the arc $\alpha_{2}$ closer to $\partial U$.
Define the arc $\alpha_{2}' = \ell_{\theta_0, \varepsilon} \cup T_{\varepsilon} \cup \ell_{2 \pi - \theta_{0}, \varepsilon}$ in the middle, which overlaps with the critical value set $S$.
When $x \in \ell_{\theta_0, \varepsilon} \cup \ell_{2 \pi-\theta_0, \varepsilon} \setminus \{\varepsilon e^{i \theta_{0}}, \varepsilon e^{i (2 \pi -\theta_0)}  \}$, the number of connected components of $(\pi_{0} |_{\partial M_0})^{-1} (x) \cap \partial \nu (C)$ is two.
When $x = \varepsilon e^{i \theta} \in T_{\varepsilon}$, $(\pi_0 |_{\partial M_0})^{-1} (x) \cap \partial \nu (C)$ is the union of two lobes centered at $\pm \varepsilon e^{i \theta}$ which are intersecting at the origin.
When the point $x = \varepsilon e^{i \theta}$ runs along $T_{\varepsilon}$, the centers of these lobes in the fiber also rotate.
Thus, in the fiber $(\pi_0 |_{\partial M_0})^{-1} (\alpha_{2}') \cap \partial \nu (C)$, we see a family of two lobes, whose centers are rotating, as in ($b$) of Figure \ref{fig:precisemove}.

Finally, consider the arc $\alpha_3 = T_{2 \varepsilon}$.
For $x = 2 \varepsilon e^{i \theta} \in \alpha_{3}$, $(\pi_{0} |_{\partial M_0})^{-1} (x) \cap \partial \nu (C)$ is the disjoint union of two lobes centered at $\pm 2\varepsilon e^{i \theta}$.
As the point $x \in \alpha_{3}$ runs along $\alpha_3$, the centers of these lobes rotate, and the trajectory draws the full twist as in ($c$) of Figure \ref{fig:precisemove}.
The centers of these lobes appears as the intersection $\pi^{-1}(x) \cap C$ for each $x \in \alpha_{3} \subset \partial U $, so they coincide with the curves $q_{i,1} (t), q_{i,2} (t)$ which are defined in Remark \ref{rmk:braid} (in this example, $n=m_i =2$). 
Thus, the centers of these lobes form the braid representing the local monodromy.

In general, $\pi^{-1} (\alpha_3) \cap C$ is the disjoint union of $n$ curves written as $q_{i,1}(t), \ldots, q_{i,n} (t)$, which are defined in Remark \ref{rmk:braid}.
They form the braid which represents the local monodromy, as noted in Remark \ref{rmk:braid}. Since $(\pi_0 |_{\partial M_0})^{-1} (\alpha_{3}) \cap \nu (C)$ is the boundary of the neighborhood of these curves, the thickened braid appears as (iv) in Figure \ref{fig:attachingcircle} and ($c$) in Figure \ref{fig:precisemove}.

\end{example}

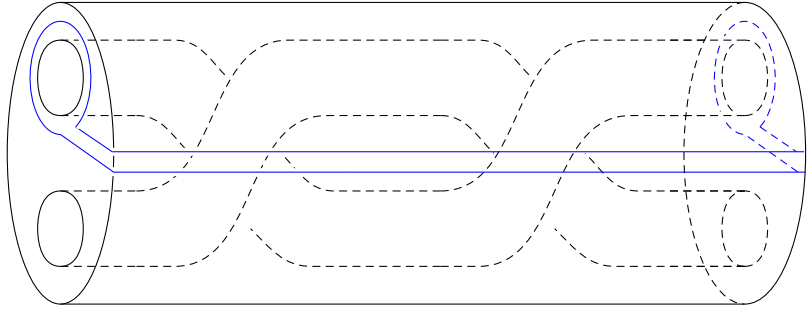
\begin{figure}
\centering
\begin{tikzpicture}
\coordinate (P5) at (0,0);

\draw (P5) ++ (-1.5,0) circle [x radius =0.7, y radius = 2];
\draw (P5) ++ (-1.5,2) --++ (9,0);
\draw (P5) ++ (-1.5,-2) --++ (9,0);

\draw (P5) ++ (7.5,-2) arc [x radius =0.8, y radius = 2, start angle =-90, end angle =90];
\draw [densely dashed](P5) ++ (7.5,2) arc [x radius =0.8, y radius = 2, start angle =90, end angle =270];
%%外円周

\draw [densely dashed](P5)++ (-1.5,1.5) --++ (1,0);
\draw[densely dashed] (P5) ++ (-1.5,0.5) --++(1,0);
\draw[densely dashed] (P5) ++ (7.5,0.5) --++(-1,0);
\draw[densely dashed] (P5) ++ (7.5,1.5) --++(-1,0);
%%上筒
\draw(P5) ++ (-1.5,0.5) [out =0,in=-90] to++ (0.3,0.5) [out=90,in=0] to++ (-0.3,0.5);
\draw (P5) ++ (-1.5,0.5) [out =180,in=-90] to++ (-0.3,0.5) [out=90,in=180] to++ (0.3,0.5);
%%左
\draw[densely dashed] (P5) ++ (-1.5,0.5) ++ (9,0) [out =0,in=-90] to++ (0.3,0.5) [out=90,in=0] to++ (-0.3,0.5);
\draw [densely dashed](P5) ++ (-1.5,0.5)++ (9,0) [out =180,in=-90] to++ (-0.3,0.5) [out=90,in=180] to++ (0.3,0.5);
%%右
%P2ファイバーの上側

%\draw (P5) ++ (-1.5,-1) --++ (3,0);
\draw [densely dashed](P5) ++ (-1.5,-1.5) --++ (1,0);
\draw [densely dashed](P5) ++ (-1.5,-0.5) --++(1,0);
\draw [densely dashed](P5) ++ (7.5,-0.5) --++(-1,0);
\draw [densely dashed](P5) ++ (7.5,-1.5) --++(-1,0);
%下筒
\draw (P5) ++ (-1.5,-1.5) [out =0,in=-90] to++ (0.3,0.5) [out=90,in=0] to++ (-0.3,0.5);
\draw(P5) ++ (-1.5,-1.5) [out =180,in=-90] to++ (-0.3,0.5) [out=90,in=180] to++ (0.3,0.5);
%%左
\draw [densely dashed](P5) ++ (-1.5,-1.5) ++ (9,0) [out =0,in=-90] to++ (0.3,0.5) [out=90,in=0] to++ (-0.3,0.5);
\draw [densely dashed](P5) ++ (-1.5,-1.5)++ (9,0) [out =180,in=-90] to++ (-0.3,0.5) [out=90,in=180] to++ (0.3,0.5);
%%右
%P2ファイバーの下側

\draw [densely dashed](P5)++ (-1.5,1.5) ++ (1,0) --++(0.5,0) to[out=0,in=135] ++(0.7,-0.5);
\draw[densely dashed] (P5) ++ (-1.5,0.5) ++(1,0) to[out=0,in=135] ++(0.7,-0.5);

\draw[densely dashed] (P5) ++ (-1.5,-1.5) ++(5,0) --++(-1.7,0) to[out=180,in=-45] ++(-0.8,0.5);
\draw[densely dashed] (P5) ++ (-1.5,-0.5) ++(5,0) --++(-1.4,0) to[out=180,in=-45] ++(-0.7,0.5);

\draw [densely dashed](P5) ++ (-1.5,-1.5) ++ (1,0) --++(0.5,0) to[out=0,in=180] ++(2,2) --++(1.5,0);
\draw [densely dashed](P5) ++ (-1.5,-0.5) ++(1,0) --++(0,0) to[out=0,in=180] ++(2,2) --++(2,0);
%%前半

\draw [densely dashed](P5) ++(4,0)++ (-1.5,1.5) ++ (1,0) --++(0.5,0) to[out=0,in=135] ++(0.7,-0.5);
\draw[densely dashed] (P5)++(4,0) ++ (-1.5,0.5) ++(1,0) to[out=0,in=135] ++(0.7,-0.5);

\draw[densely dashed] (P5)++(4,0) ++ (-1.5,-1.5) ++(5,0) --++(-1.7,0) to[out=180,in=-45] ++(-0.8,0.5);
\draw[densely dashed] (P5)++(4,0) ++ (-1.5,-0.5) ++(5,0) --++(-1.4,0) to[out=180,in=-45] ++(-0.7,0.5);

\draw [densely dashed](P5) ++(4,0)++ (-1.5,-1.5) ++ (1,0) --++(0.5,0) to[out=0,in=180] ++(2,2) --++(1.5,0);
\draw [densely dashed](P5) ++(4,0)++ (-1.5,-0.5) ++(1,0) --++(0,0) to[out=0,in=180] ++(2,2) --++(2,0);

%%後半
%中央のブレイド

\draw [blue] (P5) ++(-1.3,0.35) --++(0.48,-0.33) --++(9.1,0);
\draw [blue] (P5) ++ (-1.3,0.35) arc [x radius =0.4, y radius =0.75, start angle = -60, end angle =270] --++(0.7,-0.5) --++(9.1,0);
\draw [blue, densely dashed] (P5) ++(-1.3,0.35)++(9,0) --++(0.48,-0.33) ;
\draw [blue, densely dashed] (P5) ++ (-1.3,0.35)++(9,0) arc [x radius =0.4, y radius =0.75, start angle = -60, end angle =270] --++(0.7,-0.5);
%P2attaching circle

%ぬり直しコーナー
\draw[preaction={draw = white, line width = 3 pt}, densely dashed] (P5) ++ (7.5,0.5) --++(-0.5,0);
\draw[preaction={draw = white, line width = 3 pt}, densely dashed] (P5) ++ (7.5,1.5) --++(-0.5,0);
\draw [preaction={draw = white, line width = 3 pt},blue] (P5) ++(-1.3,0.35) ++(0.48,-0.33) --++(3.1,0);
\draw [preaction={draw = white, line width = 3 pt},blue] (P5) ++ (-1.3,0.35) arc [x radius =0.4, y radius =0.75, start angle = -60, end angle =270] --++(0.7,-0.5) --++(3.1,0);
%ぬり直しコーナー

\end{tikzpicture}
\caption{Attaching $2$-handles around $\{y^2 = x^2\}$}
\label{fig:attach22}
\end{figure}

\begin{figure}
\centering
\begin{tikzpicture}
\coordinate (P5) at (0,0);

\draw (P5) ++ (-1.5,0) circle [x radius =0.7, y radius = 3];
\draw (P5) ++ (-1.5,3) --++ (9,0);
\draw (P5) ++ (-1.5,-3) --++ (9,0);

\draw (P5) ++ (7.5,-3) arc [x radius =0.8, y radius = 3, start angle =-90, end angle =90];
\draw [densely dashed](P5) ++ (7.5,3) arc [x radius =0.8, y radius = 3, start angle =90, end angle =270];
%%外円周

\draw [densely dashed](P5)++ (-1.5,1.2) --++ (1,0);
\draw[densely dashed] (P5) ++ (-1.5,2.2) --++(1,0);
\draw[densely dashed] (P5) ++ (7.5,1.2) --++(-0.7,0);
\draw[densely dashed] (P5) ++ (7.5,2.2) --++(-1,0);
%%上筒
\draw(P5) ++ (-1.5,1.2) [out =0,in=-90] to++ (0.3,0.5) [out=90,in=0] to++ (-0.3,0.5);
\draw (P5) ++ (-1.5,1.2) [out =180,in=-90] to++ (-0.3,0.5) [out=90,in=180] to++ (0.3,0.5);
%%左
\draw[densely dashed] (P5) ++ (-1.5,1.2) ++ (9,0) [out =0,in=-90] to++ (0.3,0.5) [out=90,in=0] to++ (-0.3,0.5);
\draw [densely dashed](P5) ++ (-1.5,1.2)++ (9,0) [out =180,in=-90] to++ (-0.3,0.5) [out=90,in=180] to++ (0.3,0.5);
%%右
%P2ファイバーの上側

\draw [densely dashed](P5)++ (-1.5,-0.5) --++ (1,0);
\draw[densely dashed] (P5) ++ (-1.5,0.5) --++(1,0);
\draw[densely dashed] (P5) ++ (7.5,-0.5) --++(-1,0);
\draw[densely dashed] (P5) ++ (7.5,0.5) --++(-1,0);
%%左右筒
\draw(P5) ++ (-1.5,-0.5) [out =0,in=-90] to++ (0.3,0.5) [out=90,in=0] to++ (-0.3,0.5);
\draw (P5) ++ (-1.5,-0.5) [out =180,in=-90] to++ (-0.3,0.5) [out=90,in=180] to++ (0.3,0.5);
%%左
\draw[densely dashed] (P5) ++ (-1.5,-0.5) ++ (9,0) [out =0,in=-90] to++ (0.3,0.5) [out=90,in=0] to++ (-0.3,0.5);
\draw [densely dashed](P5) ++ (-1.5,-0.5)++ (9,0) [out =180,in=-90] to++ (-0.3,0.5) [out=90,in=180] to++ (0.3,0.5);
%%右

%P2ファイバーの中央

%\draw (P5) ++ (-1.5,-1) --++ (3,0);
\draw [densely dashed](P5) ++ (-1.5,-2.2) --++ (1,0);
\draw [densely dashed](P5) ++ (-1.5,-1.2) --++(1,0);
\draw [densely dashed](P5) ++ (7.5,-2.2) --++(-0.7,0);
\draw [densely dashed](P5) ++ (7.5,-1.2) --++(-1,0);
%下筒
\draw (P5) ++ (-1.5,-2.2) [out =0,in=-90] to++ (0.3,0.5) [out=90,in=0] to++ (-0.3,0.5);
\draw(P5) ++ (-1.5,-2.2) [out =180,in=-90] to++ (-0.3,0.5) [out=90,in=180] to++ (0.3,0.5);
%%左
\draw [densely dashed](P5) ++ (-1.5,-2.2) ++ (9,0) [out =0,in=-90] to++ (0.3,0.5) [out=90,in=0] to++ (-0.3,0.5);
\draw [densely dashed](P5) ++ (-1.5,-2.2)++ (9,0) [out =180,in=-90] to++ (-0.3,0.5) [out=90,in=180] to++ (0.3,0.5);
%%右
%P2ファイバーの下側

\draw [densely dashed](P5)++ (-1.5,2.2) ++ (1,0) --++(0.7,0) to[out=0,in=150] ++(1,-0.5);
\draw[densely dashed] (P5) ++ (-1.5,1.2) ++(1,0)--++(0.3,0) to[out=0,in=150] ++(1,-0.5);
\draw[densely dashed] (P5) ++ (-1.5,0.5) ++(5,0) --++(-0.7,0) to[out=180,in=150]++(-0.4,0.2) ;
\draw[densely dashed] (P5) ++ (-1.5,-0.5) ++(5,0) --++(-1,0) to[out=180,in=150]++(-0.6,0.3) ;
%%%1本目

\draw[densely dashed] (P5) ++ (-1.5,0.5) ++ (1,0) --++(0,0) to[out=0,in=150] ++(1,-0.5);

\draw[densely dashed] (P5) ++ (-1.5,-0.5) ++ (1,0) --++(0,0) to[out=0,in=150] ++(0.6,-0.3);

\draw[densely dashed] (P5) ++ (-1.5,-1.2) ++(5,0) --++(-1.3,0) to[out=180,in=150]++(-0.4,0.2) ;
\draw[densely dashed] (P5) ++ (-1.5,-2.2) ++(5,0) --++(-1.6,0) to[out=180,in=150]++(-0.6,0.3) ;
%%%2本目

\draw [densely dashed](P5) ++ (-1.5,-2.2) ++ (1,0) --++(0.5,0) to[out=0,in=180] ++(3,3.4) --++(0.5,0);
\draw [densely dashed](P5) ++ (-1.5,-1.2) ++(1,0) --++(0,0) to[out=0,in=180] ++(3,3.4) --++(1,0);
%%%3本目

%%前半

\draw [densely dashed](P5)++(4,0)++ (-1.5,2.2) ++ (1,0) --++(0.7,0) to[out=0,in=150] ++(1,-0.5);
\draw[densely dashed] (P5)++(4,0) ++ (-1.5,1.2) ++(1,0)--++(0.3,0) to[out=0,in=150] ++(1,-0.5);
\draw[densely dashed] (P5)++(4,0) ++ (-1.5,0.5) ++(5,0) --++(-0.7,0) to[out=180,in=150]++(-0.4,0.2) ;
\draw[densely dashed] (P5) ++(4,0)++ (-1.5,-0.5) ++(5,0) --++(-1,0) to[out=180,in=150]++(-0.6,0.3) ;
%%%1本目

\draw[densely dashed] (P5) ++(4,0)++ (-1.5,0.5) ++ (1,0) --++(0,0) to[out=0,in=150] ++(1,-0.5);

\draw[densely dashed] (P5) ++(4,0)++ (-1.5,-0.5) ++ (1,0) --++(0,0) to[out=0,in=150] ++(0.6,-0.3);

\draw[densely dashed] (P5) ++(4,0)++ (-1.5,-1.2) ++(5,0) --++(-1.3,0) to[out=180,in=150]++(-0.4,0.2) ;
\draw[densely dashed] (P5) ++(4,0)++ (-1.5,-2.2) ++(5,0) --++(-1.6,0) to[out=180,in=150]++(-0.6,0.3) ;
%%%2本目

\draw [densely dashed](P5)++(4,0) ++ (-1.5,-2.2) ++ (1,0) --++(0.5,0) to[out=0,in=180] ++(3,3.4) --++(0.5,0);
\draw [densely dashed](P5)++(4,0) ++ (-1.5,-1.2) ++(1,0) --++(0,0) to[out=0,in=180] ++(3,3.4) --++(1,0);
%%%3本目
%%後半(++(4,0)する)
%中央のブレイド

\draw [blue] (P5)++(0,0.7) ++(-1.3,0.35) --++(0.48,-0.33) --++(9.1,0);
\draw [blue] (P5)++(0,0.7) ++ (-1.3,0.35) arc [x radius =0.4, y radius =0.75, start angle = -60, end angle =270] --++(0.7,-0.5) --++(9.1,0);
\draw [blue, densely dashed] (P5)++(0,0.7) ++(-1.3,0.35)++(9,0) --++(0.48,-0.33) ;
\draw [blue, densely dashed] (P5) ++(0,0.7)++ (-1.3,0.35)++(9,0) arc [x radius =0.4, y radius =0.75, start angle = -60, end angle =270] --++(0.7,-0.5);
%P2attaching circle1

\draw [blue] (P5)++(0,-1) ++(-1.3,0.35) --++(0.48,-0.33) --++(9.1,0);
\draw [blue] (P5)++(0,-1)++ (-1.3,0.35) arc [x radius =0.4, y radius =0.75, start angle = -60, end angle =270] --++(0.7,-0.5) --++(9.1,0);
\draw [blue, densely dashed] (P5)++(0,-1) ++(-1.3,0.35)++(9,0) --++(0.48,-0.33) ;
\draw [blue, densely dashed] (P5) ++(0,-1)++ (-1.3,0.35)++(9,0) arc [x radius =0.4, y radius =0.75, start angle = -60, end angle =270] --++(0.7,-0.5);
%P2attaching circle2

%ぬり直しコーナー
\draw[preaction={draw = white, line width = 3 pt}, densely dashed] (P5) ++ (7.5,0.5) --++(-0.5,0);
%\draw[preaction={draw = white, line width = 3 pt}, densely dashed] (P5) ++ (7.5,1.5) --++(-0.5,0);
\draw [preaction={draw = white, line width = 3 pt},blue] (P5)++(0,0.7) ++(-1.3,0.35) ++(0.48,-0.33) --++(9.1,0);
\draw [preaction={draw = white, line width = 3 pt},blue] (P5)++(0,0.7) ++ (-1.3,0.35) arc [x radius =0.4, y radius =0.75, start angle = -60, end angle =270] --++(0.7,-0.5) --++(9.1,0);
\draw [preaction={draw = white, line width = 3 pt},blue] (P5)++(0,-1) ++(-1.3,0.35) --++(0.48,-0.33) --++(9.1,0);
\draw [preaction={draw = white, line width = 3 pt},blue] (P5)++(0,-1)++ (-1.3,0.35) arc [x radius =0.4, y radius =0.75, start angle = -60, end angle =270] --++(0.7,-0.5) --++(9.1,0);
%ぬり直しコーナー

\end{tikzpicture}
\caption{Attaching $2$-handles around $\{y^3 = x^2\}$}
\label{fig:attach23}
\end{figure}
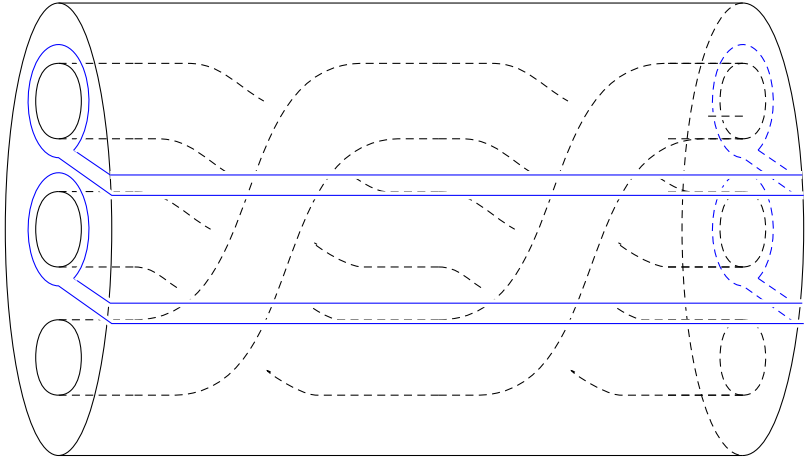

If the local equations are $\{y^2 = x^2\}$ and $\{y^3 = x^2\}$, the attaching circles of $2$-handles are described as in Figure \ref{fig:attach22}, \ref{fig:attach23} respectively.
Later, we will consider $M_{1}$ as the $1$-handlebody of a handle decomposition of $M_{0}$ and use these perturbed attaching circles when we describe the Kirby diagrams.

\section{Kirby diagrams and examples}\label{sec:kirby}
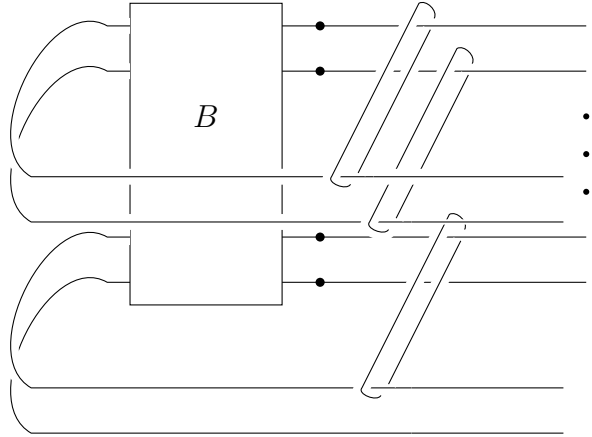
\begin{figure}
\centering

\begin{tikzpicture}
\coordinate (A) at (0,0);
\coordinate (B) at (-2,0);
\coordinate (C) at (6,0);

\draw (A) --++ (4,0);
\draw (A)++(0,-0.6) --++ (4,0);
\draw (A)++ (0,-2.8) --++ (4,0);
\draw (A)++ (0,-3.4) --++ (4,0);
\fill (A) ++ (4,-1.2) circle (0.04);
\fill (A) ++ (4,-1.7) circle (0.04);
\fill (A) ++ (4,-2.2) circle (0.04);

%手前側

\draw (A)++(0,0.3) --++ (0,-4) --++ (-2,0) --++(0,4)-- cycle;
%箱B_{p,q}

\draw (A) ++ (-1,-1.2) node{$B$};
%B_{p,q}名前。Bの上が右を向くように時計回りに90回転させる。

\draw (B) ++ (0,-3.4) --++ (-0.3,0) to[out=150,in=150] ++(-1,-2) --++ (5,0);
\draw [preaction={draw=white,line width=6pt} ](B) ++ (0,-2.8) --++ (-0.3,0) to[out=150,in=150] ++(-1,-2) --++ (5,0);
\draw [preaction={draw=white,line width=6pt} ](B) ++ (0,-0.6) --++ (-0.3,0) to[out=150,in=150] ++(-1,-2) --++ (5,0);
\draw [preaction={draw=white,line width=6pt} ](B)  --++ (-0.3,0) to[out=150,in=150] ++(-1,-2)--++ (4.5,0);
%上下つなぐやつ

%\draw (A) ++ (1.7,0.5) --++ (-1.5,-3);
\draw [preaction={draw=white,line width=6pt} ](A) ++ (2,0.1) to[out=60,in=30] ++(-0.2,0.2) --++ (-1.15,-2.3) to[out=-130, in =-120] ++ (0.25,-0.1);
\draw [preaction={draw=white,line width=6pt} ](A) ++ (1.95,-0.05) --++ (-0.95,-1.9);

\draw [preaction={draw=white,line width=6pt} ](A) ++ (2.5,-0.5) to[out=60,in=30] ++(-0.2,0.2) --++ (-1.15,-2.3) to[out=-130, in =-120] ++ (0.25,-0.1);
\draw [preaction={draw=white,line width=6pt} ](A) ++ (2.45,-0.65) --++ (-0.95,-1.9);

\draw [preaction={draw=white,line width=6pt} ](A) ++(-0.3,0)++ (2.7,-2.7) to[out=60,in=30] ++(-0.2,0.2) --++ (-1.15,-2.3) to[out=-130, in =-120] ++ (0.25,-0.1);
\draw [preaction={draw=white,line width=6pt} ](A) ++(-0.3,0)++ (2.65,-2.85) --++ (-0.95,-1.9);
%2-ハンドル

%\draw (C)++(0,0.3) --++ (0,-3.4) --++ (2,0) --++(0,3.4)-- cycle;
%右の箱B_{p,q}

%\draw (C) ++ (1,-1.2) node{$B^{-1}_{p,q}$};
\draw [preaction={draw=white,line width=6pt} ](B) ++(0,-3.4) ++ (-0.3,0) ++(-1,-2) ++ (5,0) --++(2,0);
\draw [preaction={draw=white,line width=6pt} ](B) ++(0,-2.8) ++ (-0.3,0) ++(-1,-2) ++ (5,0) --++(2,0);
\draw [preaction={draw=white,line width=2pt} ](B)  ++(0,-0.6)++ (-0.3,0) ++(-1,-2) ++ (5,0) --++(2,0);
\draw [preaction={draw=white,line width=6pt} ](B)  ++ (-0.3,0) ++(-1,-2) ++ (4.5,0) --++(2.5,0);
%下の線続き+右側

\fill (A) ++ (0.5,0) circle (0.06);
\fill (A) ++ (0.5,-0.6) circle (0.06);
\fill (A) ++ (0.5,-2.8) circle (0.06);
\fill (A) ++ (0.5,-3.4) circle (0.06);
%点付き円の点

%\draw (A) ++ (0.2,0) node[above]{$1$};
%\draw (A) ++ (0.2,-0.6) node[above]{$2$};
%\draw (A) ++ (0.2,-2.8) node[below]{$q$};
%番号

\end{tikzpicture}
\caption{Local description of the Kirby diagram. There are $n$ $1$-handles and $m_{i}-1$ $2$-handles and $B$ is the braid representing the local monodromy.}
\label{fig:localkirby}
\end{figure}

In this section, we see how to obtain Kirby diagrams and give some examples.
First, a $1$-handle is presented as a dotted circle in a Kirby diagram.
Removing the tubular neighborhood of a boundary-parallel and properly embedded $D^2$ in a $4$-ball $B^4$ is equivalent to attaching a $1$-handle. 
The dotted circle of the $1$-handle obtained in this way is the boundary of the embedded disk in $S^3$.
The $1$-handlebody $M_1$ of $M_0$ is obtained from $\pi^{-1} (D \setminus \nu(\Gamma)) \cong B^4$ by removing $\pi^{-1} (D \setminus \nu(\Gamma)) \cap \nu (C)$. 
The Kirby diagram is drawn in $\partial ( \pi^{-1} (D \setminus \nu (\Gamma))) \cong S^3$.
Since the removed part $\pi^{-1} (D \setminus \nu(\Gamma)) \cap \nu (C)$ is a disjoint union of tubular neighborhood of $n$ properly embedded $D^2$, the 
dotted circles of these $1$-handles are $\pi^{-1} (\partial (D \setminus \nu (\Gamma))) \cap C$ (which is a trivial link of $n$-components).
As noted in Remark \ref{rmk:braid}, around the fiber of tangent points and singularities of $C$, $\pi^{-1} (\partial U_{i}) \cap C$ is a braid which represents the local monodromy, (and, since $\nu (s_i)$ is small, $\pi^{-1} (\partial U_i \setminus \nu (s_{i})) \cap C$ forms the same braid). 
Since $\partial U_{i} \setminus \nu (s_{i}) \subset \partial (D \setminus \nu (\Gamma))$, this braid also appears as the dotted circle of the $1$-handles.
In Remark \ref{rmk:attachingcircle}, we moved the attaching circle of the $2$-handle so that it is on $\partial M_{1}$ and this final step implies that the braid $\pi^{-1}(\partial U_i \setminus \nu (s_i)) \cap C$ and the attaching circle links (the attaching circle looks like the ``meridian loop'' of the braid locally).
Thus, the Kirby diagram is locally described as in Figure \ref{fig:localkirby}.
To describe the global Kirby diagram, we need the relationship among local monodromies. These are determined by the trivialization along each path $s_{i}$. The braid defined by the trivialization appears as the dotted circles of $1$-handles, and this connects each local monodromy.

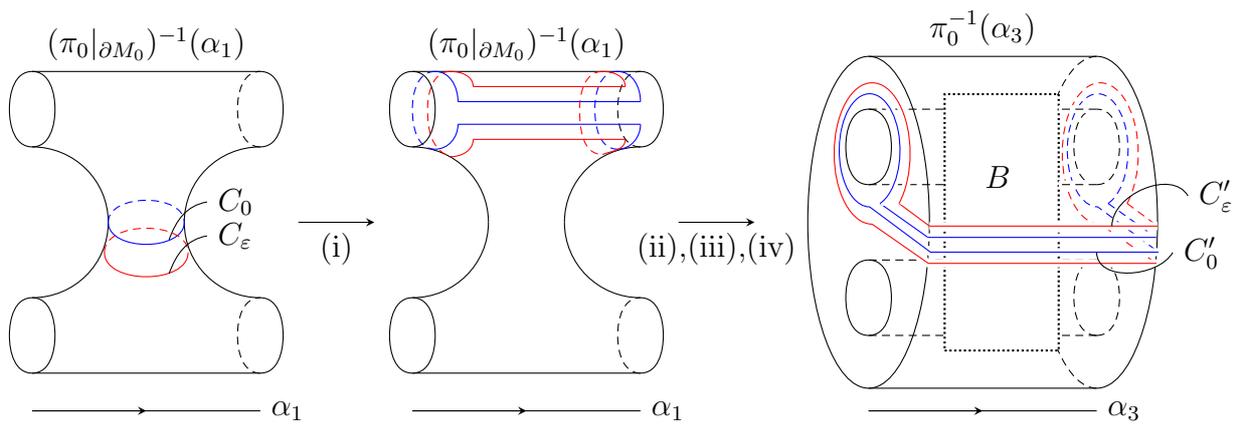
\begin{figure}
\centering
\begin{tikzpicture}

\coordinate (P1) at (-10,0);
\coordinate (P2) at (-5,0);
\coordinate (P5) at (1,0);

\draw (P1) ++ (-1.5,2) --++ (3,0);
%\draw (P1) ++ (-1.5,1) --++ (3,0);
\draw (P1) ++ (-1.5,1) [out=0,in=90] to++(1,-1);
\draw (P1) ++ (1.5,1) [out=180,in=90] to++(-1,-1);
%%上筒
\draw (P1) ++ (-1.5,1) [out =0,in=-90] to++ (0.3,0.5) [out=90,in=0] to++ (-0.3,0.5);
\draw (P1) ++ (-1.5,1) [out =180,in=-90] to++ (-0.3,0.5) [out=90,in=180] to++ (0.3,0.5);
%%左
\draw (P1) ++ (-1.5,1) ++ (3,0) [out =0,in=-90] to++ (0.3,0.5) [out=90,in=0] to++ (-0.3,0.5);
\draw [densely dashed](P1) ++ (-1.5,1)++ (3,0) [out =180,in=-90] to++ (-0.3,0.5) [out=90,in=180] to++ (0.3,0.5);
%%右
%P1ファイバーの上側

%\draw (P1) ++ (-1.5,-1) --++ (3,0);
\draw (P1) ++ (-1.5,-2) --++ (3,0);
\draw (P1) ++ (-1.5,-1) [out=0,in=-90] to++(1,1);
\draw (P1) ++ (1.5,-1) [out=180,in=-90] to++(-1,1);
%%下筒
\draw (P1) ++ (-1.5,-2) [out =0,in=-90] to++ (0.3,0.5) [out=90,in=0] to++ (-0.3,0.5);
\draw (P1) ++ (-1.5,-2) [out =180,in=-90] to++ (-0.3,0.5) [out=90,in=180] to++ (0.3,0.5);
%%左
\draw (P1) ++ (-1.5,-2) ++ (3,0) [out =0,in=-90] to++ (0.3,0.5) [out=90,in=0] to++ (-0.3,0.5);
\draw [densely dashed](P1) ++ (-1.5,-2)++ (3,0) [out =180,in=-90] to++ (-0.3,0.5) [out=90,in=180] to++ (0.3,0.5);
%%右
%P1ファイバーの下側

\draw[blue] (P1) ++(-0.5,0) to[out=-90,in=-90] ++(1,0);
\draw[densely dashed, blue] (P1) ++(-0.5,0) to[out=90,in=90] ++(1,0);
\draw (P1) ++ (0.3,-0.23) to[out=70,in=-170] ++(0.5,0.5) node[right]{$C_{0}$};
%P1attaching circle

\draw[red] (P1) ++(-0.55,-0.4) to[out=-90,in=-90] ++(1.1,0);
\draw[densely dashed, red] (P1) ++(-0.55,-0.4) to[out=90,in=90] ++(1.1,0);
\draw (P1) ++ (0.3,-0.68) to[out=70,in=-170] ++(0.5,0.5) node[right]{$C_{\varepsilon}$};
%Perturbed attaching circle

\draw (P1) ++(0,2)node[above]{$(\pi_{0} |_{\partial M_{0}})^{-1} (\alpha_{1})$};

\draw[->,>=stealth] (P1) ++ (-1.5,-2.5) --++(1.5,0);
\draw (P1) ++ (-1.5,-2.5) --++(3,0) node[right]{$\alpha_{1}$};

\draw[->,>=stealth] (P1) ++ (2,0) --++(1,0);
\draw (P1) ++(2.5,0) node[below]{(i)};

%%%%↑ここまで(P1)、(Q1)%%%%

%%%%↓ここから(P2)、(Q2)%%%%
\draw (P2) ++ (-1.5,2) --++ (3,0);
%\draw (P2) ++ (-1.5,1) --++ (3,0);
\draw (P2) ++ (-1.5,1) [out=0,in=90] to++(1,-1);
\draw (P2) ++ (1.5,1) [out=180,in=90] to++(-1,-1);
%%上筒
\draw (P2) ++ (-1.5,1) [out =0,in=-90] to++ (0.3,0.5) [out=90,in=0] to++ (-0.3,0.5);
\draw (P2) ++ (-1.5,1) [out =180,in=-90] to++ (-0.3,0.5) [out=90,in=180] to++ (0.3,0.5);
%%左
\draw (P2) ++ (-1.5,1) ++ (3,0) [out =0,in=-90] to++ (0.3,0.5) [out=90,in=0] to++ (-0.3,0.5);
\draw [densely dashed](P2) ++ (-1.5,1)++ (3,0) [out =180,in=-90] to++ (-0.3,0.5) [out=90,in=180] to++ (0.3,0.5);
%%右
%P2ファイバーの上側

%\draw (P2) ++ (-1.5,-1) --++ (3,0);
\draw (P2) ++ (-1.5,-2) --++ (3,0);
\draw (P2) ++ (-1.5,-1) [out=0,in=-90] to++(1,1);
\draw (P2) ++ (1.5,-1) [out=180,in=-90] to++(-1,1);
%%下筒
\draw (P2) ++ (-1.5,-2) [out =0,in=-90] to++ (0.3,0.5) [out=90,in=0] to++ (-0.3,0.5);
\draw (P2) ++ (-1.5,-2) [out =180,in=-90] to++ (-0.3,0.5) [out=90,in=180] to++ (0.3,0.5);
%%左
\draw (P2) ++ (-1.5,-2) ++ (3,0) [out =0,in=-90] to++ (0.3,0.5) [out=90,in=0] to++ (-0.3,0.5);
\draw [densely dashed](P2) ++ (-1.5,-2)++ (3,0) [out =180,in=-90] to++ (-0.3,0.5) [out=90,in=180] to++ (0.3,0.5);
%%右
%P2ファイバーの下側

\draw [blue] (P2) ++ (-0.9,1.6) to[out=90,in=0]++(-0.3,0.4);
\draw [blue] (P2) ++ (-0.9,1.3) to[out=-90,in=0]++(-0.3,-0.33);
\draw [blue, densely dashed] (P2) ++ (-1.2,2) to[out=180,in=90] ++(-0.3,-0.5) to[out=-90,in=180] ++(0.3,-0.53);
%%左側のからみ

\draw [blue] (P2) ++(-0.9,1.6) --++(2.4,0);
\draw [blue] (P2) ++(-0.9,1.3) --++(2.4,0);
%%中央

\draw [blue] (P2) ++ (1.5,1.6) to[out=90,in=0]++(-0.3,0.4);
\draw [blue] (P2) ++ (1.5,1.3) to[out=-90,in=0]++(-0.3,-0.33);
\draw [blue, densely dashed] (P2) ++ (1.2,2) to[out=180,in=90] ++(-0.3,-0.5) to[out=-90,in=180] ++(0.3,-0.53);
%%左側のからみ
%P2attaching circle

\draw [red] (P2) ++ (-0.7,1.8) to[out=90,in=0]++(-0.3,0.2);
\draw [red] (P2) ++ (-0.7,1.1) to[out=-90,in=0]++(-0.3,-0.23);
\draw [red, densely dashed] (P2) ++ (-1.0,2) to[out=180,in=90] ++(-0.3,-0.5) to[out=-90,in=160] ++(0.3,-0.63);
%%左側のからみ

\draw [red] (P2) ++(-0.7,1.8) --++(2.0,0);
\draw [red] (P2) ++(-0.7,1.1) --++(2.0,0);
%%中央

\draw [red] (P2) ++ (1.3,1.8) to[out=90,in=0]++(-0.3,0.2);
\draw [red] (P2) ++ (1.3,1.1) to[out=-90,in=20]++(-0.3,-0.2);
\draw [red, densely dashed] (P2) ++ (1.0,2) to[out=180,in=90] ++(-0.3,-0.5) to[out=-90,in=180] ++(0.3,-0.6);
%%左側のからみ
%Perturbed attaching circle

\draw (P2) ++(0,2)node[above]{$(\pi_{0} |_{\partial M_{0}})^{-1} (\alpha_{1})$};

\draw[->,>=stealth] (P2) ++ (-1.5,-2.5) --++(1.5,0);
\draw (P2) ++ (-1.5,-2.5) --++(3,0) node[right]{$\alpha_{1}$};

\draw[->,>=stealth] (P2) ++ (2,0) --++(1,0);
\draw (P2) ++(2.5,0) node[below]{(ii),(iii),(iv)};

%%%%↑ここまで(P2)、(Q2)%%%%

%%%%↓ここから(P5)、(Q5)%%%%
\draw (P5) ++ (-1.5,0) circle [x radius =0.8, y radius = 2.2];
\draw (P5) ++ (-1.5,2.2) --++ (3,0);
\draw (P5) ++ (-1.5,-2.2) --++ (3,0);

\draw (P5) ++ (1.5,-2.2) arc [x radius =0.8, y radius = 2.2, start angle =-90, end angle =90];
\draw [densely dashed](P5) ++ (1.5,2.2) arc [x radius =0.8, y radius = 2.2, start angle =90, end angle =130];
\draw [densely dashed](P5) ++ (1.5,-2.2) arc [x radius =0.8, y radius = 2.2, start angle =-90, end angle =-130];
%%外円周

\draw [densely dashed](P5)++ (-1.5,1.5) --++ (1,0);
\draw[densely dashed] (P5) ++ (-1.5,0.5) --++(1,0);
\draw[densely dashed] (P5) ++ (1.5,0.5) --++(-0.5,0);
\draw[densely dashed] (P5) ++ (1.5,1.5) --++(-0.5,0);
%%上筒
\draw(P5) ++ (-1.5,0.5) [out =0,in=-90] to++ (0.3,0.5) [out=90,in=0] to++ (-0.3,0.5);
\draw (P5) ++ (-1.5,0.5) [out =180,in=-90] to++ (-0.3,0.5) [out=90,in=180] to++ (0.3,0.5);
%%左
\draw[densely dashed] (P5) ++ (-1.5,0.5) ++ (3,0) [out =0,in=-90] to++ (0.3,0.5) [out=90,in=0] to++ (-0.3,0.5);
\draw [densely dashed](P5) ++ (-1.5,0.5)++ (3,0) [out =180,in=-90] to++ (-0.3,0.5) [out=90,in=180] to++ (0.3,0.5);
%%右
%P2ファイバーの上側

%\draw (P5) ++ (-1.5,-1) --++ (3,0);
\draw [densely dashed](P5) ++ (-1.5,-1.5) --++ (1,0);
\draw [densely dashed](P5) ++ (-1.5,-0.5) --++(1,0);
\draw [densely dashed](P5) ++ (1.5,-0.5) --++(-0.5,0);
\draw [densely dashed](P5) ++ (1.5,-1.5) --++(-0.5,0);
%下筒
\draw (P5) ++ (-1.5,-1.5) [out =0,in=-90] to++ (0.3,0.5) [out=90,in=0] to++ (-0.3,0.5);
\draw(P5) ++ (-1.5,-1.5) [out =180,in=-90] to++ (-0.3,0.5) [out=90,in=180] to++ (0.3,0.5);
%%左
\draw [densely dashed](P5) ++ (-1.5,-1.5) ++ (3,0) [out =0,in=-90] to++ (0.3,0.5) [out=90,in=0] to++ (-0.3,0.5);
\draw [densely dashed](P5) ++ (-1.5,-1.5)++ (3,0) [out =180,in=-90] to++ (-0.3,0.5) [out=90,in=180] to++ (0.3,0.5);
%%右
%P2ファイバーの下側

\draw[densely dotted, thick] (P5) ++(-0.5,1.7) --++(0,-3.4) --++ (1.5,0) --++(0,3.4) --cycle;
\draw (P5) ++(0.2,0.6) node{$B$};
%ブレイドの箱

\draw [blue] (P5) ++(-1.35,0.25) --++(0.63,-0.45) --++(3,0);
\draw [blue] (P5) ++ (-1.35,0.25) arc [x radius =0.4, y radius =0.75, start angle = -70, end angle =270] --++(0.8,-0.6) --++(3,0);
\draw [blue, densely dashed] (P5) ++(-1.3,0.25)++(2.95,0) --++(0.63,-0.45) ;
\draw [blue, densely dashed] (P5) ++ (-1.3,0.25)++(2.95,0) arc [x radius =0.4, y radius =0.75, start angle = -70, end angle =270] --++(0.8,-0.6);
%P2attaching circle

\draw [red] (P5) ++(-1.1,0.25) --++(0.4,-0.3) --++(3,0);
\draw [red] (P5) ++ (-1.1,0.25) arc [x radius =0.5, y radius =0.9375, start angle = -45, end angle =270] --++(0.73,-0.52) --++(3,0);
\draw [red, densely dashed] (P5) ++(-1.1,0.25)++(3,0) --++(0.4,-0.3) ;
\draw [red, densely dashed] (P5) ++ (-1.1,0.25)++(3,0) arc [x radius =0.5, y radius =0.9375, start angle = -45, end angle =270] --++(0.74,-0.52);
%Perturbed attaching circle

%ぬり直しコーナー
\draw[preaction={draw = white, line width = 3 pt}, densely dashed] (P5) ++ (1.5,0.5) --++(-0.5,0);
\draw[preaction={draw = white, line width = 3 pt}, densely dashed] (P5) ++ (1.5,1.5) --++(-0.5,0);
\draw [preaction={draw = white, line width = 3 pt},blue] (P5) ++(-1.35,0.25) ++(0.63,-0.45) --++(3,0);
\draw [preaction={draw = white, line width = 3 pt},blue] (P5) ++ (-1.35,0.25) arc [x radius =0.4, y radius =0.75, start angle = -70, end angle =270] --++(0.8,-0.6) --++(3,0);
\draw [preaction={draw = white, line width = 3 pt},red] (P5) ++(-1.1,0.25) ++(0.4,-0.3) --++(3,0);
\draw [preaction={draw = white, line width = 3 pt},red] (P5) ++ (-1.1,0.25) arc [x radius =0.5, y radius =0.9375, start angle = -45, end angle =270] --++(0.74,-0.52) --++(3,0);
%ぬり直しコーナー

\draw (P5) ++(1.5,-0.4) to[out=-90,in=-120] ++(1,0)node[right]{$C'_{0}$};
\draw (P5) ++(1.7,-0.05) to[out=90,in=120] ++(1,0.4)node[right]{$C'_{\varepsilon}$};

\draw (P5) ++(0,2.2)node[above]{$\pi_{0}^{-1} (\alpha_{3})$};

\draw[->,>=stealth] (P5) ++ (-1.5,-2.5) --++(1.5,0);
\draw (P5) ++ (-1.5,-2.5) --++(3,0) node[right]{$\alpha_{3}$};

\end{tikzpicture}
\caption{Perturbation of the attaching circle. The procedures (i), (ii), (iii), and (iv) are the same as Figure \ref{fig:attachingcircle}.}
\label{fig:perturb}
\end{figure}

Next, we have to consider the framing coefficient of each $2$-handle.
The attaching circle of the $2$-handle is written as $C_{0} = \{(X_0, Y, r_{0}, v) \mid Y^2 + v^2 = X_{0}^{2} + \eta\} \subset N_{3} \subset N_{2}$ for appropriately fixed $X_{0}$ and $r_{0}$. 
A slightly perturbed circle is expressed as $C_{\varepsilon} = \{(X_0 + \varepsilon, Y, r_{0}, v) \mid Y^2 + v^2 = (X_{0}+\varepsilon)^{2} + \eta \}$ for sufficiently small $\varepsilon > 0$.
%Since the $2$-handle $N_{2}$ is attached to $\partial M_{2}$ trivially, two circles $C_{0}$ and $C_{\varepsilon}$ do not link. 
In Remark \ref{rmk:attachingcircle}, we move the attaching circle into $\partial M_{1}$ and let $C'_{0}$ be the perturbed attaching circle here.
Similarly, we can move the perturbed circle $C_{\varepsilon}$ into $\partial M_{1}$ and denote it by $C'_{\varepsilon}$ (see Figure \ref{fig:perturb}). 
The perturbed circles $C'_{0}$ and $C'_{\varepsilon}$ are both located in $\partial \left ( \pi^{-1}_{0} (\alpha_{3}) \cup (\pi^{-1} (\alpha_{3}) \cap \nu(C) ) \right ) \cong \partial B^3 \cong S^2$. This $2$-sphere is contained in $\pi^{-1} (\partial (D \setminus \nu (\Gamma))) \cong S^3$, in which the Kirby diagram is drawn.
Therefore, $C_{0}$ and $C_{\varepsilon}$ are located in the same $2$-sphere embedded in $3$-sphere. 
Thus, we can separate $C_{0}$ and $C_{\varepsilon}$, and the linking number is equal to $0$.
Therefore, the framing coefficient is equal to $0$ for every $2$-handle.

To summarize, we have the following:
\begin{theorem}
The dotted circles of the $1$-handles in the Kirby diagram of $M_{0}$ is $\pi^{-1} (\partial (D \setminus \nu (\Gamma))) \cap C$, and the attaching circles of $2$-handles are locally described as in Figure \ref{fig:localkirby} around each $p_{i}$.
In particular, each attaching circle of a $2$-handle is a $0$-framed trivial knot and does not link any other attaching circle.
\end{theorem}

Finally, let us conclude this paper with some examples.

\begin{example} (Brieskorn polynomials)
Let $f (x,y) = x^{p}-y^{q}$ ($2 \leq p \leq q$). Then $X = \{0\}$. Take a graph $\Gamma$ as in Figure \ref{fig:gamma1}. It suffices to consider the monodromy around $x = 0$. When $x$ goes counterclockwise around $0$, the $y$-coordinate in the fiber rotates by $2p\pi/ q$. 
This monodromy can be presented by the following braid $B_{p,q} = (\sigma_{1} \sigma_{2} \cdots \sigma_{q-1})^{p}$ as in Figure \ref{fig:braidpq}, where $\sigma_{j}$ is a canonical generator, presenting the half twist of the $j$-th point and the $(j+1)$-st point  (see also \cite{oka} for the presentation of the fundamental group of the complement of the curve $x^p = y^q$).
Therefore, the Kirby diagram of the complement is described as in Figure \ref{fig:kirby1}.
There are $q$ $1$-handles and $(q-1)$ 2-handles.
For example, let us suppose $p=q=2$.
Then, the Kirby diagram is as in Figure \ref{fig:kirby1'}, so we can see that it is diffeomorphic to $T^2 \times D^2$.
This proves the well-known fact that the complement of $\{x^2-y^2 = 0\}$ is diffeomorphic to the interior of $T^2 \times D^2$.

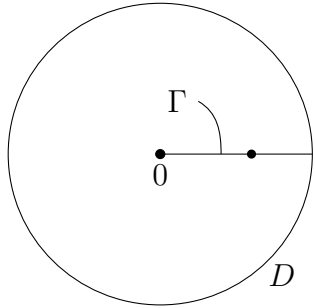
\begin{figure}[htbp]
\centering
\begin{tikzpicture}
\draw (0,0) circle (2);
\draw (1.6,-1.6) node{$D$};
%D
\filldraw (0,0) circle (0.06) node[below]{$0$};
\fill (1.2,0) circle (0.06);
\draw (0,0) -- (2,0);
%グラフ\Gamma

\draw (0.8,0) to[out=90,in=-30] ++ (-0.3,0.7) node[left]{$\Gamma$};
%\Gamma名前

\end{tikzpicture}\caption{A graph $\Gamma$ for $f(x,y) = x^p-y^q$}
\label{fig:gamma1}
\end{figure}

%\vspace{3cm}

\begin{figure}[htbp]
\centering
\begin{tikzpicture}
\coordinate (P) at (0,0);
\coordinate (Q) at (4.5,0);

\draw (P) --++ (3,0) --++ (0,-3) --++ (-3,0) --cycle;
\draw (P) ++ (1.5,-1.5) node{$B_{p,q}$};
\draw (P) ++ (0.3,0.3)node[above]{$1$} --++(0,-0.3);
\draw (P) ++ (0.7,0.3)node[above]{$2$} --++(0,-0.3);
\draw (P) ++ (1.1,0.3)node[above]{$3$} --++(0,-0.3);
\draw (P) ++ (2.7,0.3)node[above]{$q$} --++(0,-0.3);
\fill (P) ++ (1.5,0.15) circle (0.03);
\fill (P) ++ (1.9,0.15) circle (0.03);
\fill (P) ++ (2.3,0.15) circle (0.03);
\draw (P) ++(0,-3.3)++ (0.3,0.3) --++(0,-0.3);
\draw (P) ++(0,-3.3)++ (0.7,0.3) --++(0,-0.3);
\draw (P) ++(0,-3.3)++ (1.1,0.3) --++(0,-0.3);
\draw (P) ++(0,-3.3)++ (2.7,0.3) --++(0,-0.3);
\fill (P) ++(0,-3.3)++ (1.5,0.15) circle (0.03);
\fill (P) ++(0,-3.3)++ (1.9,0.15) circle (0.03);
\fill (P) ++(0,-3.3)++ (2.3,0.15) circle (0.03);
%B_{p,q}ハコ

\draw (P) ++ (3.75,-1.5) node{$:=$};

\fill (Q) ++ (1.5,0.15) circle (0.03);
\fill (Q) ++ (1.9,0.15) circle (0.03);
\fill (Q) ++ (2.3,0.15) circle (0.03);
\draw (Q) ++ (0.3,0.3)node[above]{$1$} --++(0,-0.3);
\draw (Q) ++ (0.7,0.3)node[above]{$2$} --++(0,-0.3);
\draw (Q) ++ (1.1,0.3)node[above]{$3$} --++(0,-0.3);
\draw (Q) ++ (2.7,0.3)node[above]{$q$} --++(0,-0.3);
%上部分
\draw (Q) ++ (0.7,0) to[out=-90,in=90] ++ (-0.4,-0.8);
\draw (Q) ++ (1.1,0) to[out=-90,in=90] ++ (-0.4,-0.8);
\draw (Q) ++ (2.7,0) to[out=-90,in=90] ++ (-0.4,-0.8);
\draw[preaction={draw=white,line width=6pt} ] (Q) ++(0.3,0) to[out=-90,in=90] ++ (2.4,-0.8);
%1回目の捻り

\draw (Q) ++ (0.7,-2.2) to[out=-90,in=90] ++ (-0.4,-0.8);
\draw (Q) ++ (1.1,-2.2) to[out=-90,in=90] ++ (-0.4,-0.8);
\draw (Q) ++ (2.7,-2.2) to[out=-90,in=90] ++ (-0.4,-0.8);
\draw[preaction={draw=white,line width=6pt} ] (Q) ++(0.3,-2.2) to[out=-90,in=90] ++ (2.4,-0.8);
%最後のひねり

\fill (Q) ++ (1.5,-1.1) circle (0.03);
\fill (Q) ++ (1.5,-1.4) circle (0.03);
\fill (Q) ++ (1.5,-1.7) circle (0.03);
%真ん中てんてん

\draw (Q) ++(0,-3.3)++ (0.3,0.3) --++(0,-0.3);
\draw (Q) ++(0,-3.3)++ (0.7,0.3) --++(0,-0.3);
\draw (Q) ++(0,-3.3)++ (2.3,0.3) --++(0,-0.3);
\draw (Q) ++(0,-3.3)++ (2.7,0.3) --++(0,-0.3);
\fill (Q) ++(0,-3.3)++ (1.5,0.15) circle (0.03);
\fill (Q) ++(0,-3.3)++ (1.9,0.15) circle (0.03);
\fill (Q) ++(0,-3.3)++ (1.1,0.15) circle (0.03);
%下部分

\draw (Q) ++ (3,0) to[out=0,in=180] ++(0.5,-1.5);
\draw (Q) ++ (3,-3) to[out=0,in=180] ++(0.5,1.5);
\draw (Q) ++ (3,0) ++ (0.5,-1.5) node[right] {$p$ times};

\end{tikzpicture}
\caption{A braid $B_{p,q}$}
\label{fig:braidpq}
\end{figure}
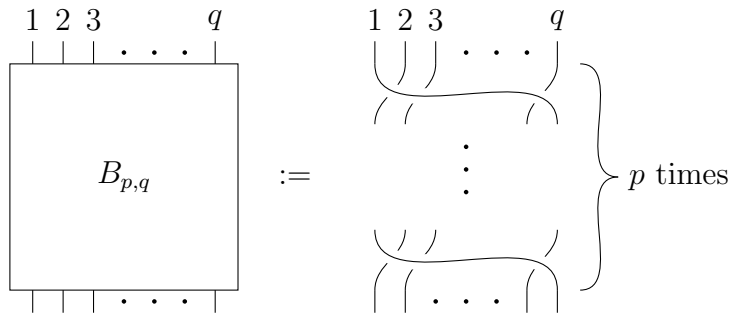

%\vspace{3cm}

\begin{figure}[htbp]
\centering
\begin{tikzpicture}
\coordinate (A) at (0,0);
\coordinate (B) at (-2,0);
\coordinate (C) at (6,0);

\draw (A) --++ (6,0);
\draw (A)++(0,-0.6) --++ (6,0);
\draw (A)++ (0,-2.8) --++ (6,0);
\fill (A) ++ (4,-1.2) circle (0.04);
\fill (A) ++ (4,-1.7) circle (0.04);
\fill (A) ++ (4,-2.2) circle (0.04);
%手前側

\draw (A)++(0,0.3) --++ (0,-3.4) --++ (-2,0) --++(0,3.4)-- cycle;
%箱B_{p,q}

\draw (A) ++ (-1,-1.2) node{\rotatebox{-90}{$B_{p,q}$}};
%B_{p,q}名前。Bの上が右を向くように時計回りに90回転させる。

\draw (B) ++ (0,-2.8) --++ (-0.3,0) to[out=150,in=150] ++(-1,-2) --++ (5,0);
\draw [preaction={draw=white,line width=6pt} ](B) ++ (0,-0.6) --++ (-0.3,0) to[out=150,in=150] ++(-1,-2) --++ (5,0);
\draw [preaction={draw=white,line width=6pt} ](B)  --++ (-0.3,0) to[out=150,in=150] ++(-1,-2)--++ (4.5,0);
%上下つなぐやつ

%\draw (A) ++ (1.7,0.5) --++ (-1.5,-3);
\draw [preaction={draw=white,line width=6pt} ](A) ++ (2,0.1) to[out=60,in=30] ++(-0.2,0.2) --++ (-1.15,-2.3) to[out=-130, in =-120] ++ (0.25,-0.1);
\draw [preaction={draw=white,line width=6pt} ](A) ++ (1.95,-0.05) --++ (-0.95,-1.9);

\draw [preaction={draw=white,line width=6pt} ](A) ++ (2.5,-0.5) to[out=60,in=30] ++(-0.2,0.2) --++ (-1.15,-2.3) to[out=-130, in =-120] ++ (0.25,-0.1);
\draw [preaction={draw=white,line width=6pt} ](A) ++ (2.45,-0.65) --++ (-0.95,-1.9);

%\draw [preaction={draw=white,line width=6pt} ](A) ++ (2.7,-2.7) to[out=60,in=30] ++(-0.2,0.2) --++ (-1.15,-2.3) to[out=-130, in =-120] ++ (0.25,-0.1);
%\draw [preaction={draw=white,line width=6pt} ](A) ++ (2.65,-2.85) --++ (-0.95,-1.9);
%2-ハンドル

\draw (C)++(0,0.3) --++ (0,-3.4) --++ (2,0) --++(0,3.4)-- cycle;
%右の箱B_{p,q}

\draw (C) ++ (1,-1.2) node{\rotatebox{-90}{$B^{-1}_{p,q}$}};
%B_{p,q}名前。回転させる。Bの上が右を向くように時計回りに90回転させる。

\draw [preaction={draw=white,line width=6pt} ](B) ++(0,-2.8) ++ (-0.3,0) ++(-1,-2) ++ (5,0) --++(6,0)to[out=30,in=-30] ++(1,2) --++ (-0.7,0);
\draw [preaction={draw=white,line width=2pt} ](B)  ++(0,-0.6)++ (-0.3,0) ++(-1,-2) ++ (5,0) --++(6,0)to[out=30,in=-30] ++(1,2) --++ (-0.7,0);
\draw [preaction={draw=white,line width=6pt} ](B)  ++ (-0.3,0) ++(-1,-2) ++ (4.5,0) --++(6.5,0) to[out=30,in=-30] ++(1,2) --++ (-0.7,0);
%下の線続き+右側

\fill (A) ++ (0.5,0) circle (0.06);
\fill (A) ++ (0.5,-0.6) circle (0.06);
\fill (A) ++ (0.5,-2.8) circle (0.06);
%点付き円の点

\draw (A) ++ (0.2,0) node[above]{$1$};
\draw (A) ++ (0.2,-0.6) node[above]{$2$};
\draw (A) ++ (0.2,-2.8) node[below]{$q$};
%番号

\end{tikzpicture}
\caption{The Kirby diagram obtained from the graph $\Gamma$}
\label{fig:kirby1}
\end{figure}

%\vspace{3cm}

\begin{figure}[htbp]
\centering
\begin{tikzpicture}
\coordinate (A) at (0,0);
\coordinate (B) at (-2,0);
\coordinate (C) at (6,0);
\coordinate (D) at (8,0);

\draw (A) --++ (6,0);
\draw (A)++(0,-1) --++ (6,0);
%手前側

\draw (B) ++ (0,-1) --++ (-0.3,0) to[out=150,in=150] ++(-1,-2) --++ (5,0);
\draw [preaction={draw=white,line width=6pt} ](B)  --++ (-0.3,0) to[out=150,in=150] ++(-1,-2)--++ (4.5,0);
%上下つなぐやつ

%\draw (A) ++ (1.7,0.5) --++ (-1.5,-3);
\draw [preaction={draw=white,line width=6pt} ](A) ++ (2,0.1) to[out=60,in=30] ++(-0.2,0.2) --++ (-1.15,-2.3) to[out=-130, in =-120] ++ (0.25,-0.1);
\draw [preaction={draw=white,line width=6pt} ](A) ++ (1.95,-0.05) --++ (-0.95,-1.9);

%\draw [preaction={draw=white,line width=6pt} ](A) ++ (2.5,-0.9) to[out=60,in=30] ++(-0.2,0.2) --++ (-1.15,-2.3) to[out=-130, in =-120] ++ (0.25,-0.1);
%\draw [preaction={draw=white,line width=6pt} ](A) ++ (2.45,-1.05) --++ (-0.95,-1.9);

%2-ハンドル

%\draw (C)++(0,0.3) --++ (0,-3.4) --++ (2,0) --++(0,3.4)-- cycle;
%右の箱B_{p,q}

%\draw (C) ++ (1,-1.2) node{$B^{-1}_{p,q}$};
%B_{p,q}名前。回転させる。Bの上が右を向くように時計回りに90回転させる。

\draw [preaction={draw=white,line width=2pt} ](B)  ++(0,-1)++ (-0.3,0) ++(-1,-2) ++ (5,0) --++(6,0)to[out=30,in=-30] ++(1,2) --++ (-0.7,0);
\draw [preaction={draw=white,line width=6pt} ](B)  ++ (-0.3,0) ++(-1,-2) ++ (4.5,0) --++(6.5,0) to[out=30,in=-30] ++(1,2) --++ (-0.7,0);
%下の線続き+右側

\fill (A) ++ (0.5,0) circle (0.06);
\fill (A) ++ (0.5,-1) circle (0.06);
%点付き円の点

\draw (A) ++ (0,-1) to[out=180,in=0] ++(-1,1);
\draw (A) ++ (-1,-1) to[out=180,in=0] ++(-1,1);
\draw [preaction={draw=white,line width=6pt} ] (A) to[out=180,in=0] ++(-1,-1);
\draw [preaction={draw=white,line width=6pt} ] (A) ++ (-1,0) to[out=180,in=0] ++(-1,-1);
%捻り左側

\draw(D) to[out=180,in=0] ++(-1,-1);
\draw (D) ++ (-1,0) to[out=180,in=0] ++(-1,-1);
\draw  [preaction={draw=white,line width=6pt} ] (D) ++ (0,-1) to[out=180,in=0] ++(-1,1);
\draw[preaction={draw=white,line width=6pt} ]  (D) ++ (-1,-1) to[out=180,in=0] ++(-1,1);

\end{tikzpicture}
\caption{The case $p=q=2$}
\label{fig:kirby1'}
\end{figure}
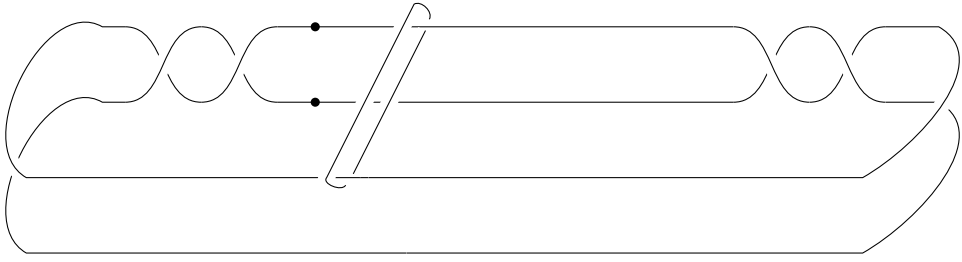

\end{example}

\begin{example} (An arrangement of three generic lines)
Let $f(x,y) = (x+y) (x-y) (y-1)$. The curve $C =\{f(x,y) = 0\}$ is the union of three generic lines. 
Then, $X = \{-1,0,1\}$. Each point in $X$ corresponds to the intersection of two lines. 
Take a graph $\Gamma$ as in Figure \ref{fig:gamma2}.
Each local monodromy around the point in $X$ is isomorphic to the one at the origin of $\{x^2 - y^2 = 0\}$. 
Since $\Gamma$ runs parallel to the real axis close enough except around the points in $X$, non-trivial braidings do not appear over this part.
When $\Gamma$ makes a half turn around a point in $X$, a half twist appears in the Kirby diagram as the half of the monodromy.
Therefore, the Kirby diagram is as in Figure \ref{fig:kirby2}.
Here, $B$ is the braid representing the inverse of all the products of local monodromy.
%Note that the trivialization along the s_{i} is "trivial".

Also, we can take a graph $\Gamma'$ as in Figure \ref{fig:gamma2'}. 
From this graph, we can obtain the handle decomposition as follows, although in a slightly different way from the one described in the main result.
First, attach a $2$-handle corresponding to the fiber around $-1$. Next, after trivializing along the path connecting $-1$ and $0$, attach a $2$-handle corresponding to $0 \in X$. 
At this time, dotted circles look different from the ones obtained from the graph $\Gamma$, that is, the local monodromy around $0$ splits into two half rotations. 
Similarly, we can attach the remaining $2$-handle corresponding to $1 \in X$. 
Finally, the resulting Kirby diagram is described as in Figure \ref{fig:kirby2'}.

\begin{remark}
In \cite{sug-yos}, Kirby diagrams for complements of complexified real line arrangements are obtained. 
We can see that the Kirby diagram described in the example above is obtained from the one in \cite{sug-yos} by handle sliding.
\end{remark}

\begin{figure}[htbp]
\centering
\begin{tikzpicture}
\draw (0,0) circle (2);
\draw (1.6,-1.6) node{$D$};
%円D
\fill (1.6,0) circle (0.06);
\filldraw (-1,0) circle (0.06) node[below]{$-1$};
\filldraw (0,0) circle (0.06) node[below]{$0$};
\filldraw (1,0) circle (0.06) node[below]{$1$};
%p_iたち

\draw (1,0) -- (2,0);
\draw (1.6,0) -- (1.2,0.15) to [out=120,in=90] (0.7,0) -- (0,0);
\draw (-1,0) -- (-0.3,0) to[out=90,in=120] (0.2,0.15)--(0.6,0.15)to[out=60,in=150]  (1.2,0.3) -- (1.6,0);
 
%\draw (-1,0) -- (2,0);
%グラフ\Gamma

\draw (0.8,0.35) to[out=90,in=-30] ++ (-0.3,0.7) node[left]{$\Gamma$};
%\Gamma名前

\end{tikzpicture}
\caption{A graph $\Gamma$ for $f(x,y) = (x+y)(x-y)(y-1)$}
\label{fig:gamma2}
\end{figure}
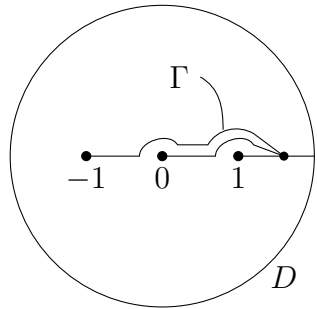

%\vspace{3cm}

\begin{figure}[htbp]
\centering
\begin{tikzpicture}
\coordinate (A) at (0,0);
\coordinate (B) at (-2,0);
\coordinate (B1) at (-1,-2);
\coordinate (C) at (6,0);
\coordinate (D) at (1.5,0);
\coordinate (E) at (1.5,-2);
\coordinate (F) at  (6.7,0);
\coordinate (G) at  (4.7,0);
\coordinate (H) at (6.7,-2);  
\coordinate (I) at  (11.4,0);
\coordinate (J) at  (9.4,0);
\coordinate (K) at (11.4,-2);  
%手前側

\draw (B) ++ (0,-1.8)  to[out=150,in=150] ++(-1,-2) --++ (2,0);
\draw [preaction={draw=white,line width=6pt} ](B) ++ (0,-1) to[out=150,in=150] ++(-1,-2) --++ (2,0);
\draw [preaction={draw=white,line width=6pt} ](B)  to[out=150,in=150] ++(-1,-2)--++ (2,0);
%最左の上下つなぐやつ

\draw (A) --++(-1,0);
\draw (A) ++ (0,-1) --++(-1,0);

\draw (A) ++ (-1,-1) to[out=180,in=0] ++(-0.5,1);
\draw (A) ++ (-1.5,-1) to[out=180,in=0] ++(-0.5,1);
\draw [preaction={draw=white,line width=4pt} ] (A) ++ (-1,0) to[out=180,in=0] ++(-0.5,-1);
\draw [preaction={draw=white,line width=4pt} ] (A) ++ (-1.5,0) to[out=180,in=0] ++(-0.5,-1);
\draw (A) ++ (0,-1.8) --++ (-2,0);
%捻り最左側

%\draw (A) ++ (1.7,0.5) --++ (-1.5,-3);
\draw [preaction={draw=white,line width=6pt} ](A) ++(-2,0) ++ (2,0.1) to[out=60,in=30] ++(-0.2,0.2) --++ (-1.15,-2.3) to[out=-130, in =-120] ++ (0.25,-0.1);
\draw [preaction={draw=white,line width=6pt} ](A)++(-2,0)  ++ (1.95,-0.05) --++ (-0.95,-1.9);
%\draw [preaction={draw=white,line width=6pt} ](A) ++(-0.5,0) ++ (2.5,-0.9) to[out=60,in=30] ++(-0.2,0.2) --++ (-1.15,-2.3) to[out=-130, in =-120] ++ (0.25,-0.1);
%\draw [preaction={draw=white,line width=6pt} ](A) ++(-0.5,0) ++ (2.45,-1.05) --++ (-0.95,-1.9);
%最左の2-ハンドル

%\draw [preaction={draw=white,line width=6pt} ](B)  ++(0,-1) ++(-1,-2) ++ (4,0) to[out=0,in=180] ++(1,-0.8) --++ (1,0);
%\draw [preaction={draw=white,line width=6pt} ](B)  ++(0,-1.8) ++(-1,-2) ++ (4,0) to[out=0,in=180] ++(1,0.8) --++(1,0)  ;
%\draw [preaction={draw=white,line width=6pt} ](B)   ++(-1,-2) ++ (4,0) --++ (2,0) ; 
%to[out=30,in=-30] ++(1,2) --++ (-0.7,0);

\draw (B1) --++ (2.5,0);
\draw (B1) ++ (0,-1.8) to[out=0,in=180] ++ (0.5,0.8) --++(2,0);
\draw [preaction={draw=white,line width=4pt} ](B1) ++ (0,-1)to[out=0,in=180] ++(0.5,-0.8)--++(2,0) ;
%捻り左の右側（下）

%\draw (D) ++ (0,-1) to[out=180,in=0] ++(-1,1);
%\draw (D) ++ (-1,-1) to[out=180,in=0] ++(-1,1);
%\draw [preaction={draw=white,line width=6pt} ] (D) to[out=180,in=0] ++(-1,-1);
%\draw [preaction={draw=white,line width=6pt} ] (D) ++ (-1,0) to[out=180,in=0] ++(-1,-1);

\draw (A) --++ (0.5,0) --++ (0.5,0) to[out=0,in=180] ++(0.5,-1);
\draw (A) ++(0,-1)--++ (0.5,0)  to[out=0,in=180] ++(0.5,-0.8) --++(0.5,0) ;
\draw[preaction={draw=white,line width=4pt} ] (A)++(0,-1.8) --++ (0.5,0) to[out=0,in=180] ++(0.5,0.8) to[out=0,in=180] ++(0.5,1);
%捻り左の右側（上）

\draw(D)  to[out=0,in=180] ++(1.2,2.5)--++(9,0) to[out=0,in=120] ++(0.6,-0.5);
\draw(D) ++(0,-1) to[out=0,in=180] ++(1.2,2.5)--++(9,0)to[out=0,in=120] ++(0.6,-0.5);
\draw(D) ++ (0,-1.8)  to[out=0,in=180] ++(1.2,2.5)--++(9,0)to[out=0,in=120] ++(0.6,-0.5);
%左→中央→右→上下つなぎの上側（上）
\draw(E) --++(0.5,0) to[out=0,in=180] ++(1.2,2.5) --++ (3.5,0) to[out=0,in=90] ++ (0.5,-0.25);
\draw(E)++(0,-1) --++(0.5,0)  to[out=0,in=180] ++(1.2,2.5)  --++ (3.5,0) to[out=0,in=90] ++ (0.5,-0.25);
\draw(E) ++ (0,-1.8)--++(0.5,0)   to[out=0,in=180] ++(1.2,2.5)  --++ (3.5,0) to[out=0,in=90] ++ (0.5,-0.25);
%左→中央（下）

\draw (F)  to[out=0,in=-90] ++(0.5,0.25);
\draw (F)++(0,-1)  to[out=0,in=-90] ++(0.5,0.25);
\draw (F)++(0,-1.8)  to[out=0,in=-90] ++(0.5,0.25);
\draw (F) --++(-2,0);
\draw (F) [preaction={draw=white,line width=6pt} ] ++(0,-1.8) --++(-1,0)to[out=180,in=0] ++(-0.5,0.8);
\draw (F) [preaction={draw=white,line width=6pt} ]++(0,-1)--++(-1,0)to[out=180,in=0] ++(-0.5,-0.8);
\draw (F) [preaction={draw=white,line width=6pt} ]++(0,-1)++(-1.5,-0.8)to[out=180,in=0] ++(-0.5,0.8) ;
\draw (F) [preaction={draw=white,line width=6pt} ]++(0,-1.8)++(-1.5,0.8)to[out=180,in=0] ++(-0.5,-0.8)  ;
%中央上

\draw (G) ++ (0,-1.8)  to[out=150,in=150] ++(-1,-2) --++ (3,0);
\draw [preaction={draw=white,line width=6pt} ](G) ++ (0,-1) to[out=150,in=150] ++(-1,-2) --++ (3,0);
 \draw [preaction={draw=white,line width=6pt} ](G)  to[out=150,in=150] ++(-1,-2)--++ (3,0);
%中央の上下つなぐやつ〜下

\draw [preaction={draw=white,line width=4pt} ](G) ++(0,-1) ++ (2,0.1) to[out=60,in=30] ++(-0.2,0.2) --++ (-1.15,-2.3) to[out=-130, in =-120] ++ (0.25,-0.1);
\draw [preaction={draw=white,line width=4pt} ](G)++(0,-1)  ++ (1.95,-0.05) --++ (-0.95,-1.9);
%中央2ハンドル

\draw [preaction={draw=white,line width=6pt} ](G)  to[out=150,in=150] ++(-1,-2)--++ (3,0);
%中央の上下つなぐやつ〜下

\draw(H)++(0,-1)to[out=0,in=180] ++ (0.5,1)  to[out=0,in=180] ++(1.2,2.5)  --++ (3,0) to[out=0,in=90] ++ (0.5,-0.25);
\draw[preaction={draw=white,line width=4pt} ](H) to[out=0,in=180]++(0.5,-1) to[out=0,in=180] ++(1.2,2.5) --++ (3,0) to[out=0,in=90] ++ (0.5,-0.25);
\draw(H) ++ (0,-1.8)--++(0.5,0)   to[out=0,in=180] ++(1.2,2.5)  --++ (3,0) to[out=0,in=90] ++ (0.5,-0.25);
%中央→右（下）

\draw (I)  to[out=0,in=-90] ++(0.5,0.25);
\draw (I)++(0,-1)  to[out=0,in=-90] ++(0.5,0.25);
\draw (I)++(0,-1.8)  to[out=0,in=-90] ++(0.5,0.25);
\draw (I)++(0,-1.8) --++(-2,0);
\draw (I) [preaction={draw=white,line width=6pt} ] ++(0,-1) --++(-1,0)to[out=180,in=0] ++(-0.5,1);
\draw (I) [preaction={draw=white,line width=6pt} ]++(0,0)--++(-1,0)to[out=180,in=0] ++(-0.5,-1);
\draw (I) [preaction={draw=white,line width=6pt} ]++(0,0)++(-1.5,-1)to[out=180,in=0] ++(-0.5,1) ;
\draw (I) [preaction={draw=white,line width=6pt} ]++(0,-1)++(-1.5,1)to[out=180,in=0] ++(-0.5,-1)  ;
%右上

\draw (J) ++ (0,-1.8)  to[out=150,in=150] ++(-1,-2) --++ (3,0);
\draw [preaction={draw=white,line width=6pt} ](J) ++ (0,-1) to[out=150,in=150] ++(-1,-2) --++ (3,0);
\draw [preaction={draw=white,line width=6pt} ](J)  to[out=150,in=150] ++(-1,-2)--++ (3,0);
%右の上下つなぐやつ〜下

\draw [preaction={draw=white,line width=4pt} ](J) ++(0.3,0) ++ (2,0.1) to[out=60,in=30] ++(-0.2,0.2) --++ (-1.15,-2.3) to[out=-130, in =-120] ++ (0.25,-0.1);
\draw [preaction={draw=white,line width=4pt} ](J)++(0.3,0)  ++ (1.95,-0.05) --++ (-0.95,-1.9);
%右2ハンドル

%\draw(K) --++(0.5,0) to[out=0,in=180] ++(1.2,2.5) --++ (2.5,0) to[out=0,in=90] ++ (0.5,-0.25);
%\draw(K)++(0,-1) --++(0.5,0)  to[out=0,in=180] ++(1.2,2.5)  --++ (2.5,0) to[out=0,in=90] ++ (0.5,-0.25);
%\draw(K) ++ (0,-1.8)--++(0.5,0)   to[out=0,in=180] ++(1.2,2.5)  --++ (2.5,0) to[out=0,in=90] ++ (0.5,-0.25);
%中央→右（下）

\draw (K) ++(0,-1.8) --++(0.5,0) to[out=0,in=-60] ++(1.5,0.5) --++(-0.7,2.1);
\draw [preaction={draw=white,line width=4pt} ](K) ++(0,-1)--++(0.5,0) to[out=0,in=-60] ++(1.5,0.5)--++(-0.7,2.1);
\draw [preaction={draw=white,line width=4pt} ](K) --++(0.5,0) to[out=0,in=-60] ++(1.5,0.5)--++(-0.7,2.1);
%右側上下のつなぎ

\filldraw[fill=white, thick] (I) ++(1.4,0.7) --++(0,-2.6) --++(-0.7,2) --++ (0,2.6) --cycle;
%右の箱

\draw (I) ++ (1,1) to[out=60,in=180] ++(1,1) node[right]{$B$};

%\fill[red] (A) circle (0.06);
%\fill[red] (B) circle (0.06);
%\fill[red] (B1) circle (0.06);
%\fill[red] (C) circle (0.06);
%\fill[red] (D) circle (0.06);
%\fill[red] (E) circle (0.06);
%\fill[red] (F) circle (0.06);
%\fill[red] (G) circle (0.06);
%\fill[red] (H) circle (0.06);
%\fill[red] (I) circle (0.06);
%\fill[red] (J) circle (0.06);
%\fill[red] (K) circle (0.06);
%目印たち

\fill (A) ++ (0.25,0) circle (0.06);
\fill (A) ++ (0.25,-1) circle (0.06);
\fill (A) ++ (0.25,-1.8) circle (0.06);
%点付き円の点

\end{tikzpicture}
\caption{The Kirby diagram obtained from the graph $\Gamma$}
\label{fig:kirby2}
\end{figure}

%\vspace{3cm}

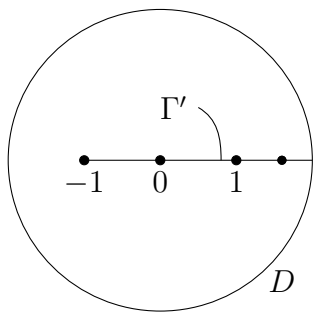
\begin{figure}[htbp]
\centering
\begin{tikzpicture}

\draw (0,0) circle (2);
\draw (1.6,-1.6) node{$D$};
%円D
\fill (1.6,0) circle (0.06);
\filldraw (-1,0) circle (0.06) node[below]{$-1$};
\filldraw (0,0) circle (0.06) node[below]{$0$};
\filldraw (1,0) circle (0.06) node[below]{$1$};
%p_iたち

\fill (1.6,0) circle (0.06);
\draw (-1,0) -- (2,0);
%グラフ\Gamma

\draw (0.8,0) to[out=90,in=-30] ++ (-0.3,0.7) node[left]{$\Gamma'$};
%\Gamma名前
\end{tikzpicture}
\caption{Another graph $\Gamma'$ for $f(x,y) = (x+y)(x-y)(y-1)$}
\label{fig:gamma2'}
\end{figure}

%\vspace{3cm}

\begin{figure}[htbp]
\centering
\begin{tikzpicture}

\coordinate (A) at (0,0);
\coordinate (B) at (-2,0);
\coordinate (B1) at (-1,-2);
\coordinate (C) at (6,0);
\coordinate (D) at (1.5,0);
\coordinate (E) at (1.5,-2);
\coordinate (F) at  (6.5,0);
\coordinate (G) at  (4.5,0);
\coordinate (H) at (6.5,-2);  
\coordinate (I) at  (11,0);
\coordinate (J) at  (9,0);
\coordinate (K) at (11,-2);  
%手前側

\draw (B) ++ (0,-1.8)  to[out=150,in=150] ++(-1,-2) --++ (2,0);
\draw [preaction={draw=white,line width=6pt} ](B) ++ (0,-1) to[out=150,in=150] ++(-1,-2) --++ (2,0);
\draw [preaction={draw=white,line width=6pt} ](B)  to[out=150,in=150] ++(-1,-2)--++ (2,0);
%最左の上下つなぐやつ

\draw (A) --++(-1,0);
\draw (A) ++ (0,-1) --++(-1,0);

\draw (A) ++ (-1,-1) to[out=180,in=0] ++(-0.5,1);
\draw (A) ++ (-1.5,-1) to[out=180,in=0] ++(-0.5,1);
\draw [preaction={draw=white,line width=4pt} ] (A) ++ (-1,0) to[out=180,in=0] ++(-0.5,-1);
\draw [preaction={draw=white,line width=4pt} ] (A) ++ (-1.5,0) to[out=180,in=0] ++(-0.5,-1);
\draw (A) ++ (0,-1.8) --++ (-2,0);
%捻り最左側

%\draw (A) ++ (1.7,0.5) --++ (-1.5,-3);
\draw [preaction={draw=white,line width=6pt} ](A) ++(-2,0) ++ (2,0.1) to[out=60,in=30] ++(-0.2,0.2) --++ (-1.15,-2.3) to[out=-130, in =-120] ++ (0.25,-0.1);
\draw [preaction={draw=white,line width=6pt} ](A)++(-2,0)  ++ (1.95,-0.05) --++ (-0.95,-1.9);
%\draw [preaction={draw=white,line width=6pt} ](A) ++(-0.5,0) ++ (2.5,-0.9) to[out=60,in=30] ++(-0.2,0.2) --++ (-1.15,-2.3) to[out=-130, in =-120] ++ (0.25,-0.1);
%\draw [preaction={draw=white,line width=6pt} ](A) ++(-0.5,0) ++ (2.45,-1.05) --++ (-0.95,-1.9);
%最左の2-ハンドル

%\draw [preaction={draw=white,line width=6pt} ](B)  ++(0,-1) ++(-1,-2) ++ (4,0) to[out=0,in=180] ++(1,-0.8) --++ (1,0);
%\draw [preaction={draw=white,line width=6pt} ](B)  ++(0,-1.8) ++(-1,-2) ++ (4,0) to[out=0,in=180] ++(1,0.8) --++(1,0)  ;
%\draw [preaction={draw=white,line width=6pt} ](B)   ++(-1,-2) ++ (4,0) --++ (2,0) ; 
%to[out=30,in=-30] ++(1,2) --++ (-0.7,0);

\draw (B1) --++ (2.5,0);
\draw (B1) ++ (0,-1)--++(2.5,0) ;
\draw (B1) ++ (0,-1.8)  --++(2.5,0);
%捻り左の右側（下）

%\draw (D) ++ (0,-1) to[out=180,in=0] ++(-1,1);
%\draw (D) ++ (-1,-1) to[out=180,in=0] ++(-1,1);
%\draw [preaction={draw=white,line width=6pt} ] (D) to[out=180,in=0] ++(-1,-1);
%\draw [preaction={draw=white,line width=6pt} ] (D) ++ (-1,0) to[out=180,in=0] ++(-1,-1);

\draw (A) --++ (1.5,0) ;
\draw (A) ++(0,-1)--++ (1.5,0) ;
\draw(A)++(0,-1.8) --++ (1.5,0);
%捻り左の右側（上）

\draw (D) --++(1.5,0)--++(0.5,0)--++(3,0);

\draw (D) ++(0,-1)--++(1.5,0)to[out=0,in=180]++(0.5,-0.8)--++(3,0);
\draw [preaction={draw=white,line width=6pt} ] (D) ++(0,-1.8)--++(1.5,0)to[out=0,in=180]++(0.5,0.8)--++(3,0);
\draw (E) --++(1.5,0)--++(0.5,0)--++(3,0);%←あとで上塗りする
\draw (E) ++(0,-1)--++(1.5,0)to[out=0,in=180]++(0.5,-0.8)--++(3,0);
\draw [preaction={draw=white,line width=6pt} ] (E) ++(0,-1.8)--++(1.5,0)to[out=0,in=180]++(0.5,0.8)--++(3,0);
%中央分

\draw (F) ++(0,-1.8)--++(1.5,0)--++(0.5,0)--++(3,0);
\draw (F) ++(0,0)--++(1.5,0)to[out=0,in=180]++(0.5,-1)--++(3,0);
\draw [preaction={draw=white,line width=6pt} ] (F) ++(0,-1)--++(1.5,0)to[out=0,in=180]++(0.5,1)--++(3,0);
\draw (H) ++(0,-1.8) --++(1.5,0)--++(0.5,0)--++(3,0);
\draw (H) ++(0,0)--++(1.5,0)to[out=0,in=180]++(0.5,-1)--++(3,0);
\draw [preaction={draw=white,line width=6pt} ] (H) ++(0,-1)--++(1.5,0)to[out=0,in=180]++(0.5,1)--++(3,0);
%中央分

\draw [preaction={draw=white,line width=4pt} ](G) ++(-1,-1) ++ (2,0.1) to[out=60,in=30] ++(-0.2,0.2) --++ (-1.15,-2.3) to[out=-130, in =-120] ++ (0.25,-0.1);
\draw [preaction={draw=white,line width=4pt} ](G)++(-1,-1)  ++ (1.95,-0.05) --++ (-0.95,-1.9);
%中央2ハンドル

\draw [preaction={draw=white,line width=4pt}](E) --++(1.5,0)--++(0.5,0)--++(3,0);%←上塗り
%中央塗り直し

\draw [preaction={draw=white,line width=4pt} ](J) ++(-1,0) ++ (2,0.1) to[out=60,in=30] ++(-0.2,0.2) --++ (-1.15,-2.3) to[out=-130, in =-120] ++ (0.25,-0.1);
\draw [preaction={draw=white,line width=4pt} ](J)++(-1,0)  ++ (1.95,-0.05) --++ (-0.95,-1.9);
%右2ハンドル

\draw (K) ++(0,-1.8) ++(0.5,0) to[out=0,in=-60] ++(1.5,0.5) --++(-0.7,2.1);
\draw [preaction={draw=white,line width=4pt} ](K) ++(0,-1)++(0.5,0) to[out=0,in=-60] ++(1.5,0.5)--++(-0.7,2.1);
\draw [preaction={draw=white,line width=4pt} ](K) ++(0.5,0) to[out=0,in=-60] ++(1.5,0.5)--++(-0.3,0.9);
%右側上下のつなぎ（下側）

\draw (I) --++(0.5,0) to[out=0,in=150] ++(0.8,-0.2);
\draw (I) ++(0,-1)--++(0.5,0) to[out=0,in=150] ++(0.8,-0.2);
\draw (I) ++(0,-1.8)--++(0.5,0) to[out=0,in=150] ++(0.8,-0.2);
%右側上下のつなぎ（上側）

\filldraw[fill=white,thick] (I) ++(1.8,-0.6) --++(0,-2.6) --++(-0.7,1.2) --++ (0,2.6) --cycle;
%右の箱

\draw (K) ++(0,-1.8) --++(0.5,0)to[out=0,in=-60] ++(1.5,0.5);
\draw [preaction={draw=white,line width=4pt} ](K) ++(0,-1)--++(0.5,0)to[out=0,in=-60] ++(1.5,0.5);
\draw [preaction={draw=white,line width=4pt} ](K) --++(0.5,0)to[out=0,in=-60] ++(1.5,0.5);
%右側上下のつなぎ（下側）２

\draw (I) ++ (1.5,0) to[out=60,in=180] ++(0.5,1) node[right]{$B$};
%K名前

%\fill[red] (A) circle (0.06);
%\fill[red] (B) circle (0.06);
%\fill[red] (B1) circle (0.06);
%\fill[red] (C) circle (0.06);
%\fill[red] (D) circle (0.06);
%\fill[red] (E) circle (0.06);
%\fill[red] (G) circle (0.06);
%\fill[red] (H) circle (0.06);
%\fill[red] (I) circle (0.06);
%\fill[red] (J) circle (0.06);
%\fill[red] (K) circle (0.06);
%目印たち

\fill (A) ++ (0.25,0) circle (0.06);
\fill (A) ++ (0.25,-1) circle (0.06);
\fill (A) ++ (0.25,-1.8) circle (0.06);
%点付き円の点

\end{tikzpicture}
\caption{The Kirby diagram obtained from the graph $\Gamma'$}
\label{fig:kirby2'}
\end{figure}

\end{example}

\begin{example} (Fermat curves)
Let $f (x,y) = x^n + y^n -1$. Then, $X = \{1, \zeta, \zeta^{2}, \dots, \zeta^{n-1}\}$, where $\zeta^{k} = e^{2 k\pi  \sqrt{-1}/n}$. 
%Since the curve $C =\{f(x,y) = 0\}$ is smooth, every fiber $\pi^{-1} (\zeta^{k}) $ tangents to $C$. 
%Remark that these tangent points are not locally isomorphic to $\{y=x^2\}$, however, we can describe explicit monodromy as follows.
Take the graph $\Gamma$ as in Figure \ref{fig:gamma3}. 
The complement $D \setminus \nu (\Gamma)$ of the neighborhood of the graph can be deformed into $D'$ defined as follows.
First, define subsets $U_0, s_{0}^{\pm}, \gamma_{0}$ in $\bC$ as follows:
\begin{eqnarray*}
U_{0} &=& \{x \in \bC \mid |x-1| \leq \varepsilon, \arg(x) \notin (-\delta, \delta)\}, \\
s_{0}^{\pm} &=& \{x = u+ \sqrt{-1} v \mid 1 - \varepsilon \cos \delta \leq u \leq 2, v=\pm \sin \delta\}, \\
\gamma_{0} &=&\{x \in \bC \mid |x| =2, \arg (x) \in [\delta , (2 \pi /n) - \delta]\}. \\
\end{eqnarray*}
Let $U_{k}, s_{k}^{\pm}, \gamma_{k}$ ($k=1, \dots, n-1$) be the sets obtained from $U_{0}, s_{0}^{\pm}, \gamma_{0}$ by rotating $2k\pi/n$ degrees. 
At this time, define $D'$ as the domain containing the origin rounded by the curve $\bigcup_{k=0}^{n-1} (\partial U_{k} \cup s^{+}_{k} \cup s^{-}_{k} \cup \gamma_{k})$ (see Figure \ref{fig:gamma3'}).
Each $2$-handle is attached at the fiber of the area $2k\pi/n - \delta \leq \theta \leq 2k\pi/n + \delta$ around $\zeta^{k} \in X$.
By symmetry, we will describe the handle attachment corresponding to just $1 \in X$. 
For each other $\zeta^{k} \in X$, we can similarly obtain the handle attachment.

Since the defining equation is $y^{n} = x^{n} -1$, $L_{1} \cap C$ is the set consisting of one point, to which $n$ points of the intersection between the fiber of a point in the neighborhood of $1 \in X$ and $C$ degenerate. %英文？
Around the tangent point $(1,0)$ in $L_1 \cap C$, $C$ is locally isomorphic to $\{x=y^n\}$.
From easy computation, when $x$ turn around $\partial U_{0}$, the fiber coordinate $y$ rotates $2\pi/n$.
Also, when $x$ moves along the path $\gamma_{0}$, the fiber coordinate $y$ also rotates $2\pi/n$.
Therefore, we can obtain the picture of $\pi^{-1} (\partial D') \cap C$ and the Kirby diagram of the complement. (See Figure \ref{fig:kirby3} for $n=3$ case).
For general $n$, we can describe a similar diagram, and moving by an isotopy, we obtain the Kirby diagram as in Figure \ref{fig:kirby3'}.  
\end{example}

\begin{remark}
All non-singular curves that intersect the line at infinity transversally for a fixed degree have diffeomorphic complements.
Therefore, the Kirby diagram obtained in the above example gives a Kirby diagram of the complement of a non-singular curve intersecting the line at infinity transversally.
\end{remark}

\begin{figure}[htbp]
\centering
\begin{tikzpicture}

\draw (0,0) circle (2);
\draw (1.6,-1.6) node{$D$};
%円D
\fill (1.6,0) circle (0.06);
\filldraw (1,0) circle (0.06) node[below]{$1$};
\filldraw (0.8,0.6) circle (0.06) node[below]{$\zeta$};
\filldraw (0.2,0.9) circle (0.06) node[below]{$\zeta^2$};
\filldraw (0.8,-0.6) circle (0.06) node[left]{$\zeta^{n-1}$};
%p_iたち

\fill (-1,0) circle(0.03);
\fill (-0.9,0.4) circle(0.03);
\fill (-0.9,-0.4) circle(0.03);
%てんてん

\draw (1.6,0) to[out=130,in=0] (0.8,0.6);
\draw (1.6,0) to[out=100,in=30] (0.2,0.9);
\draw (1.6,0) to[out=60,in=0]  (0,1.8)to[out=180,in=90] (-1.8,0)to[out=-90,in=180](0,-1.8)to[out=0,in=-65] (0.8,-0.6);
\draw (1,0) -- (2,0);
%グラフ\Gamma

\draw (1.5,1) to[out=90,in=-130] ++ (0.3,0.7) node[right]{$\Gamma$};
%\Gamma名前
\end{tikzpicture}
\caption{A graph $\Gamma$ for $f(x,y) =y^n - x^n-1$}
\label{fig:gamma3}
\end{figure}

%\vspace{3cm}

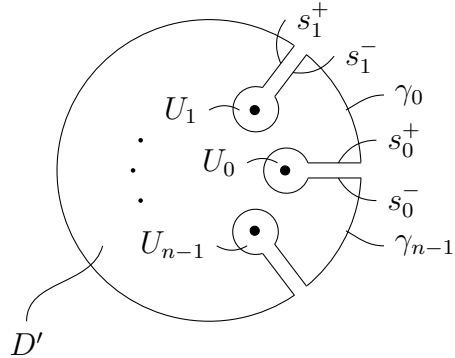
\begin{figure}[htbp]
\centering
\begin{tikzpicture}

\draw (0,0) circle (2);
\draw (-1.4,-1) to[out=-105,in=90]++(-1,-1) node[below]{$D'$};
%円D
%\fill (1.6,0) circle (0.06);
\filldraw (1,0) circle (0.06); %node[below]{$1$}
\filldraw (0.6,0.8) circle (0.06); %node[below]{$\zeta$};
%\filldraw (0.2,0.9) circle (0.06) node[below]{$\zeta^2$};
\filldraw (0.6,-0.8) circle (0.06); %node[left]{$\zeta^{n-1}$};
%p_iたち

\draw [white,very thick] (1.99,-0.1) arc (-4:2:2);
\draw (1.3,0.1) --++(0.7,0);
\draw (1.3,-0.1) --++(0.7,0);
\draw (1.3,0.1) arc (20:340:0.3);
%1まわり

\draw [white,very thick,rotate=53] (1.99,-0.1) arc (-4:2:2);
\draw[rotate=53] (1.3,0.1) --++(0.7,0);
\draw[rotate=53] (1.3,-0.1) --++(0.7,0);
\draw[rotate=53] (1.3,0.1) arc (20:340:0.3);
%zetaまわり

\draw [white,very thick,rotate=-53] (1.99,-0.1) arc (-4:2:2);
\draw[rotate=-53] (1.3,0.1) --++(0.7,0);
\draw[rotate=-53] (1.3,-0.1) --++(0.7,0);
\draw[rotate=-53] (1.3,0.1) arc (20:340:0.3);
%zeta^n-1まわり

\draw (1.7,0.1) to[out=60,in=180] ++(0.5,0.3) node[right]{$s_{0}^{+}$};
\draw  (1.7,-0.1) to[out=-60,in=180] ++(0.5,-0.3) node[right]{$s_{0}^{-}$};
\draw [rotate=53](1.7,0.1) to[out=60,in=180] ++(0.5,0.3) node[right]{$s_{1}^{+}$};
\draw [rotate=53](1.7,-0.1) to[out=-60,in=180] ++(0.5,-0.3) node[right]{$s_{1}^{-}$};
\draw (1.8,0.8) to[out=30,in=180] ++(0.5,0.2) node[right]{$\gamma_{0}$};
\draw (1.8,-0.8) to[out=-30,in=180] ++(0.5,-0.2) node[right]{$\gamma_{n-1}$};
\draw (0.9,0.1) to[out=150,in=30] ++ (-0.4,-0) node[left]{$U_{0}$};
\draw (0.4,0.8) to[out=150,in=30] ++ (-0.4,-0) node[left]{$U_{1}$};
\draw (0.5,-1) to[out=-150,in=-30] ++ (-0.4,-0) node[left]{$U_{n-1}$};

\fill (-1,0) circle(0.03);
\fill (-0.9,0.4) circle(0.03);
\fill (-0.9,-0.4) circle(0.03);
%てんてん

%\draw (1.6,0) to[out=130,in=0] (0.8,0.6);
%\draw (1.6,0) to[out=100,in=30] (0.2,0.9);
%\draw (1.6,0) to[out=60,in=0]  (0,1.8)to[out=180,in=90] (-1.8,0)to[out=-90,in=180](0,-1.8)to[out=0,in=-65] (0.8,-0.6);
%\draw (1,0) -- (2,0);
%グラフ\Gamma

%\draw (1.5,1) to[out=90,in=-130] ++ (0.3,0.7) node[right]{$\Gamma$};
%\Gamma名前
\end{tikzpicture}
\caption{The disk $D'$}
\label{fig:gamma3'}
\end{figure}

%\vspace{3cm}

\begin{figure}[htbp]
\centering
\begin{tikzpicture}

\coordinate (A) at (5,1.8);
\coordinate (B) at (3,1.8);
\coordinate (B1) at (4,-0.2);
\coordinate (O) at (0,0);

\draw (B) ++ (0,-1.8)  to[out=150,in=150] ++(-1,-2) --++ (0.5,0);
\draw [preaction={draw=white,line width=6pt} ](B) ++ (0,-1) to[out=150,in=150] ++(-1,-2) --++ (0.5,0);
\draw [preaction={draw=white,line width=6pt} ](B)  to[out=150,in=150] ++(-1,-2)--++ (0.5,0);
%最左の上下つなぐやつ

\draw (B1) --++(-0.5,0) to[out=180,in=0] ++(-0.5,-1)to[out=180,in=0] ++ (-0.5,-0.8);
\draw[preaction={draw=white,line width=4pt} ] (B1) ++(0,-1)--++(-0.5,0) to[out=180,in=0] ++(-0.5,1) --++(-0.5,0) ;
\draw [preaction={draw=white,line width=4pt} ](B1) ++(0,-1.8) --++(-0.5,0)--++(-0.5,0) to[out=180,in=0] ++ (-0.5,0.8);
%下のひねり

\draw (A) --++(-1,0);
\draw (A) ++ (0,-1) --++(-1,0);

\draw (A) ++ (0,-1.8) --++ (-2,0);
\draw (A) ++ (0,-1) --++ (-2,0);
\draw (A) ++ (0,0) --++ (-2,0);
%捻り最左側

%\draw (A) ++ (1.7,0.5) --++ (-1.5,-3);
\draw [preaction={draw=white,line width=6pt} ](A) ++(-1.5,0) ++ (2,0.1) to[out=60,in=30] ++(-0.2,0.2) --++ (-1.15,-2.3) to[out=-130, in =-120] ++ (0.25,-0.1);
\draw [preaction={draw=white,line width=6pt} ](A)++(-1.5,0)  ++ (1.95,-0.05) --++ (-0.95,-1.9);
%\draw [preaction={draw=white,line width=6pt} ](A) ++(-0.5,0) ++ (2.5,-0.9) to[out=60,in=30] ++(-0.2,0.2) --++ (-1.15,-2.3) to[out=-130, in =-120] ++ (0.25,-0.1);
%\draw [preaction={draw=white,line width=6pt} ](A) ++(-0.5,0) ++ (2.45,-1.05) --++ (-0.95,-1.9);
%最左の2-ハンドル

%\draw [preaction={draw=white,line width=6pt} ](B)  ++(0,-1) ++(-1,-2) ++ (4,0) to[out=0,in=180] ++(1,-0.8) --++ (1,0);
%\draw [preaction={draw=white,line width=6pt} ](B)  ++(0,-1.8) ++(-1,-2) ++ (4,0) to[out=0,in=180] ++(1,0.8) --++(1,0)  ;
%\draw [preaction={draw=white,line width=6pt} ](B)   ++(-1,-2) ++ (4,0) --++ (2,0) ; 
%to[out=30,in=-30] ++(1,2) --++ (-0.7,0);

\draw (B1) --++ (2.5,0);
\draw (B1) ++ (0,-1)--++(0.5,0)--++(2,0) ;
\draw [preaction={draw=white,line width=4pt} ](B1) ++ (0,-1.8) --++ (0.5,0) --++(2,0);
%捻り左の右側（下）

%\draw (D) ++ (0,-1) to[out=180,in=0] ++(-1,1);
%\draw (D) ++ (-1,-1) to[out=180,in=0] ++(-1,1);
%\draw [preaction={draw=white,line width=6pt} ] (D) to[out=180,in=0] ++(-1,-1);
%\draw [preaction={draw=white,line width=6pt} ] (D) ++ (-1,0) to[out=180,in=0] ++(-1,-1);

\draw (A) --++ (0.5,0) --++ (1,0);
\draw (A) ++(0,-1)--++ (1.5,0) ;
\draw (A)++(0,-1.8) --++ (0.5,0) --++(1,0);
%捻り左の右側（上）

%\draw (A) ++ (1.7,0.5) --++ (-1.5,-3);
\draw [preaction={draw=white,line width=4pt} ](A) ++(-1.5,0) ++ (2,0.1) to[out=60,in=30] ++(-0.2,0.2) --++ (-1.15,-2.3) to[out=-130, in =-120] ++ (0.25,-0.1);
\draw [preaction={draw=white,line width=4pt} ](A)++(-1.5,0)  ++ (1.95,-0.05) --++ (-0.95,-1.9);
\draw [preaction={draw=white,line width=4pt} ](A) ++(-1.8,0) ++ (2.5,-0.9) to[out=60,in=30] ++(-0.2,0.2) --++ (-1.15,-2.3) to[out=-130, in =-120] ++ (0.25,-0.1);
\draw [preaction={draw=white,line width=4pt} ](A) ++(-1.8,0) ++ (2.45,-1.05) --++ (-0.95,-1.9);
%最左の2-ハンドル

\draw [preaction={draw=white,line width=4pt} ](B1) ++(0.7,0) --++ (1.8,0);
%

%\fill[red] (A) circle (0.06);
%\fill[red] (B) circle (0.06);
%\fill[red] (B1) circle (0.06);
%\fill[red] (C) circle (0.06);
%\fill[red] (D) circle (0.06);
%\fill[red] (E) circle (0.06);
%\fill[red] (F) circle (0.06);
%\fill[red] (G) circle (0.06);
%\fill[red] (H) circle (0.06);
%\fill[red] (I) circle (0.06);
%\fill[red] (J) circle (0.06);
%\fill[red] (K) circle (0.06);
%\fill[red] (O) circle(0.06);
%目印たち

%%%%%%%%%%%%%%%%%%%%ここから120度回転%%%%%%%%%%%%%%%%%%%%
\coordinate (A-1) at (-4.06,3.4301);
\coordinate (B-1) at (-3.05,1.698);
\coordinate (B1-1) at (-1.8,3.564);
\coordinate (O) at (0,0);

\draw[rotate=120] (B-1) ++ (0,-1.8)  to[out=150,in=150] ++(-1,-2) --++ (0.5,0);
\draw [rotate=120] [preaction={draw=white,line width=6pt} ](B-1) ++ (0,-1) to[out=150,in=150] ++(-1,-2) --++ (0.5,0);
\draw [rotate=120][preaction={draw=white,line width=6pt} ](B-1)  to[out=150,in=150] ++(-1,-2)--++ (0.5,0);
%最左の上下つなぐやつ

\draw [rotate=120](B1-1) --++(-0.5,0) to[out=180,in=0] ++(-0.5,-1)to[out=180,in=0] ++ (-0.5,-0.8);
\draw[rotate=120][preaction={draw=white,line width=4pt} ] (B1-1) ++(0,-1)--++(-0.5,0) to[out=180,in=0] ++(-0.5,1) --++(-0.5,0) ;
\draw [rotate=120][preaction={draw=white,line width=4pt} ](B1-1) ++(0,-1.8) --++(-0.5,0)--++(-0.5,0) to[out=180,in=0] ++ (-0.5,0.8);
%下のひねり

\draw[rotate=120] (A-1) --++(-1,0);
\draw [rotate=120](A-1) ++ (0,-1) --++(-1,0);

\draw [rotate=120](A-1) ++ (0,-1.8) --++ (-2,0);
\draw[rotate=120] (A-1) ++ (0,-1) --++ (-2,0);
\draw[rotate=120] (A-1) ++ (0,0) --++ (-2,0);
%捻り最左側

%\draw (A) ++ (1.7,0.5) --++ (-1.5,-3);
\draw[rotate=120] [preaction={draw=white,line width=6pt} ](A-1) ++(-1.5,0) ++ (2,0.1) to[out=60,in=30] ++(-0.2,0.2) --++ (-1.15,-2.3) to[out=-130, in =-120] ++ (0.25,-0.1);
\draw[rotate=120] [preaction={draw=white,line width=6pt} ](A-1)++(-1.5,0)  ++ (1.95,-0.05) --++ (-0.95,-1.9);
%\draw [preaction={draw=white,line width=6pt} ](A-1) ++(-0.5,0) ++ (2.5,-0.9) to[out=60,in=30] ++(-0.2,0.2) --++ (-1.15,-2.3) to[out=-130, in =-120] ++ (0.25,-0.1);
%\draw [preaction={draw=white,line width=6pt} ](A-1) ++(-0.5,0) ++ (2.45,-1.05) --++ (-0.95,-1.9);
%最左の2-ハンドル

%\draw [preaction={draw=white,line width=6pt} ](B-1)  ++(0,-1) ++(-1,-2) ++ (4,0) to[out=0,in=180] ++(1,-0.8) --++ (1,0);
%\draw [preaction={draw=white,line width=6pt} ](B-1)  ++(0,-1.8) ++(-1,-2) ++ (4,0) to[out=0,in=180] ++(1,0.8) --++(1,0)  ;
%\draw [preaction={draw=white,line width=6pt} ](B-1)   ++(-1,-2) ++ (4,0) --++ (2,0) ; 
%to[out=30,in=-30] ++(1,2) --++ (-0.7,0);

\draw [rotate=120](B1-1) --++ (2.5,0);
\draw [rotate=120](B1-1) ++ (0,-1)--++(0.5,0)--++(2,0) ;
\draw[rotate=120] [preaction={draw=white,line width=4pt} ](B1-1) ++ (0,-1.8) --++ (0.5,0) --++(2,0);
%捻り左の右側（下）

%\draw (D) ++ (0,-1) to[out=180,in=0] ++(-1,1);
%\draw (D) ++ (-1,-1) to[out=180,in=0] ++(-1,1);
%\draw [preaction={draw=white,line width=6pt} ] (D) to[out=180,in=0] ++(-1,-1);
%\draw [preaction={draw=white,line width=6pt} ] (D) ++ (-1,0) to[out=180,in=0] ++(-1,-1);

\draw [rotate=120](A-1) --++ (0.5,0) --++ (1,0);
\draw [rotate=120](A-1) ++(0,-1)--++ (1.5,0) ;
\draw[rotate=120] (A-1)++(0,-1.8) --++ (0.5,0) --++(1,0);
%捻り左の右側（上）

%\draw (A-1) ++ (1.7,0.5) --++ (-1.5,-3);
\draw [rotate=120][preaction={draw=white,line width=4pt} ](A-1) ++(-1.5,0) ++ (2,0.1) to[out=60,in=30] ++(-0.2,0.2) --++ (-1.15,-2.3) to[out=-130, in =-120] ++ (0.25,-0.1);
\draw [rotate=120][preaction={draw=white,line width=4pt} ](A-1)++(-1.5,0)  ++ (1.95,-0.05) --++ (-0.95,-1.9);
\draw[rotate=120] [preaction={draw=white,line width=4pt} ](A-1) ++(-1.8,0) ++ (2.5,-0.9) to[out=60,in=30] ++(-0.2,0.2) --++ (-1.15,-2.3) to[out=-130, in =-120] ++ (0.25,-0.1);
\draw [rotate=120][preaction={draw=white,line width=4pt} ](A-1) ++(-1.8,0) ++ (2.45,-1.05) --++ (-0.95,-1.9);
%最左の2-ハンドル

\draw [rotate=120][preaction={draw=white,line width=4pt} ](B1-1) ++(0.7,0) --++ (1.8,0);
%

%\fill[red] (A-1) circle (0.06);
%\fill[red] (B-1) circle (0.06);
%\fill[red] (B1-1) circle (0.06);
%\fill[red] (C) circle (0.06);
%\fill[red] (D) circle (0.06);
%\fill[red] (E) circle (0.06);
%\fill[red] (F) circle (0.06);
%\fill[red] (G) circle (0.06);
%\fill[red] (H) circle (0.06);
%\fill[red] (I) circle (0.06);
%\fill[red] (J) circle (0.06);
%\fill[red] (K) circle (0.06);
%\fill[red] (O) circle(0.06);
%目印たち

%%%%%%%%%%%%%%%%%%%ここから240度回転%%%%%%%%%%%%%%%%%%%%
\coordinate (A-2) at (-0.941,-5.230);
\coordinate (B-2) at (0.058,-3.498);
\coordinate (B1-2) at (-2.1732,-3.3641);
\coordinate (O) at (0,0);

\draw[rotate=240] (B-2) ++ (0,-1.8)  to[out=150,in=150] ++(-1,-2) --++ (0.5,0);
\draw [rotate=240] [preaction={draw=white,line width=6pt} ](B-2) ++ (0,-1) to[out=150,in=150] ++(-1,-2) --++ (0.5,0);
\draw [rotate=240][preaction={draw=white,line width=6pt} ](B-2)  to[out=150,in=150] ++(-1,-2)--++ (0.5,0);
%最左の上下つなぐやつ

\draw [rotate=240](B1-2) --++(-0.5,0) to[out=180,in=0] ++(-0.5,-1)to[out=180,in=0] ++ (-0.5,-0.8);
\draw[rotate=240][preaction={draw=white,line width=4pt} ] (B1-2) ++(0,-1)--++(-0.5,0) to[out=180,in=0] ++(-0.5,1) --++(-0.5,0) ;
\draw [rotate=240][preaction={draw=white,line width=4pt} ](B1-2) ++(0,-1.8) --++(-0.5,0)--++(-0.5,0) to[out=180,in=0] ++ (-0.5,0.8);
%下のひねり

\draw[rotate=240] (A-2) --++(-1,0);
\draw [rotate=240](A-2) ++ (0,-1) --++(-1,0);

\draw [rotate=240](A-2) ++ (0,-1.8) --++ (-2,0);
\draw[rotate=240] (A-2) ++ (0,-1) --++ (-2,0);
\draw[rotate=240] (A-2) ++ (0,0) --++ (-2,0);
%捻り最左側

%\draw (A) ++ (1.7,0.5) --++ (-1.5,-3);
\draw[rotate=240] [preaction={draw=white,line width=6pt} ](A-2) ++(-1.5,0) ++ (2,0.1) to[out=60,in=30] ++(-0.2,0.2) --++ (-1.15,-2.3) to[out=-130, in =-120] ++ (0.25,-0.1);
\draw[rotate=240] [preaction={draw=white,line width=6pt} ](A-2)++(-1.5,0)  ++ (1.95,-0.05) --++ (-0.95,-1.9);
%\draw [preaction={draw=white,line width=6pt} ](A-2) ++(-0.5,0) ++ (2.5,-0.9) to[out=60,in=30] ++(-0.2,0.2) --++ (-1.15,-2.3) to[out=-130, in =-120] ++ (0.25,-0.1);
%\draw [preaction={draw=white,line width=6pt} ](A-2) ++(-0.5,0) ++ (2.45,-1.05) --++ (-0.95,-1.9);
%最左の2-ハンドル

%\draw [preaction={draw=white,line width=6pt} ](B-2)  ++(0,-1) ++(-1,-2) ++ (4,0) to[out=0,in=180] ++(1,-0.8) --++ (1,0);
%\draw [preaction={draw=white,line width=6pt} ](B-2)  ++(0,-1.8) ++(-1,-2) ++ (4,0) to[out=0,in=180] ++(1,0.8) --++(1,0)  ;
%\draw [preaction={draw=white,line width=6pt} ](B-2)   ++(-1,-2) ++ (4,0) --++ (2,0) ; 
%to[out=30,in=-30] ++(1,2) --++ (-0.7,0);

\draw [rotate=240](B1-2) --++ (2.5,0);
\draw [rotate=240](B1-2) ++ (0,-1)--++(0.5,0)--++(2,0) ;
\draw[rotate=240] [preaction={draw=white,line width=4pt} ](B1-2) ++ (0,-1.8) --++ (0.5,0) --++(2,0);
%捻り左の右側（下）

%\draw (D) ++ (0,-1) to[out=180,in=0] ++(-1,1);
%\draw (D) ++ (-1,-1) to[out=180,in=0] ++(-1,1);
%\draw [preaction={draw=white,line width=6pt} ] (D) to[out=180,in=0] ++(-1,-1);
%\draw [preaction={draw=white,line width=6pt} ] (D) ++ (-1,0) to[out=180,in=0] ++(-1,-1);

\draw [rotate=240](A-2) --++ (0.5,0) --++ (1,0);
\draw [rotate=240](A-2) ++(0,-1)--++ (1.5,0) ;
\draw[rotate=240] (A-2)++(0,-1.8) --++ (0.5,0) --++(1,0);
%捻り左の右側（上）

%\draw (A-2) ++ (1.7,0.5) --++ (-1.5,-3);
\draw [rotate=240][preaction={draw=white,line width=4pt} ](A-2) ++(-1.5,0) ++ (2,0.1) to[out=60,in=30] ++(-0.2,0.2) --++ (-1.15,-2.3) to[out=-130, in =-120] ++ (0.25,-0.1);
\draw [rotate=240][preaction={draw=white,line width=4pt} ](A-2)++(-1.5,0)  ++ (1.95,-0.05) --++ (-0.95,-1.9);
\draw[rotate=240] [preaction={draw=white,line width=4pt} ](A-2) ++(-1.8,0) ++ (2.5,-0.9) to[out=60,in=30] ++(-0.2,0.2) --++ (-1.15,-2.3) to[out=-130, in =-120] ++ (0.25,-0.1);
\draw [rotate=240][preaction={draw=white,line width=4pt} ](A-2) ++(-1.8,0) ++ (2.45,-1.05) --++ (-0.95,-1.9);
%最左の2-ハンドル

\draw [rotate=240][preaction={draw=white,line width=4pt} ](B1-2) ++(0.7,0) --++ (1.8,0);
%

%\fill[red] (A-2) circle (0.06);
%\fill[red] (B-2) circle (0.06);
%\fill[red] (B1-2) circle (0.06);
%\fill[red] (C) circle (0.06);
%\fill[red] (D) circle (0.06);
%\fill[red] (E) circle (0.06);
%\fill[red] (F) circle (0.06);
%\fill[red] (G) circle (0.06);
%\fill[red] (H) circle (0.06);
%\fill[red] (I) circle (0.06);
%\fill[red] (J) circle (0.06);
%\fill[red] (K) circle (0.06);
%\fill[red] (O) circle(0.06);
%目印たち

%%%少し1-ハンドルを消す%%%
\draw [white,line width=6pt] (A) ++ (0,-2) ++ (0,-1) ++(0.78,0) --++(1,0);
\draw [white,line width=6pt] (A) ++ (0,-2) ++ (0,-1.8) ++(-0.12,0) --++(1.8,0);

\draw [rotate=120][white,line width=6pt] (A-1) ++ (0,-2) ++ (0,-1) ++(0.78,0) --++(1,0);
\draw [rotate=120][white,line width=6pt] (A-1) ++ (0,-2) ++ (0,-1.8) ++(-0.12,0) --++(1.8,0);

\draw [rotate=240][white,line width=6pt] (A-2) ++ (0,-2) ++ (0,-1) ++(0.78,0) --++(1,0);
\draw [rotate=240][white,line width=6pt] (A-2) ++ (0,-2) ++ (0,-1.8) ++(-0.12,0) --++(1.8,0);

%%%%%%%%%%%%%%%%%%%%%%%%ここからつなぐやつ

%%%0度回転％％％
\draw  (A) ++(0,-1.8)++ (0.5,0)--++(1,0) to[out=0,in=-30] ++(-2,3);
\draw  [preaction={draw=white,line width=6pt} ](A) ++(0,-1)++ (0.5,0)--++(1,0) to[out=0,in=-30] ++(-2,3);
\draw [preaction={draw=white,line width=6pt} ](A) ++ (0.5,0)--++(1,0) to[out=0,in=-30] ++(-2,3);
%前半
\draw (A) ++ (0,-1.8) ++ (1.5,0) ++ (-2,3) to[out=150,in=-30] ++(-0.6,1.3) to[out=150,in=-30] ++(-0.6,1.3);
\draw [preaction={draw=white,line width=6pt} ](A) ++ (0,-1) ++ (1.5,0) ++ (-2,3) to[out=150,in=-30] ++(-0.6,-0.4) --++(-0.6,0.5);
\draw [preaction={draw=white,line width=6pt} ](A) ++ (0,0) ++ (1.5,0) ++ (-2,3) --++(-0.6,0.4) to[out=150,in=-30] ++(-0.6,-0.5);
%捻り
\draw (A) ++ (1.5,0) ++ (-2,3) ++ (-0.6,0.4) ++ (-0.6,-0.5) to [out=150,in=70] ++(-5.13,0.9);
\draw (A) ++ (0,-1) ++(1.5,0) ++ (-2,3) ++ (-0.6,-0.4) ++ (-0.6,0.5) to[out=150,in=70] ++(-4,1.3);
\draw (A) ++ (0,-1.8) ++ (1.5,0) ++ (-2,3) ++(-0.6,1.3) ++(-0.6,1.3) to[out=150,in=70] ++(-6.34,0.13);
%後半

%%%120度回転％％％
\draw  [rotate=120](A-1) ++(0,-1.8)++ (0.5,0)--++(1,0) to[out=0,in=-30] ++(-2,3);
\draw[rotate=120]  [preaction={draw=white,line width=6pt} ](A-1) ++(0,-1)++ (0.5,0)--++(1,0) to[out=0,in=-30] ++(-2,3);
\draw[rotate=120] [preaction={draw=white,line width=6pt} ](A-1) ++ (0.5,0)--++(1,0) to[out=0,in=-30] ++(-2,3);
%前半
\draw [rotate=120](A-1) ++ (0,-1.8) ++ (1.5,0) ++ (-2,3) to[out=150,in=-30] ++(-0.6,1.3) to[out=150,in=-30] ++(-0.6,1.3);
\draw [rotate=120][preaction={draw=white,line width=6pt} ](A-1) ++ (0,-1) ++ (1.5,0) ++ (-2,3) to[out=150,in=-30] ++(-0.6,-0.4) --++(-0.6,0.5);
\draw[rotate=120] [preaction={draw=white,line width=6pt} ](A-1) ++ (0,0) ++ (1.5,0) ++ (-2,3) --++(-0.6,0.4) to[out=150,in=-30] ++(-0.6,-0.5);
%捻り
\draw[rotate=120] (A-1) ++ (1.5,0) ++ (-2,3) ++ (-0.6,0.4) ++ (-0.6,-0.5) to [out=150,in=70] ++(-5.13,0.9);
\draw[rotate=120] (A-1) ++ (0,-1) ++(1.5,0) ++ (-2,3) ++ (-0.6,-0.4) ++ (-0.6,0.5) to[out=150,in=70] ++(-4,1.3);
\draw [rotate=120](A-1) ++ (0,-1.8) ++ (1.5,0) ++ (-2,3) ++(-0.6,1.3) ++(-0.6,1.3) to[out=150,in=70] ++(-6.37,0.13);
%後半
%%%240度回転％％％
\draw  [rotate=240](A-2) ++(0,-1.8)++ (0.5,0)--++(1,0) to[out=0,in=-30] ++(-2,3);
\draw[rotate=240]  [preaction={draw=white,line width=6pt} ](A-2) ++(0,-1)++ (0.5,0)--++(1,0) to[out=0,in=-30] ++(-2,3);
\draw[rotate=240] [preaction={draw=white,line width=6pt} ](A-2) ++ (0.5,0)--++(1,0) to[out=0,in=-30] ++(-2,3);
%前半
\draw [rotate=240](A-2) ++ (0,-1.8) ++ (1.5,0) ++ (-2,3) to[out=150,in=-30] ++(-0.6,1.3) to[out=150,in=-30] ++(-0.6,1.3);
\draw [rotate=240][preaction={draw=white,line width=6pt} ](A-2) ++ (0,-1) ++ (1.5,0) ++ (-2,3) to[out=150,in=-30] ++(-0.6,-0.4) --++(-0.6,0.5);
\draw[rotate=240] [preaction={draw=white,line width=6pt} ](A-2) ++ (0,0) ++ (1.5,0) ++ (-2,3) --++(-0.6,0.4) to[out=150,in=-30] ++(-0.6,-0.5);
%捻り
\draw[rotate=240] (A-2) ++ (1.5,0) ++ (-2,3) ++ (-0.6,0.4) ++ (-0.6,-0.5) to [out=150,in=70] ++(-5.13,0.9);
\draw[rotate=240] (A-2) ++ (0,-1) ++(1.5,0) ++ (-2,3) ++ (-0.6,-0.4) ++ (-0.6,0.5) to[out=150,in=70] ++(-4,1.3);
\draw [rotate=240](A-2) ++ (0,-1.8) ++ (1.5,0) ++ (-2,3) ++(-0.6,1.3) ++(-0.6,1.3) to[out=150,in=70] ++(-6.37,0.13);
%後半

%%%%%%ドット%%%%
\fill (B1) ++ (-0.3,0) circle (0.06);
\fill (B1) ++ (-0.3,-1) circle (0.06);
\fill (B1) ++ (-0.3,-1.8) circle (0.06);

\end{tikzpicture}
\caption{The Kirby diagram obtained from the graph $\Gamma'$ ($n=3$)}
\label{fig:kirby3}
\end{figure}

%\vspace{3cm}

\begin{figure}[htbp]
\centering
\begin{tikzpicture}
\coordinate(A) at (0,0);
\coordinate (A1) at (-0.5,-1.5);
\coordinate (A2) at (-0.5,-3);
\coordinate (A3) at (-0.5,-5.5);

\coordinate(B) at (2.5,0);
\coordinate (B1) at (2,-2.1);
\coordinate (B2) at (2,-3.6);
\coordinate (B3) at (2,-6.1);

\coordinate(B') at (5,0);

\coordinate(C) at (8.5,0);
\coordinate (C1) at (8,-1.5);
\coordinate (C2) at (8,-3);
\coordinate (C3) at (8,-5.5);

\coordinate(D) at (11,0);
\coordinate (A1) at (-0.5,-1.5);

\draw (A) to[out=180,in=90] ++ (-0.5,-1)--++(0,-5)to[out=-90,in=180] ++ (0.5,-1);
\draw (A) to[out=0,in=90] ++ (0.5,-1) --++(0,-5)to[out=-90,in=0]++(-0.5,-1);
\draw (B) to[out=180,in=90] ++ (-0.5,-1)--++(0,-5)to[out=-90,in=180] ++ (0.5,-1);
\draw (B) to[out=0,in=90] ++ (0.5,-1) --++(0,-5)to[out=-90,in=0]++(-0.5,-1);
\draw (B') to[out=180,in=90] ++ (-0.5,-1)--++(0,-5)to[out=-90,in=180] ++ (0.5,-1);
\draw (B') to[out=0,in=90] ++ (0.5,-1) --++(0,-5)to[out=-90,in=0]++(-0.5,-1);

\draw (C) to[out=180,in=90] ++ (-0.5,-1)--++(0,-5)to[out=-90,in=180] ++ (0.5,-1);
\draw (C) to[out=0,in=90] ++ (0.5,-1) --++(0,-5)to[out=-90,in=0]++(-0.5,-1);
\draw (D) to[out=180,in=90] ++ (-0.5,-1)--++(0,-5)to[out=-90,in=180] ++ (0.5,-1);
\draw (D) to[out=0,in=90] ++ (0.5,-1) --++(0,-5)to[out=-90,in=0]++(-0.5,-1);
%点つき円

\draw [preaction={draw=white,line width=6pt} ](A1) ++ (-0.2,0) to[out=180,in=180] ++(0,0.4) --++(2.5,0);
\draw [preaction={draw=white,line width=6pt} ](A1) ++(0.2,0) --++(2.5,0) to[out=0,in=0] ++(0,0.4);

\draw [preaction={draw=white,line width=6pt} ](A2) ++ (-0.2,0) to[out=180,in=180] ++(0,0.4) --++(2.5,0);
\draw [preaction={draw=white,line width=6pt} ](A2) ++(0.2,0) --++(2.5,0) to[out=0,in=0] ++(0,0.4);

\draw [preaction={draw=white,line width=6pt} ](A3) ++ (-0.2,0) to[out=180,in=180] ++(0,0.4) --++(2.5,0);
\draw [preaction={draw=white,line width=6pt} ](A3) ++(0.2,0) --++(2.5,0) to[out=0,in=0] ++(0,0.4);

\draw [preaction={draw=white,line width=6pt} ](B1) ++ (-0.2,0) to[out=180,in=180] ++(0,0.4) --++(2.5,0);
\draw [preaction={draw=white,line width=6pt} ](B1) ++(0.2,0) --++(2.5,0) to[out=0,in=0] ++(0,0.4);

\draw [preaction={draw=white,line width=6pt} ](B2) ++ (-0.2,0) to[out=180,in=180] ++(0,0.4) --++(2.5,0);
\draw [preaction={draw=white,line width=6pt} ](B2) ++(0.2,0) --++(2.5,0) to[out=0,in=0] ++(0,0.4);

\draw [preaction={draw=white,line width=6pt} ](B3) ++ (-0.2,0) to[out=180,in=180] ++(0,0.4) --++(2.5,0);
\draw [preaction={draw=white,line width=6pt} ](B3) ++(0.2,0) --++(2.5,0) to[out=0,in=0] ++(0,0.4);

\draw [preaction={draw=white,line width=6pt} ](C1) ++ (-0.2,0) to[out=180,in=180] ++(0,0.4) --++(2.5,0);
\draw [preaction={draw=white,line width=6pt} ](C1) ++(0.2,0) --++(2.5,0) to[out=0,in=0] ++(0,0.4);

\draw [preaction={draw=white,line width=6pt} ](C2) ++ (-0.2,0) to[out=180,in=180] ++(0,0.4) --++(2.5,0);
\draw [preaction={draw=white,line width=6pt} ](C2) ++(0.2,0) --++(2.5,0) to[out=0,in=0] ++(0,0.4);

\draw [preaction={draw=white,line width=6pt} ](C3) ++ (-0.2,0) to[out=180,in=180] ++(0,0.4) --++(2.5,0);
\draw [preaction={draw=white,line width=6pt} ](C3) ++(0.2,0) --++(2.5,0) to[out=0,in=0] ++(0,0.4);

\fill (A) ++ (1.1,0) ++ (0,-4.4)circle(0.04);
\fill (A) ++ (1.1,0) ++ (0,-3.5)circle(0.04);
\fill (A) ++ (1.1,0) ++ (0,-3.9)circle(0.04);
\fill (B) ++ (0,-3.5) ++ (4.5,0)circle (0.04);
\fill (B) ++ (0,-3.5) ++ (3.5,0)circle (0.04);
\fill (B) ++ (0,-3.5) ++ (4,0)circle (0.04);
%てんてん

\fill(A) circle (0.06);
\fill(B) circle (0.06);
\fill(B') circle (0.06);
\fill(C) circle (0.06);
\fill(D) circle (0.06);
%点つき円

\draw (A) ++ (-1,-1) to[out=-150,in=30] ++(-0.7,-2.5);
\draw (A) ++ (-1,-6) to[out=150,in=-30] ++ (-0.7,2.5)node[left]{$n$};

\draw (A) ++ (0,-7.5) to[out=-30,in=120] ++(5.5,-0.5);
\draw (A) ++ (11,-7.5) to[out=-150,in=60] ++(-5.5,-0.5)node[below]{$n$};
\end{tikzpicture}
\caption{The Kirby diagram for general $n$}
\label{fig:kirby3'}
\end{figure}
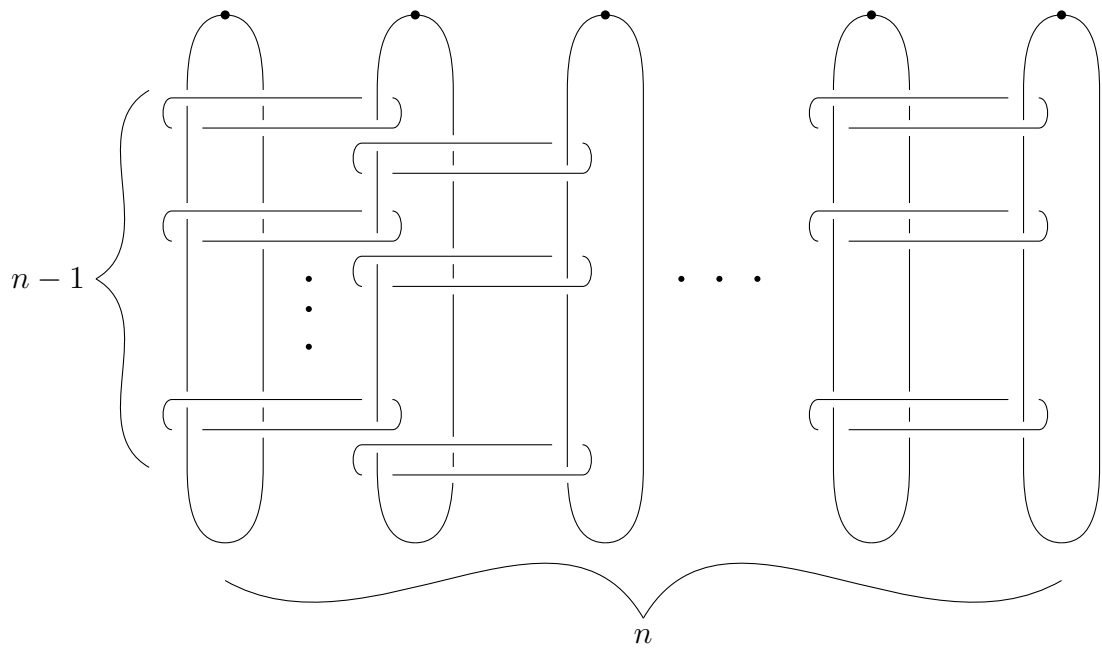

%\begin{figure}[htbp]
%\centering
%\includegraphics[width=5cm]{bitmap.png}

%\end{figure}

\end{document}